\documentclass[11pt]{article}

\usepackage{amssymb}
\usepackage{amsthm}
\usepackage{amsmath}
\usepackage{amscd}

\allowdisplaybreaks
    
\newcommand{\R}{\mathbb{R}}

\newcommand{\Sph}{\mathbb{S}}
\newcommand{\T}{\mathbb{T}}
\newcommand{\eps}{\varepsilon}

\newcommand{\dif}{\mathrm{d}} 
\newcommand{\Dif}{\mathrm{D}}

\newcommand{\linspan}{\mathrm{span}}

\newcommand{\cobg}{C^{0,\beta}_{\delta - 2}}
\newcommand{\ctbg}{C^{2,\beta}_\delta}
\newcommand{\clbg}{C^{k,\beta}_\delta}

\newcommand{\btheta}{\boldsymbol{\Theta}}
\newcommand{\cliff}{\T^{n_1,n_2}_ \alpha}

\newcommand{\approxsol}{\tilde \Lambda_{\alpha, \tau}}
\newcommand{\difop}{\Phi_{\alpha, \tau}}
\newcommand{\linop}{\tilde{\mathcal L}_{\alpha, \tau}}
\newcommand{\secapprox}{\tilde \Lambda_{\alpha, \tau, \sigma}}
\newcommand{\secdifop}{\Phi_{a, \tau, \sigma }}
\newcommand{\seclinop}{\tilde{\mathcal L}_{\alpha, \tau, \sigma}}
\newcommand{\handleapprox}{\tilde {\Lambda }_{\alpha, \tau, \sigma}}

\newcommand{\doubleapprox}{\tilde {\Lambda }_{\alpha, \tau, \sigma}}

\DeclareMathOperator{\arccosh}{arccosh}

\DeclareMathOperator{\arccot}{arccot}

\newtheorem{mainthm}{Main Theorem}
\newtheorem{thm}{Theorem}
\newtheorem{lemma}[thm]{Lemma}

\newtheorem{prop}[thm]{Proposition}
\newtheorem*{nonumthm}{Theorem}

\theoremstyle{definition}
\newtheorem{defn}[thm]{Definition}

\newtheoremstyle{rmk}{5pt}{5pt}{}{}{\scshape}{:}{.5em}{}
\theoremstyle{rmk}
\newtheorem*{rmk}{Remark}

\newcommand{\mylabel}
	{\label}
\begin{document}

\title{Gluing Constructions Amongst \\ Constant Mean Curvature Hypersurfaces in $\Sph^{n+1}$}

\author{
\begin{minipage}{3.125in}
	\begin{center}
		Adrian Butscher \\ University of Toronto at Scarborough \\ email: \ttfamily butscher@utsc.utoronto.ca
	\end{center}
\end{minipage}
}

\maketitle

\begin{abstract}
	Four constructions of constant mean curvature (CMC) hypersurfaces in $\Sph^{n+1}$ are given, which should be considered analogues of `classical' constructions that are possible for CMC hypersurfaces in Euclidean space. First, Delaunay-like hypersurfaces, consisting roughly of a chain of hyperspheres winding multiple times around an equator, are shown to exist for all values of the mean curvature.  Second, a hypersurface is constructed which consists of two chains of spheres winding around a pair of orthogonal equators, showing that Delaunay-like hypersurfaces can be fused together in a symmetric manner.  Third, a Delaunay-like handle can be attached to a generalized Clifford torus of the same mean curvature.  Finally, two generalized Clifford tori of equal but opposite mean curvature of any magnitude can be attached to each other by symmetrically positioned Delaunay-like `arms'.  This last result extends Butscher and Pacard's doubling construction for generalized Clifford tori of small mean curvature.
\end{abstract}

\tableofcontents

\renewcommand{\baselinestretch}{1.25}
\normalsize

\section{Introduction}

\subsection{Gluing Constructions of Constant Mean Curvature Hypersurfaces}

A constant mean curvature (CMC) hypersurface $\Lambda$ contained in an ambient Riemannian manifold $X$ of dimension $n+1$ has the property that its mean curvature with respect to the induced metric is constant.  This property ensures that $n$-dimensional area of $\Lambda$ is a critical value of the area functional for hypersurfaces of $X$ subject to an enclosed-volume constraint.  Constant mean curvature hypersurfaces have been objects of great interest since the beginnings of modern differential geometry.  Classical examples of non-trivial CMC surfaces in three-dimensional Euclidean space $\R^3$ are the sphere, the cylinder and the Delaunay surfaces, and for a long while these were the only known CMC surfaces.  In fact, a result of Alexandrov \cite{alexandrov} states that the only compact, connected, embedded CMC surfaces in $\R^3$ are spheres.

The theory of CMC surfaces in $\R^3$ took a leap forward with Wente's discovery of a family of compact, though immersed, CMC tori in 1986 \cite{wente}, and has progressed considerably since then.   For instance, the techniques used by Wente have culminated in a representation for CMC surfaces in $\R^3$ akin to the classically-known Weierstra\ss\ representation of minimal surfaces in which a harmonic but not anti-conformal map from a Riemann surface to the unit sphere becomes the Gau\ss\ map of a CMC immersion into $\R^3$ from which the immersion can be determined \cite{dpw,kenmotsu}.  Most important for this paper, however, has been the development of the \emph{gluing technique} for constructing CMC hypersurfaces in $\R^3$ from simple building blocks. This technique was pioneered by Kapouleas and uses perturbation arguments from the theory of geometric partial differential equations to construct many new CMC surfaces: e.g.~compact, immersed,  genus 2 surfaces by fusing Wente tori together \cite{kapouleas5}; as well as compact surfaces of higher genus and non-compact surfaces with arbitrary numbers of ends by fusing together spheres and Delaunay surfaces \cite{kapouleas2,kapouleas3}. Kapouleas' discoveries have since been complemented by much research into gluing and other constructions from geometric analysis that can be performed in the class of CMC surfaces, most notably in the work of Mazzeo, Pacard and others \cite{jlelipacard,mazzeopacard1,mazzeopacardpollack,mazzeo}, where such constructions as the attachment of Delaunay ends to CMC surfaces and other types of connected sums are realized.

The corresponding picture amongst CMC hypersurfaces of higher dimension or in other ambient manifolds is not nearly as rich as it is in $\R^3$.  There is a certain amount of literature on CMC hypersurfaces in hyperbolic space \cite{bryanthyper,pacardpimentel,saerptoubiana,umeharayamada}; but due to the non-compactness of hyperbolic space, this theory can be considered not such a vast departure from the theory of CMC hypersurfaces in $\R^{n+1}$.  Much less is known when the ambient space is the sphere.  The classically known examples in $\Sph^{n+1}$ are the hyperspheres obtained from intersecting $\Sph ^{n+1}$ with affine hyperplanes, and the so-called generalized Clifford tori which are products of lower-dimensional spheres of the form $\T^{p,q}_\alpha : = \Sph^p \bigl( \cos (\alpha) \bigr) \times \Sph^q \bigl( \sin (\alpha) \bigr)$  for $p+q = n$ and $\alpha \in (0 , \pi /2 )$.  These are embedded hypersurfaces in $\Sph^{n+1}$ with constant mean curvature equal to $H_\alpha : =  q \, \cot (\alpha) - p \tan (\alpha)$.  There are few other examples; and up to now, no truly general methods for the construction of CMC surfaces in $\Sph ^{n+1}$ have been developed.  However, the gluing technique is well-adapted to generalization to other settings, for philosophical as well as practical reasons: gluing is meant to assemble complex objects from simple ones, and this principle should hold equally well in $\Sph^{n+1}$ as in $\R^{n+1}$ because one has simple building blocks in $\Sph^{n+1}$ and a large group of isometries with which to move them around; and many of the operations involved in a gluing construction --- such as forming connected sums using small bridging surfaces near a point of mutual tangency --- are all \emph{local} and thus have straightforward generalizations to $\Sph^{n+1}$.   

This paper adapts the gluing technique for constructing CMC hypersurfaces in $\Sph^{n+1}$ and uses it to realize various analogues of the `classical' constructions that are possible in Euclidean space.  The central idea is to position hyperspheres and/or generalized Clifford tori of the same mean curvature throughout $\Sph^{n+1}$ in various ways such that each building blocks is separated from its neighbours by a small amount.  After appropriate small modifications of this initial configuration, catenoidal necks are then inserted between the building blocks at the points where they come closest to each other.  This \emph{approximate solution} of the gluing problem now has approximately constant mean curvature, where the error between its mean curvature and the mean curvature of the original building blocks is small, except in the neck regions where the approximate solution is almost singular.   One then attempts to perturb the approximate solution by a small normal deformation until it has exactly constant mean curvature.

As is often the case in gluing constructions, there are obstructions to the solvability of the constant mean curvature equation that come from the \emph{Jacobi fields} of the approximate solution.  These are eigenfunctions of the linearized mean curvature operator with zero or small eigenvalue tending to zero as the approximate solution becomes singular.  As is well known, their origin is geometric: every one-parameter family of isometries of the approximate solution generates a Jacobi field with zero eigenvalue; and every one parameter family of isometries, modified by a cut-off function so that it only moves one of the building blocks of the approximate solution while fixing all the others, generates a space of functions which approximates the small eigenspaces.

The key to overcoming these obstructions and realizing the CMC gluing construction in a particular case is to combine two techniques: \emph{symmetry conditions} and \emph{force balancing arguments}.  Both of these techniques play an important role in the gluing constructions in Euclidean space and continue to do so here, though it turns out that the symmetries of the initial configurations play a more prominent role here than they do in Euclidean space.  

The role of the symmetry conditions is as follows.    If the approximate solution is invariant under a group of isometries of $\Sph^{n+1}$, then one can attempt to perturb the approximate solution in an equivariant manner.  Any Jacobi fields that are \emph{not} invariant under the isometry group are thus automatically excluded from consideration.  If it turns out that the isometry group is sufficiently large so that \emph{all} Jacobi fields are excluded, then the obstructions to the solvability of the constant mean curvature equation disappear and the perturbation to an exactly CMC hypersurface can be realized.

However, if it turns out that the approximate solution possesses Jacobi fields that are invariant under all of its isometries, then the method above can not succeed in eliminating all the obstructions to the solvability of the constant mean curvature equations.   This is where the so-called force balancing arguments come into play.  The way to proceed is to derive the \emph{balancing condition} that must be satisfied by the approximate solution if a perturbation to an exactly CMC hypersurface is to exist.   The idea is that the obstructions provided by the invariant Jacobi fields can be avoided exactly when the approximate solution satisfies this balancing condition. 

The balancing condition is best explained in the more general context found in Korevaar-Kusner-Solomon's work \cite{kks}.   First, suppose that $\Lambda$ is a hypersurface with constant mean curvature $h$ in an $(n+1)$-dimensional Riemannian manifold $(X, g)$ possessing a Killing field $V$.  Let $\mathcal U$ be an open set in $\Lambda$ and $\bar{\mathcal U}$ be an open set in $X$ such that $\partial \bar{ \mathcal U} = \partial \mathcal U \cup Q$ where $Q$ is a bounded $n$-dimensional hypersurface.  Then the first variation formula for the $n$-volume of $\mathcal U$ subject to the constraint of constant enclosed $(n+1)$-volume of $\bar{\mathcal U}$ in the direction of the variation determined by $V$ implies 
\begin{equation}
	\mylabel{eqn:balancing}
	\int_{\partial \mathcal U} g(\nu, V) - h \int_Q g(N, V) = 0
\end{equation}
where $\nu$ is the unit normal vector field of $\partial \mathcal U$ in $\Lambda$ and $N$ is the unit normal vector field of $Q$ in $X$.  One can now apply this formula to the approximate solution of the CMC perturbation problem, having mean curvature approximately equal to $h$, in the following way.  Choose the open set $\mathcal U$ as one of the building blocks of the approximate solution.  Then $\partial \mathcal U$ consists of a disjoint union of small $(n-1)$-spheres at the centres of the necks attaching $\mathcal U$ to its neighbours, and $Q$ is the disjoint union of the small disks that cap these spheres off.  The left hand side of \eqref{eqn:balancing} now encodes information about the width and location of the neck regions of $\mathcal U$.  If the left hand side of \eqref{eqn:balancing} is sufficiently close to zero, then one says that $\mathcal U$ is in \emph{balanced position}.  The idea is now that in order to be able to overcome the obstructions to the solvability of the constant mean curvature equations,  the approximate solution must be constructed in such a way that all its building blocks are in balanced position.

It is important to realize that the balancing condition amounts to a form of local symmetry satisfied by each building block with respect to its nearest neighbours.  This is similar to what happens in Euclidean space.  However, because compact hypersurfaces of the sphere must close up on themselves, it turns out that certain global symmetries are a consequence of the totality of local symmetries imposed by the balancing condition.  It is for this reason that symmetry conditions play a more prominent role in the determination of initial configurations of building blocks that can be glued together and perturbed to have constant mean curvature.

Force balancing in itself is not the end of the story; and a balanced approximate solution can not necessarily be perturbed to an exactly CMC hypersurface.  It is in addition necessary to be able to re-position the various building blocks with respect to each other so as to maintain the force balancing condition even if these are deformed by a small perturbation.  Technically speaking, this amounts to the requirement that the mapping taking a re-positioned approximate solution to a set of small real numbers via the integrals on the left hand side of \eqref{eqn:balancing} be invertible.  This requirement is well known and makes an appearance in the CMC gluing problem in Euclidean space.  It also arises here, and is in fact more difficult to achieve than in Euclidean space.  However, the imposition of additional symmetry conditions on the initial configuration of building blocks is enough to overcome this problem.

\subsection{Statement of Results}

The following analogues in $\Sph^{n+1}$ of `classical' constructions amongst CMC hypersurfaces in Euclidean space will be realized in this paper.  First, hypersurfaces will be constructed which consist of a chain of equally spaced hyperspheres winding a number of times around an equator.  These should be construed as analogues of the Delaunay hypersurfaces, and will be called \emph{Delaunay-like}.  Depending on the choices made in the construction, these can be embedded, immersed or non-compact.  Next, one can attempt to fuse a number of Delaunay-like hypersurfaces together.  This fusion will be realized in one particular instance, thereby illustrating the method that could be used for more ambitious situations.  Next, a more general problem along the lines of the fusion of two Delaunay-like hypersurfaces is adding Delaunay-like handles to a given CMC hypersurface.  This construction will be realized again in one particular instance, where the CMC hypersurface in question is a generalized Clifford torus.  Finally, it is possible to \emph{double} two generalized Clifford tori having equal but opposite mean curvature of any magnitude by attaching them to each other with a number of symmetrically placed Delaunay-like `arms'. This is not so much the analogue of a `classical' construction amongst CMC hypersurfaces but rather an extension of the result obtained by Butscher and Pacard in \cite{mepacard1,mepacard2} in which two generalized Clifford tori of equal but opposite mean curvature very close to zero are doubled by attaching them to each other with a number of symmetrically placed small catenoidal necks.  

The spirit of this paper is to give a collection of examples rather than fully general constructions.  The reasons for this are two-fold.  First, it is beneficial for one's intuition to have readily visualizable concrete examples of CMC hypersurfaces of $\Sph^{n+1}$ at hand.  Second, the details of the construction vary greatly depending upon the symmetry group of the initial configuration, the invariant Jacobi fields of the approximate solution, the number of windings needed for a configuration to close up, and so on. Therefore it is best to illustrate the methods used in the gluing construction in a few special cases, which can then be readily adapted to other settings.

\paragraph*{Delaunay-like CMC hypersurfaces of \boldmath $\Sph^{n+1}$.} Let $S_\alpha$ be the  hypersphere obtained by intersecting $\Sph^{n+1}$ by an affine hyperplane passing a distance $\cos(\alpha) \in (0,1)$ from the origin.  The mean curvature of $S_\alpha$ is the constant $H_\alpha:= n \cot(\alpha)$.

\begin{mainthm} 
	\mylabel{result1}
	Suppose $\gamma$ is the great circle in $\Sph^{n+1}$ generated by the one-parameter family of rotations $R_\theta \in SO(n+2)$.  For every $\alpha \in (0,\pi/2)$ and sufficiently small $\tau > 0$ there exists a CMC hypersurface $\Lambda_{\alpha, \tau}$ of mean curvature $H_\alpha := n \cot(\alpha)$ which is approximately equal to a union of  hyperspheres of the form $R_{\theta}^k(S_\alpha)$ that are separated by a distance $\tau$ from each other and connected by small catenoidal necks.  These hypersurfaces are either non-compact and immersed, compact and immersed, or compact and embedded depending on the values of $\alpha$ and $\tau$.  Finally, as $\tau \rightarrow 0$, then $\Lambda_{\alpha, \tau}$ converges in the $C^\infty$ topology to the union of  hyperspheres of mean curvature $H_\alpha$ positioned end-to-end along $\gamma$.
\end{mainthm}





\paragraph*{Gluing two Delaunay-like hypersurfaces together.}  One can imagine a large variety of ways in which two or more Delaunay-like hypersurfaces can be pasted together.  One simple way is to glue a Delaunay-like hypersurface to one whose equator has been rotated.  The resulting configuration consists of two sets of hyperspheres, winding around two different equators.  The only overlap occurs at two antipodal hyperspheres; this is where the first Delaunay-like hypersurface will be glued to the second.  It is worth noting that every such configuration is in balanced position --- this is because balancing requires only that when a hypersphere is to be glued to its neighbours at either two or four points, then these points are opposite each other.  Consequently, the balancing map can not be invertible as a function of the geometry of the initial configuration and the CMC perturbation can not be realized.  This reflects the fact that one does \emph{not} expect to find families of solutions using the gluing technique in a compact setting.  However, if the angle made by the equators is $\pi/2$ then the initial configuration is invariant under additional reflection symmetries.  Now force balancing arguments together with symmetry conditions make it possible to realize the CMC perturbation.

\begin{mainthm} 
	\mylabel{result2}
	Let $\Lambda_{\alpha, \tau}$ denote the Delaunay-like CMC hypersurface winding around the geodesic $\gamma$ constructed in Main Theorem 1.  Let $Q \in SO(n+2)$ be a rotation by an angle of $\pi/2$ about an axis perpendicular to the normal vector that defines $\gamma$.  Then for suitable choices of $\alpha$ and $\tau$, it is possible to form the connected sum of $\Lambda_{\alpha, \tau}$ and $Q(\Lambda_{\alpha, \tau})$.  Moreover, if both of these hypersurfaces are embedded, then so is their connected sum.
\end{mainthm}

\paragraph*{Attaching Delaunay-like handles to Clifford tori.}  Similar methods as in the proof of Main Theorem 2 can be used to add Delaunay-like \emph{handles} to a CMC hypersurface of $\Sph^{n+1}$ under certain circumstances.  That is, suppose $\Lambda$ is a CMC hypersurface in $\Sph^{n+1}$ of constant mean curvature $H$ such that two are two points $p$ and $q$ in $\Lambda$ that can be connected to each other by a geodesic $\gamma$ that ideally does not touch $\Lambda$ in any other points.  Then if the initial configuration consisting of $\Lambda$ and  hyperspheres of parameter $\alpha = \arccot(H/n)$ positioned along $\gamma$ possesses enough symmetry, then one might expect that the perturbation and force-balancing arguments of Main Theorem 2 carry through. 

The symmetry group of $\Lambda$ as well as its Jacobi fields clearly play a crucial role in determining if the construction outlined above is possible.  Instead of stating general conditions on an arbitrary CMC hypersurface $\Lambda$ under which the perturbation to a CMC hypersurface is can be executed, the third main theorem of this paper is to give one concrete example of the construction which can then serve as an archetype for more general constructions.  In this example, $\Lambda$ is the generalized Clifford torus $\cliff := \{ (\cos(\alpha) \Theta_1, \sin(\alpha) \Theta_2 ) : \Theta_j \in \Sph^{n_j} \} \subseteq \R^{n_1+1} \times \R^{n_2+1}$ where $n_1 + n_2 = n$, which has mean curvature $H_\alpha := n_2 \cot(\alpha) - n_1 \tan(\alpha)$.

\begin{mainthm}
	\mylabel{result3}
	Let $\cliff$ be the generalized Clifford torus of $\Sph^{n+1}$.  Then there is a geodesic $\gamma$ that meets $\cliff$ orthogonally in two points $p, p'$.  For suitable choices of $\alpha \in (0, \pi/2)$, small $\tau$ and positive integers $N$ and $m$, there exists a CMC hypersurface of $\Sph^{n+1}$ which is approximately equal to the connected sum at the points $p, p'$ of $\cliff$ and a truncated Delaunay-like hypersurface consisting of $N$ spherical regions separated a distance $\tau$ from each other and winding $m$ times around $\gamma$.  If $\alpha$ is sufficiently close to $\pi/2$ then it is possible to have $m=1$ in this construction, and the resulting hypersurface is embedded.
\end{mainthm}	

\paragraph*{Doubling Clifford tori --- revisited.}

Butscher and Pacard have constructed in \cite{mepacard1,mepacard2} constant mean curvature hypersurfaces of small but non-zero mean curvature by \emph{doubling} the unique minimal Clifford torus in the family of Clifford tori of $\Sph^{n+1}$.  That is, these hypersurfaces are small perturbations of two nearby parallel translates of $\T^{n_1,n_2}_{\alpha_\ast}$, where $\alpha_\ast = \arctan \bigl( \sqrt{n_2/n_1}\, \bigr)$, which are glued together at a number of symmetrically located points $\{ p_1, \ldots , p_J\}$ by means of small catenoidal necks.  The conditions for the existence of these hypersurfaces are phrased in terms of the location of the gluing points: if the gluing points are such that the symmetries of $\T^{n_1, n_2}_{\alpha_\ast} \setminus \{ p_1, \ldots , p_J\}$ rule out all the Jacobi fields of the linearized mean curvature operator, then it is possible to solve the constant mean curvature equations.

The methods developed in this paper allow the above construction to be extended in the following way.  One can take two Clifford tori with equal but opposite mean curvature of \emph{any} magnitude, choose a collection of gluing points on each one according to the same requirements as in \cite{mepacard1,mepacard2}, and place a sequence of hyperspheres of the same mean curvature end-to-end along the geodesics connecting the gluing points on one Clifford torus to those on the other.  The entire configuration can then be glued together by inserting catenoidal necks between the hyperspheres, and then perturbed to have exactly constant mean curvature.  The symmetry condition on the location of the points $\{ p_1, \ldots , p_J\}$ is exactly the same as the one used in \cite{mepacard1,mepacard2}.  That is, if the full group of symmetries preserving $\cliff$ is $O(n_1+1) \times O(n_2+1) \subseteq O(n+2)$ acting diagonally on $\R^{n+2} = \R^{n_1+1} \times \R^{n_2+1}$, then the finite subgroups of interest will be of the form $G \subseteq \{ (\omega^1 , \omega^2) : \omega^j \in O(n_j+1) \}$ and will contain the element $T := (T^1, T^2)$ where $T^j \in O(n_j+1)$ is the reflection symmetry across the $x^1_j$-axis given by $T^j (x_j^1, x_j^2,\ldots, x_j^{n_j+1}) := (x_j^1, -x_j^2, \ldots, -x_j^{n_j+1})$. 

\begin{mainthm}
	Consider two generalized Clifford tori $\cliff$ and $\T_{\bar \alpha}^{n_1,n_2}$ with equal but opposite mean curvature $H_\alpha = - H_{\bar \alpha}$ and a finite subgroup $G \subseteq O(n+2)$ of the type discussed above.  Assume that there are no $G$-invariant functions of the form
	$$(x_1, x_2) \longmapsto \sum_{j=1}^{n_1+1} \sum_{j'=1}^{n_2+1} a_{j j'} x_1^j x_2^{j'}
	$$
	where $a_{jj'} \in \R$.  Let $p$ and $\bar p$ be two points in $\cliff$ and $\T_{\bar \alpha}^{n_1,n_2}$, respectively, that are connected to each other by the normal geodesic $\gamma$.  Let $M = \{ p_1 = p, \ldots, p_{|G|}\}$ and $\bar M= \{ \bar p_1 = \bar p, \ldots, \bar p_{|G|} \}$ be the orbits under $G$ of $p$ and $\bar p$, and let $\gamma_s$ be the geodesic separating $p_s$ from $\bar p_s$.  For suitable choices of $\alpha \in (0, \pi/2)$, small $\tau$ and positive integers $N$ and $m$ there is a $G$-invariant constant mean curvature hypersurface consisting roughly of $\cliff$ and $\T_{\bar \alpha}^{n_1,n_2}$ attached to each other by a chain of $N$ hyperspheres of the same mean curvature winding $m$ times around each of the geodesics $\gamma_1, \ldots, \gamma_{|G|}$.  For sufficiently large values of the mean curvature, it is possible to find embedded hypersurfaces of this kind.
\end{mainthm}

As a final comment before launching into the proof of the first main theorem, it is important to realize what can and can not be achieved in the constructions above.  The gluing technique is by its very nature a singular perturbation method.  Thus the only hypersurfaces that can be created with it are nearly singular, resembling CMC building blocks separated by an extremely small amount $\tau$ and connected by extremely small necks of high principal curvatures.   Two geometric consequences of this fact are worth mentioning; and both revolve around the need to create \emph{compact} glued CMC hypersurfaces.  First, small $\tau$ and fixed mean curvature $H$ can result in extremely large windings of the geodesics used in the constructions.  The number of windings is dictated by the need to for the glued configuration to close and remain compact.  Moreover, a small change in $\tau$ necessarily results in a large change in the number of hyperspheres positioned along the geodesic for the same reason.  Second, the only way small changes in $\tau$ can produced equally small changes in the glued configuration is if the mean curvature is allowed to change.  This is because changing $\tau$ creates changes the amount of space between the building blocks.  Maintaining compactness and the number of windings forces a commensurate change in the radius parameters of the building blocks, and this amounts to a change in the overall mean curvature.

\subsection{Acknowledgements}

I would like to thank Frank Pacard for suggesting this problem to me, providing invaluable guidance to me during its completion, and showing me excellent hospitality during my visits to Paris.  I would also like to thank Rafe Mazzeo and Rick Schoen for their support.

\section{Delaunay-Like CMC Hypersurfaces in \boldmath $\Sph^{n+1}$}
\mylabel{sec:delaunay}

\subsection{Introduction}

The Delaunay surfaces in $\R^3$ are the classical examples of non-compact, periodic, axially symmetry constant mean curvature surfaces.  They have been studied carefully by many authors, for example \cite{kapouleas7,kks,mazzeopacardpollack}.  In general, one can say that the Delaunay surfaces look roughly like a chain of spherical lobes connected to each other by neck regions.  There are two parameters that naturally describe the Delaunay surfaces, namely the mean curvature and the width of the neck at its narrowest point.  (This latter parameter is related in a complicated manner to the period of the Delaunay surface.)  There is a maximal value of the width parameter for a given mean curvature; and in this limit, the Delaunay surface is just a cylinder.  When the width parameter is small, the spherical lobes approximate true spheres and the neck regions approximate truncated, re-scaled catenoids.  In the limit when the width parameter approaches zero, the Delaunay surface approaches a union of osculating spheres whose points of contact are on a common line.   

The geometric picture of the Delaunay surfaces in $\R^3$ leads to the following idea.  One can attempt to construct ``Delaunay-like" CMC surfaces in $\Sph^3$, at least in the small-neck-width regime, by positioning mutually tangent  hyperspheres along an equator, connecting them together using truncated and re-scaled catenoidal necks, and then perturbing the resulting surface until it has constant mean curvature.  This construction will be realized in this section of the paper in the more general context of CMC hypersurfaces in $\Sph^{n+1}$ for arbitrary $n \geq 2$.  

\subsection{The Building Blocks of the Construction}
\mylabel{subsec:blocks}

\paragraph*{The  hyperspheres in \boldmath $\Sph^{n+1}$.}  Constant mean curvature hyperspheres inside $\Sph^{n+1}$ are obtained by intersecting $\Sph^{n+1}$ with affine hyperplanes.  Without loss of generality, suppose that $\R^{n+2}$ is written as $\R \times \R^{n+1}$ and given the coordinates $(x^0, x^1, \ldots, x^{n+1})$, while the hyperplane in question is $\Pi_a := \{x \in \R^{n+2} \: : \: x^0 = a \}$ for some $a \in [0,1)$.  For convenience, introduce an angle $\alpha \in (0, \pi/2)$ so that $\cos(\alpha) = a$.

\begin{defn}
	\mylabel{defn:sph}
	Up to $SO(n+2)$-rotation, the  \emph{affine hypersphere} with parameter $\alpha \in (0, \pi/2]$ in $\Sph^{n+1}$ is the hypersurface $S_\alpha := \Sph^{n+1} \cap \Pi_a$.  In other words, $$S_\alpha = \{ x \in \R^{n+2} \: : \: x^0 = \cos{\alpha} \:\: \mbox{and} \:\: (x^1)^2 + \cdots + (x^{n+1})^2 = \sin^2(\alpha) \} \, .$$
\end{defn}

It is essential for the constructions that follow to have a thorough understanding of the geometry of \emph{normal graphs} over $S_\alpha$.  Let $\btheta : \Sph^n \rightarrow \R^{n+1}$ be a parametrization of the unit sphere in $\R^{n+1}$.  Then one can parametrize $S_\alpha$ via $\btheta \longmapsto (\cos(\alpha) , \sin(\alpha) \btheta)$.  Furthermore, the displacement by a distance $\sigma$ along the geodesic normal to a point on $S_\alpha$ is found using the exponential map and is given by 
$$
\exp(\sigma N_\alpha)(\btheta) = (\cos(\alpha + \sigma), \sin( \alpha + \sigma) \btheta)
$$
where $N_\alpha$ is the unit outward normal of $S_\alpha$, given explicitly by $N_\alpha := - \cot(\alpha) P + \sin(\alpha) \frac{\partial}{\partial x^0}$ and $P := \sum_{k=1}^{n+1} x^k \frac{\partial }{\partial x^k}$ is the position vector field of the $\R^{n+1}$ factor.  Suppose now that $F: \Sph^n \rightarrow \R$ is a function on $\Sph^n$.  Then one can parametrize the normal graph over $S_\alpha$ corresponding to $F$ via
$$
\btheta \longmapsto \big( \cos(\alpha + F(\btheta)), \sin(\alpha + F(\btheta)) \btheta \big) \, .
$$
The geometric features of normal graphs over $S_\alpha$ that will be relevant later on are collected in the following lemma.

\begin{lemma}
	\mylabel{lemma:normgrgeom}
	The following facts hold for the normal graph over $S_\alpha$ corresponding to the function $F :\Sph^n \rightarrow \R$.  The norms and derivatives on the right hand sides of the expressions below are those of the standard metric on $\Sph^n$.  
	\begin{enumerate}
		\item The induced metric is
		$$g = \dif F \otimes \dif F + \sin^2 ( \alpha + F) g_{\Sph^n} \, .$$
		
		\item The Laplacian of $g$ satisfies
		\begin{align*}
			\sin^2( \alpha + F) \, \Delta_g (u) &= \Delta u - \frac{\nabla^2 u ( \nabla F, \nabla F)}{A^2} \\
			&\qquad + \frac{\big( (n-2) \sin(\alpha + F) \cos(\alpha + F) - \Delta F \big) \nabla F \cdot \nabla u}{A^2} \\
			&\qquad + \frac{\big( \nabla^2 F( \nabla F, \nabla F) + \sin(\alpha + F) \cos(\alpha+F) \| \nabla F \|^2 \big) \nabla F \cdot \nabla u}{A^4} 
		\end{align*}
		where $A := \big( \sin^2(\alpha + F) + \| \nabla F \|^2 \big)^{1/2}$.
		
		\item The unit normal vector field $N_\alpha$ is chosen to be
		\begin{equation*}
			AN_\alpha :=  \sin^2( \alpha + F) \frac{\partial }{\partial x^0} - \sum_{j=1}^{n+2} \left( \cos(\alpha +F) \sin(\alpha+F) \btheta^j - [g_{\Sph^n}]^{st} \nabla_s F  \frac{\partial \btheta^j}{\partial \nu^t} \right) \frac{\partial}{\partial x^j}%
		\end{equation*}
		where $\nu \longmapsto \Theta(\nu)$ parametrizes the unit $n$-sphere.
		
		\item The second fundamental form satisfies 
		\begin{equation*}
			A\,  B_{st} = - \sin(\alpha + F) \nabla_s  \nabla_t F + 2 \cos(\alpha + F) \nabla_s F \,\nabla_t F + \sin^2(\alpha + F) \cos(\alpha + F) [g_{\Sph^n}]_{st}  \, .
		\end{equation*}
		
		\item The mean curvature satisfies
		\begin{align*}
			A \sin(\alpha + F)  \, H &= - \Delta F + n \sin(\alpha + F) \cos(\alpha+F) \\
			&\qquad - \frac{\nabla^2 F (\nabla F, \nabla F) - \cos(\alpha + F) \sin(\alpha+F) \| \nabla F \|^2}{A^2} \, .
		\end{align*}
	\end{enumerate} 
\end{lemma}

Note that when $F \equiv 0$ is substituted into the various quantities in Lemma \ref{lemma:normgrgeom}, one obtains the corresponding quantities for $S_\alpha$ itself.  That is, the induced metric is $g_\alpha := \sin^2(\alpha) g_{\Sph^{n}}$ where $g_{\Sph^n}$ is the standard metric of the $n$-sphere, the second fundamental form is $B_\alpha := \cos(\alpha) \sin(\alpha) g_{\Sph^{n}}$ and mean curvature is $H_\alpha := n \cot(\alpha)$.    Therefore one sees that the second fundamental form of $S_\alpha$ is parallel and has constant mean curvature. 

\paragraph*{The generalized catenoid in \boldmath $\R^{n+1}$.}  The generalized catenoid in $\R^{n+1}$ is the $n$-dimensional analogue of the standard catenoid, namely the unique, cylindrically symmetric, minimal hypersurface in $\R^{n+1}$.  The derivation of the generalized catenoid (hereinafter referred to as simply the catenoid) is a fairly simple exercise in Riemannian geometry and will not be carried out here.   Since the zero mean curvature condition is invariant under scaling, there is a one-parameter family of such minimal hypersurfaces.  The definition of the scaled catenoid is as follows, in which $\R^{n+1}$ is written $\R \times \R^n$ using coordinates $(y^1, \ldots, y^{n+1})$ with $\hat y := (y^2, \ldots, y^{n+1})$.

\begin{defn}
	\mylabel{defn:gencat}
	The $\eps$-\emph{scaled catenoid} in $\R \times \R^n$ is the hypersurface $\eps \Sigma$ parameterized by
	$$
	( s , \Theta) \in \RÊ\times \Sph^{n-1} \longmapsto \eps ( \psi (s), \phi(s) \,  \Theta )	$$
	where $\phi(s) :=  (\cosh (n-1) s)^{1/(n-1)}$ and $\psi (s) : = \int_{0}^s \, \phi^{2-n} (\sigma) \, \dif \sigma$.  
\end{defn}

An alternate parametrization of the $\eps$-scaled catenoid will also be used in the sequel, namely when $\eps \Sigma$ is written as the union of two graphs over the $\R^n$ factor.  That is, by inverting the equation $\| \hat y \| = \eps \phi(s)$, one finds that  $\eps \Sigma = \Sigma_\eps^+ \cup \Sigma_\eps^-$ where $\Sigma_{\eps}^\pm := \big\{ ( \pm F_\eps (\| \hat y \|), \hat y ) : \Vert \hat y  \Vert \geq \eps \big\}$.  The function $F_\eps : \{ x \in \R : x \geq \eps \} \rightarrow \R$ is defined by $F_\eps (x) =  \eps F(x /\eps)$ where 
\begin{equation}
	\mylabel{eqn:catenoidgrfn}
	F(x) := \int_1^x (\sigma^{2n-2} - 1)^{-1/2} \dif \sigma \, .
\end{equation}
Note that in dimension $n=2$ this function is simply $F(x) = \arccosh(x)$.

The geometric features of the $\eps$-scaled catenoid that will be relevant later on are as follows.  The induced metric of $\Sigma$ is $g_{\eps \Sigma} := \eps^2 \phi^2 \, (\dif s^2 + g_{\Sph^{n-1}})$.  The unit normal vector field of $\Sigma$ is chosen as
$$
N_{\eps \Sigma} : = \frac{\dot \phi}{\phi} \, \frac{\partial}{\partial y^1} - \phi^{1-n} \, P_\Theta  \, ,
$$
where $P_\Theta$ is the position vector field of the $\R^n$ factor evaluated at the point $\Theta \in \Sph^{n-1}$. Then the second fundamental form of $\Sigma$ is given by 
$$
B_{\eps \Sigma} : = \eps \phi^{2-n} \, ((1-n) \, \dif s^2 + g_{\Sph^{n-1}})
$$
and satisfies
$$\Vert B_{\eps \Sigma} \Vert = \frac{C_1}{\eps \phi^{n}} \qquad \mbox{and} \qquad \Vert \nabla B_{\eps \Sigma} \Vert = \frac{C_2 \sinh((n-1) s)}{\eps^2 \phi^{2n}}$$
where $C_1 = \sqrt{n(n-1)}$ and $C_2 =  n \sqrt{(n+2)(n-1)}$.  The mean curvature of $\Sigma$ vanishes.

\subsection{Assembling the Approximate Solution}
\mylabel{subsec:assembly}

The approximately CMC ``Delaunay-like" hypersurface will be constructed roughly as follows.  First consider a single  hypersphere $S_\alpha$ with a fixed $\alpha \in (0,\pi/2)$.  Introduce a small \emph{spacing parameter} $\tau$ satisfying $0<\tau \ll 2\pi$.  Now position rotated copies of $S_\alpha$ in a symmetric manner around a geodesic $\gamma$ in $\Sph^{n+1}$ so that each hypersphere is separated from its two nearest neighbours by a distance $\tau$.  Depending on the choice of $\alpha$ and $\tau$, either this configuration of hyperspheres winds around $\gamma$ a finite number of times before repeating, or else it never does.  Now insert between the points of closest approach between every pair of neighbouring hyperspheres a suitably re-scaled, truncated and embedded catenoid.  The truncation and scaling parameter should be such that its height is $\mathcal O(\tau)$ and it makes optimally smooth contact with the hyperspheres on either side.  However, for this to be possible it turns out that it is necessary also to perturb the hyperspheres very slightly to improve their `fit' with the catenoids.  

Here are the details of this construction. 

\paragraph*{Positioning the hyperspheres.} The hyperspheres will first be positioned in the most symmetric  way possible using a cyclic subgroup rotations.  Define the rotation 
$$R_{\theta} = \left(
\begin{array}{cc|c}
	\cos(\theta) & - \sin(\theta) & 0 \\
	\sin(\theta) & \cos(\theta) & 0 \\
	\hline
	0 & 0 &I_{n}
\end{array} 
\right)
$$
where $I_{n}$ is the $n \times n$ identity matrix.  The rotation $R_\theta$ generates the geodesic $\gamma$ formed by intersecting $\Sph^n$ with the $(x^0,x^1)$-plane.  Straightforward trigonometry shows that any pair of rotated hyperspheres $R_{2 \alpha + \tau}^k(S_\alpha)$ and $R_{2 \alpha + \tau}^{k+1}(S_\alpha)$ are separated by a distance $\tau$.  The points of closest approach between these rotated hyperspheres are $R_{2 \alpha + \tau}^k (p^+) \in R_{2 \alpha + \tau}^k(S_\alpha)$ and $R_{2 \alpha + \tau}^{k+1} (p^-) \in R_{2 \alpha + \tau}^{k+1}(S_\alpha)$ where $p^\pm := R_\alpha^{\pm 1} (p) \in S_\alpha$ and $p$ is the point $p := (1, 0, \ldots, 0)$.  Furthermore, the set 
$$\Lambda_{\alpha, \tau}^{\#} := \bigcup_{k=0}^\infty R_{2 \alpha + \tau}^k (S_\alpha)$$ 
consists of either: a finite number of distinct rotated copies of $S_\alpha$ when $2 \alpha + \tau$ is a rational multiple of $2 \pi$; and an infinite number of distinct rotated copies of $S_\alpha$ otherwise.  In the former case,  if $2 \alpha + \tau = 2 \pi m / N$ for some pair of positive and relatively prime integers $m$ and $N$, then the number of copies of $S_\alpha$ equals $N$ and these $N$ copies of wind $m$ times around $\gamma$ before repeating. If $m=1$ then they do not intersect each other at all.

\paragraph*{The normal perturbation of the hyperspheres.}  Each hypersphere in the initial configuration must be perturbed slightly in the normal direction.  In order to preserve symmetry with respect to $R_{2\alpha+ \tau}$, the same perturbation will be used for each hypersphere.  Thus it suffices to explain how $S_\alpha$ is perturbed.  

This begins with the choice of a function $G: S_\alpha \rightarrow \R$ which determines the normal perturbation.  Recall that the linearized mean curvature operator on the space of normal graphs over a hypersurface $\Lambda$ in $\Sph^{n+1}$ is $\mathcal L_\Lambda := \Delta_\Lambda + \| B_\Lambda \|^2 + n$ where $\Delta_\Lambda$ is the Laplacian of $\Lambda$ and $B_\Lambda$ is the second fundamental form of $\Lambda$. In the case $\Lambda = S_\alpha$, then  $\mathcal L_\alpha := \sin^{-2}(\alpha) \big( \Delta_{\Sph^n} + n \big)$.  Choose a parametrization of $S_\alpha \setminus \{ p^+, p^- \}$ in which rotational symmetry around the geodesic $\gamma$ is in evidence; namely the parametrization 
\begin{equation}
	\mylabel{eqn:sphparam}
	(\mu, \Theta ) \in (0, \pi) \times \Sph^{n-1} \longmapsto \big( \cos(\alpha), \sin(\alpha) \cos(\mu), \sin(\alpha) \sin(\mu) \Theta \big)
\end{equation}
where $\Theta : \Sph^{n-1} \rightarrow \R^{n}$ is a parametrization of $\Sph^{n-1}$ of the unit sphere in $\R^n$.  Now let $G$ be the $\Theta$-independent solution of the equation $\mathcal L_\alpha (G) = 0$ which is singular at $\mu = 0$ and $\mu = \pi$ and symmetric with respect to $\mu \mapsto \pi - \mu$.  Explicitly, this function is
\begin{equation}
	\mylabel{eqn:green}
	G(\mu) :=
		- \sin(\mu) - \cos(\mu) \displaystyle \int_{\pi/2}^\mu \frac{1 - \sin^{n-1}(\sigma)}{\cos^2(\sigma) \sin^{n-1}(\sigma)} \dif \sigma \, .
\end{equation}
Important for later on is the asymptotic expansion of $G$ at $\mu = 0$.  This is
\begin{equation}
	\mylabel{eqn:normgraphasym}
	G(\mu) = 
	\begin{cases}
		-1 + \log 2 -\log(\mu) + \mathcal O (\mu^2 |\log(\mu)|) &\qquad n=2 \\[1.5ex]
		\dfrac{1}{\mu} + \mathcal O (\mu) &\qquad n=3 \\[1.5ex]
		\dfrac{1}{2 \mu^2} + \mathcal O(1+|\log(\mu)|) &\qquad n=4 \\[1.5ex]
		\dfrac{1}{(n-2) \mu^{n-2}} + \mathcal O ( \mu^{4-n} ) &\qquad n \geq 5
	\end{cases}
\end{equation}
when $0 < \mu \ll \pi/2$.  

The desired normal perturbation can be defined as follows.  Choose a small parameter $\eps > 0$ (whose relation to $\tau$ and $\alpha$ will be specified in due course) and define the normal graph $$\tilde S_{a,\eps} := \exp(\eps^{n-1} G N_\alpha) (S_\alpha \setminus \{ p^+, p^- \})$$ where $N_\alpha$ is the outward unit normal vector field of $S_\alpha$.  For concreteness, this hypersurface has the parametrization
\begin{equation}
	\mylabel{eqn:greensurf}
	(\mu, \Theta) \longmapsto \big(\cos(\alpha + \eps^{n-1} G(\mu)), \sin (\alpha + \eps^{n-1} G(\mu)) \cos(\mu), \sin(\alpha + \eps^{n-1} G(\mu)) \sin(\mu) \Theta \big)
\end{equation}
as in the discussion preceding Lemma \ref{lemma:normgrgeom}.
 
\paragraph*{Stereographic coordinates adapted to a pair of hyperspheres.} Canonical coordinates that are well-adapted to the pair of hyperspheres $R_{2 \alpha + \tau}^{k} (S_\alpha)$ and $R_{2 \alpha + \tau}^{k+1} (S_\alpha)$ will be needed.  These can be defined as follows.  First, note that the point $R_{\alpha + \tau /2}^{2k+1}(p)$ lies on the geodesic $\gamma$ at the midpoint between these two hyperspheres.  Now let $K : \Sph^{n+1} \setminus \{ - p \} \rightarrow \R^{n+1}$ denote the stereographic projection centered at $p$ defined by
$$ K(x^0, x^1, \ldots, x^{n+1}) := \left( \frac{x^1}{1+x^0} , \cdots, \frac{x^{n+1}}{1+x^0} \right)\, .
$$
Then the desired adapted coordinates are given by the inverse of the mapping $K \circ R_{\alpha + \tau /2}^{-(2k+1)} : \Sph^{n+1} \setminus \{ -R_{\alpha + \tau /2}^{2k+1}(p) \} \rightarrow \R^{n+1}$. 

Recall that stereographic projection sends non-equatorial $k$-spheres in $\Sph^{n+1}$ to $k$-spheres in $\R^{n+1}$ and sends equatorial $k$-spheres to linear subspaces.  One thus expects that the coordinate image of the geodesic $\gamma$ is the $y^1$-axis.  Also, one expects that the coordinate images of the two hyperspheres $R_{2 \alpha + \tau}^{k} (S_\alpha)$ and $R_{2 \alpha + \tau}^{k+1} (S_\alpha)$ are two hyperspheres symmetrically located on either side of the origin and centered at two points on the  $y^1$-axis.  Indeed, one can check that any point of $\Sph^{n+1}$ of the form \eqref{eqn:sphparam} maps to the point $(y^1, \hat y)$ in $\R^{n+1}$ given by
\begin{equation}
	\mylabel{eqn:stereo}
	\begin{aligned}
		y^1 & = \frac{- \sin(\alpha + \tau/2) \cos(\alpha ) + \cos(\alpha + \tau/2) \sin(\alpha ) \cos(\mu)}{1 + \cos(\alpha + \tau/2) \cos(\alpha) + \sin(\alpha + \tau/2) \sin(\alpha) \cos(\mu)}\\[1ex]
		\hat y &=  \frac{\sin(\alpha ) \sin(\mu) \Theta}{1 + \cos(\alpha + \tau/2) \cos(\alpha) + \sin(\alpha + \tau/2) \sin(\alpha ) \cos(\mu)} \, .
	\end{aligned}
\end{equation}
And, it is easy to check that this point lies on the locus of points satisfying the equation
\begin{equation}
	\mylabel{eqn:sphimage}
	\left( y^1 + d \right)^2 + \Vert \hat y \Vert^2 =  r^2
\end{equation}
where 
\begin{equation}
	\mylabel{eqn:raddisp}
	\begin{aligned}
		r &= r(\alpha, \tau) := \frac{\sin(\alpha)}{ \cos(\alpha) + \cos(\alpha + \tau/2)} \\
		d &= d(\alpha, \tau) := \frac{ \sin(\alpha + \tau/2)}{\cos(\alpha) + \cos(\alpha + \tau /2)} \, .
	\end{aligned}
\end{equation} 
Observe that $d - r = \tan(\tau/4)$ is displacement from the origin of the hypersphere determined by \eqref{eqn:sphimage} in these coordinates.  This yields a geometric displacement of $\tau/2$ when measured with respect to the metric $A^{-2} g_0$, as expected.

The reason for using stereographic coordinates, other than the fact that spheres map to spheres, is that these are in geodesic normal form: at the centre of the coordinates, the metric is Euclidean and its derivatives vanish.  This can be seen by the computation of the metric $$\big( R_{\alpha + \tau/2}^{2k+1} \circ K^{-1} \big)^\ast g_{\Sph^{n+1}} = A^{-2} g_0$$ where $g_0$ is the Euclidean metric of $\R^{n+1}$ and $A(y) = \frac{1}{2} (1 + \sum_{k=1}^{n+1} (y^k)^2)$ using the coordinates $(y^1, \ldots, y^{n+1})$ for $\R^{n+1}$.   The geodesic normal form will have the effect of distorting as little as possible the geometry of objects embedded into the sphere using the stereographic coordinate map, provided one remains near the origin.  The following lemma contains the results of this kind that will be used in the sequel.

\begin{lemma}
	\mylabel{lemma:stereomeancurv}
	Let $\sigma : \Sigma \rightarrow \R^{n+1}$ be an embedding of the submanifold $\Sigma$ into $\R^{n+1}$.  Let $h_0$, $\mathring \Delta$, $N_0$, $\mathring B_{ij}$ and $H_0$ be the induced metric, Laplacian, unit outward normal vector field, components of the second fundamental form and mean curvature of $\sigma(\Sigma)$ with respect to the Euclidean metric.  The corresponding quantities with respect to the metric $g = A^{-2} g_0$ coming from stereographic projection satisfy the following identities.
	\begin{enumerate}
		\item The induced metric of $\sigma(\Sigma)$ with respect to $g$ is $h = [A \circ \sigma]^{-2} h_0$.  The Laplacian of $\sigma(\Sigma)$ with respect to $h$ is
		\begin{equation*}
			\Delta u = [A \circ \sigma]^2 \mathring \Delta u - \sum_{k=1}^{n+1} (n-2) [A \circ \sigma] \, h_0( \sigma ^k \mathring \nabla \sigma ^k , \mathring \nabla u) 
		\end{equation*}
		where $\sigma^k$ are the components of $\sigma$.
				 				
		\item The second fundamental form of $\phi(\Sigma)$ with respect to the metric $g$ is
		\begin{equation*}
			B_{ij} = [A \circ \sigma ]^{-1}  \mathring B_{ij} - [A \circ \sigma ]^{-2} g_0 (\sigma, N_0) [g_0]_{ij} \, .
		\end{equation*}
		\item The mean curvature of $\sigma(\Sigma)$ with respect to the metric $g$ is
		\begin{equation*}
			H = [A \circ \sigma]  H_0 - n g_0( \sigma, N_0 ) \, .
		\end{equation*}
	\end{enumerate}	
\end{lemma}

\paragraph*{Inserting truncated catenoids.}  The pair of perturbed hyperspheres $R_{2 \alpha + \tau}^k (\tilde S_\alpha)$ and $R_{2 \alpha + \tau}^{k+1} (\tilde S_\alpha)$ maps to a pair of hypersurfaces which are small deformations of $R_{2 \alpha + \tau}^k (S_\alpha)$ and $R_{2 \alpha + \tau}^{k+1} (S_\alpha)$.  It is necessary to have a precise description of the shape of these hypersurfaces near the $y^1$-axis in order to determine to what extent they resemble the ends of a catenoid.  This will be done by reparametrizing them as graphs over the $\hat y$-hyperplane in a small neighbourhood of the $y^1$-axis. 

From equation \eqref{eqn:greensurf} and the properties of the stereographic projection following equation \eqref{eqn:stereo}, one finds that the coordinates $y(\mu, \Theta) \in \R^{n+1}$ of a point in the stereographic projection of the hypersphere $R_{2 \alpha + \tau}^{k} (\tilde S_\alpha)$ satisfy 
\begin{equation}
	\mylabel{eqn:vgraph}
	y^1(\mu) = -D(\mu) + \sqrt{[R(\mu)]^2 - \| \hat y \|^2}
\end{equation}
where
\begin{align*}
	R(\mu) &:=  \frac{\sin(\alpha+ \eps^{n-1} G (\mu))}{ \cos(\alpha+ \eps^{n-1} G (\mu)) + \cos(\alpha + \tau/2)} \\[0.5ex]
	D(\mu) &:=  \frac{ \sin(\alpha + \tau/2)}{\cos(\alpha+ \eps^{n-1} G (\mu)) + \cos(\alpha + \tau /2)} \, .
\end{align*}
Furthermore, the relation between $\mu$ and $\| \hat y \|$ is given by 
\begin{equation}
	\mylabel{eqn:invrel}
	\| \hat y \| = \frac{\sin(\alpha+ \eps^{n-1} G (\mu) ) \sin(\mu)}{1 + \cos(\alpha + \tau/2) \cos(\alpha+ \eps^{n-1} G (\mu)) + \sin(\alpha + \tau/2) \sin(\alpha+ \eps^{n-1} G (\mu) ) \cos(\mu)} \, .
\end{equation}
By computing the derivative $\frac{\dif}{\dif \mu} \| \hat y \|$, one finds that the relation \eqref{eqn:invrel} is invertible in the region where both $\eps^{n-1} G(\mu)$ and $\mu$ are small, and moreover that $\mu(\| \hat y \|) = 2 \csc(\alpha) \cos^2(\tau/4) \| \hat y \| + \mathcal O(\| \hat y \|^3)$.  Substituting this into \eqref{eqn:vgraph} yields $y^1(\| \hat y \|) = G_\eps(\| \hat y \|) := D( \mu(\| \hat y \|)) - \sqrt{[R(\mu(\| \hat y \|))]^2 - \| \hat y \|^2}$ from which one can find the asymptotic expansion 
\begin{equation}
	\mylabel{eqn:greensurfasym}
	G_\eps(\| \hat y \|) =
	\begin{cases}
		\begin{aligned}
			&\hspace{-0.875ex} - \tan(\tau/4) - \frac{\| \hat y \|^2}{2r} + \eps c_2 - \eps C_2 \log ( \| \hat y \|) \\[-0.25ex]
			&\hspace{-0.875ex} \qquad \qquad + \mathcal O (\| \hat y \|^4) + \mathcal O(\eps \| \hat y \|^2 | \log(\| \hat y \|)| ) 
		\end{aligned} 
		&\quad \mbox{$n=2$}  \\[4ex]
		\displaystyle -\tan(\tau/4)  - \frac{\| \hat y \|^2}{2r}  + \frac{\eps^{2} C_3}{ \| \hat y \|} + \mathcal O (\| \hat y \|^4) + \mathcal O \left( \eps^2 \| \hat y \| \right) &\quad \mbox{$n=3$}\\[1.5ex]
		\displaystyle - \tan(\tau/4)  - \frac{\| \hat y \|^2}{2r}   + \frac{\eps^{3} C_4}{ 2  \| \hat y \|^{2}} + \mathcal O (\| \hat y \|^4)  +\mathcal O \left( \eps^{3} \bigl[ 1 + | \log(\| \hat y \|) | \bigr]  \right) &\quad \mbox{$n=4$} \\[1.5ex]
		\displaystyle -\tan(\tau/4)  - \frac{\| \hat y \|^2}{2r}  + \frac{\eps^{n-1} C_n}{(n-2) \| \hat y \|^{n-2}} + \mathcal O (\| \hat y \|^4) + \mathcal O \left( \frac{\eps^{n-1}}{\| \hat y \|^{n-4}} \right) &\quad \mbox{$n\geq 5$}
	\end{cases}
\end{equation}
in the region where both $\| \hat y \|$ and $\eps^{n-1} G(\mu(\| \hat y \|))$ remain small.  Here $d = d(\alpha, \tau)$ and $r = r(\alpha, \tau)$ are the quantities in \eqref{eqn:raddisp}, while 
\begin{equation}
	\label{eqn:expansionconst}
	\begin{aligned}
		c_2 &=  \frac{\cos(2 \alpha + \tau/2) \bigl(1  + \log \bigl( \csc(\alpha) \cos^2(\tau/4) \bigr) \bigr)}{\bigl( \cos(\alpha) + \cos(\alpha + \tau/2) \bigr)^2} \\[0.5ex]
		C_n &= \frac{\cos(2 \alpha + \tau/2)}{\bigl( \cos(\alpha) + \cos(\alpha + \tau/2) \bigr)^2 \bigl( 2 \csc(\alpha) \cos^2(\tau/4) \bigr)^{n-2}}
	\, .
	\end{aligned}
\end{equation}
In a similar manner, one finds that the equation satisfied by the image under stereographic projection of the other hypersphere $R_{2 \alpha + \tau}^{k+1} (\tilde S_\alpha)$ is $y^1(\| \hat y \|) = - G_\eps(\| \hat y \|)$.

The task is now to find a truncation and re-scaling of the catenoid that `fits' exactly within the gap between the two perturbed hyperspheres.  The precise dimensions of this catenoid will be found by matching the asymptotic expansion of the graphing function given in \eqref{eqn:catenoidgrfn} to the expansion \eqref{eqn:greensurfasym} in the highest-order terms.  

The upper and lower halves of the catenoid scaled by a factor $\tilde \eps$ are the graphs $y^1 = \pm \tilde \eps F( \| \hat y \| / \tilde \eps)$ where $F(x) = \int_1^x \big( \sigma^{2n-2} - 1 \big)^{-1/2} \dif \sigma$.  Hence one has the asymptotic expansion
\begin{equation}
	\mylabel{eqn:catenoidasym}
	\tilde \eps F( \| \hat y \| / \tilde \eps) = 
	\begin{cases}\displaystyle
		\tilde \eps  \log (2 / \tilde \eps)  + \tilde \eps \log ( \| \hat y \| ) - \frac{\tilde \eps^3}{4 \| \hat y \|^2 } + \mathcal O \left( \frac{ \tilde \eps^5}{\| \hat y \|^4} \right) &\quad \mbox{$n=2$} \\[1.5ex]
		\displaystyle \tilde \eps c_n - \frac{\tilde \eps^{n-1}}{(n-2) \| \hat y \|^{n-2}} - \frac{\tilde \eps^{3n-3}}{2(3n-4) \| \hat y \|^{3n-4}} + \mathcal O \left( \frac{\tilde \eps^{5n-5}}{\| \hat y \|^{5n-6}} \right) &\quad \mbox{$n \geq 3$}
	\end{cases}
\end{equation}
where $c_n$ is yet another constant.  If one  compares the asymptotic expansion \eqref{eqn:catenoidasym} (multiplied by $-1$) with the asymptotic expansion \eqref{eqn:greensurfasym}, then the matching is optimal if: $\tilde \eps = \eps$ in dimension $n=2$ and $\tilde \eps = \eps C_n^{(n-2)/(n-1)}$ in higher dimensions; and $\eps$ is chosen to satisfy the equation $\tilde \eps \log(2 /\tilde \eps) = \tan(\tau/4) + \eps c_2$ in dimension $n=2$ and $\eps c = \tan(\tau/4)$ in higher dimensions.

The scale factor $\eps := \eps(\tau)$ is determined by the considerations above. Furthermore, the location where the catenoid must be truncated in order for the gluing to take place with optimal error can be determined by similar considerations.  That is, once $\eps$ and $\tilde \eps$ have been found, then the error $|F(\| \hat y \|) - G_\eps (\| \hat y \|)|$ is smallest when $\hat y$ is chosen to lie in a range where the quantity $\frac{1}{2r} \| \hat y \|^2 + \frac{\eps^{3n-3}}{2(3n-4)} \| \hat y \|^{4-3n}$ is minimized.  This occurs when $\| \hat y \| = \mathcal O(\rho_\eps)$ where $\rho_\eps := \eps^{(3n-3)/(3n-2)}$.  Thus one must truncate the $\eps$-scaled catenoid exactly at $\| \hat y \| = \rho_\eps$ for an optimally smooth gluing.

\paragraph*{Assembling the approximate solution.}  The first step is to construct a smooth hypersurface with boundary in $\R^{n+1}$ that interpolates between the stereographic coordinate images of a pair of perturbed hyperspheres near the $y^1$-axis.   Denote these objects by $S_\pm$.  Now let $\eta : [0, \infty) \rightarrow \R$ be a smooth, monotone cut-off function satisfying $\eta(s) = 0$ for $s \in [0, 1/2]$ and $\eta(s) = 1$ for $s \in [2, \infty)$.  Define the function $\tilde F_{\alpha, \tau} : \bar B_{2 \rho_\eps}(0) \setminus B_\eps (0) \subseteq  \R^{n} \rightarrow \R$ by 
\begin{equation}
	\mylabel{eqn:mergedfn}
	\tilde F_{\alpha, \tau} (\hat y) = \eps \big( 1 - \eta(\Vert \hat y\Vert / \rho_\eps) \big) F(\Vert \hat y\Vert / \eps) + \eta(\Vert \hat y\Vert /\rho_\eps) G_\eps (\| \hat y \| ) \, .
\end{equation}
Then define the hypersurfaces $\tilde \Sigma^\pm_\eps = \{ ( \pm \tilde F_{a, \tau}(\hat y), \hat y) : \Vert\hat y\Vert \in [\eps, 2\rho_\eps] \}$ so that $\tilde \Sigma_\eps := \tilde \Sigma_\eps^+ \cup \tilde \Sigma_\eps^-$ is a smooth hypersurface connecting $S_+ \setminus ( \R \times B_{2 \rho_\eps}(0) )$ to  $S_- \setminus ( \R \times B_{2 \rho_\eps}(0) )$ through the catenoid.

Pushing forward to the sphere, one now has a small catenoidal hypersurface $R_{\alpha + \tau /2}^{2k+1} \circ K^{-1} ( \tilde \Sigma_\eps)$ that overlaps perfectly with the pair of perturbed hyperspheres $R_{2 \alpha + \tau}^{k}(\tilde S_\alpha)$ and $R_{2 \alpha + \tau}^{k+1}(\tilde S_\alpha)$ inside a small tubular annulus of the geodesic $\gamma$ near the point $R_{ \alpha + \tau/2}^{2k+1}( p)$.  

The second step is to connect the hypersurfaces above together according to the following definition.  Note that there exists a radius $\tilde \rho_\eps$ so that the boundary of 
$R_{\alpha + 2 \tau}^k(S_\alpha \setminus B_{\tilde \rho_\eps}(p^+)) \cup R_{\alpha + 2 \tau}^{k+1}(S_\alpha \setminus B_{\tilde \rho_\eps}(p^-))$ under the stereographic projection $R_{\alpha+ \tau/2}^{2k+1} \circ K$ coincides with $\partial (\tilde \Sigma_\eps^+ \cup \tilde \Sigma_\eps^-)$.  Clearly $\tilde \rho_\eps = \mathcal O(\rho_\eps)$.

\begin{defn}
\mylabel{defn:approxsol}
The \emph{approximate solution} with parameters $\alpha$ and $\tau$ is the hypersurface
\begin{equation}
	\mylabel{eqn:approxsol}
	\tilde \Lambda_{\alpha,\tau} :=  \left[ \bigcup_{k=0}^{\infty} R_{2\alpha + \tau}^k \big( \tilde S_\alpha \setminus \big( B_{\tilde \rho_\eps} (p^+) \cup B_{\tilde \rho_\eps} (p^-) \big) \big) \right] \cup  \left[ \bigcup_{k=0}^\infty R_{\alpha + \tau /2}^{2k+1} \circ K^{-1} ( \tilde \Sigma_\eps) \right] 
\end{equation}
where $\eps = \eps(\tau)$ is the scale parameter associated to $\tau$.
\end{defn}

Note that when the angle $2 \alpha + \tau$ is a rational multiple of $2 \pi$, then $\approxsol$ is a compact immersed hypersurface of $\Sph^n$ that winds around the geodesic $\gamma$ an integral  number of times.  In this case, $R_{2 \alpha + \tau}$ generates a finite cyclic subgroup of $SO(n+1)$ of order $N$ and the infinite unions in \eqref{eqn:approxsol} can be replaced by finite unions from $k=0$ to $k=N-1$.  Furthermore, $\approxsol$ is embedded if and only if it wraps around $\gamma$ exactly once.    If $2 \alpha + \tau$ is not a rational multiple of $2 \pi$, then $\approxsol$ is non-compact. It  will henceforth be assumed that $2 \alpha + \tau$ is a rational multiple of $2 \pi$.  The analysis that follows can readily be adapted to the case of general $2 \alpha + \tau$.  

This section of the paper concludes with the following identification of the various regions of $\approxsol$.  Let $\mathit{Cyl}(\rho) := \{ (y^1, \hat y ) \in \R \times \R^n : \| \hat y \| < \rho \}$ denote the cylinder of radius $\rho$ parallel to the $y^1$-axis in $\R \times \R^n$.
\begin{defn}
	The approximate solution $\approxsol$ is divided into the following three regions.
	\begin{itemize}
		\item Let $\mathcal N_\eps^k :=  R_{\alpha + \tau /2}^{2k+1} \circ K^{-1} \big( \tilde \Sigma_\eps \cap \mathit{Cyl}({\rho_\eps /2}) \big)$.  Then $\mathcal N_\eps := \bigcup_{k=0}^\infty \mathcal N_\eps^k$ is the \emph{neck region} of $\approxsol$.
		
		\item Let $\mathcal T_{\eps}^{k, \pm} :=  R_{\alpha + \tau /2}^{2k+1} \circ K^{-1} \big( \tilde \Sigma_\eps^\pm \cap\big[ \mathit{Cyl}({2\rho_\eps}) \setminus \mathit{Cyl}({\rho_\eps/2}) \big] \big)$.   Then $\mathcal T_\eps := \bigcup_{k=0}^\infty \big[ \mathcal T_\eps^{k, +} \cup \mathcal T_\eps^{k, -} \big]$ is the \emph{transition region} of $\approxsol$.
		
		\item Let $\mathcal E_\eps^k :=  R_{2\alpha + \tau}^k \big( \tilde S_\alpha \setminus \big( B_{\tilde \rho_\eps} (p^+) \cup B_{\tilde \rho_\eps} (p^-) \big) \big)$.  Then $\mathcal E_\eps := \bigcup_{k=0}^\infty \mathcal E_\eps^k= \approxsol \setminus \big[ \mathcal N_\eps \cup \mathcal T_\eps \big]$ is the \emph{exterior region} of $\approxsol$.
	\end{itemize}
\end{defn}

\subsection{Symmetries of the Approximate Solution}
\mylabel{sec:syms}
  
It is clear that $\approxsol$ is invariant under the group generated by $R_{2 \alpha + \tau}$.  The additional symmetries possessed by $\approxsol$ will play an important role in the forthcoming analysis.   The first set of symmetries is defined as follows.  Define $S^{\, 01}_{B} \in O(n+2)$ by choosing $B \in O(n)$ and then setting
$$S^{\, 01}_B = \left( 
\begin{array}{c|c} 
	\! \! \mbox{\scriptsize$\begin{array}{cc}
		1 & \\[-0.5ex]
		& 1
	\end{array}$} \! \!& 0\\ 
	\hline 0 & B 
\end{array}
\right)
$$
which is the transformation keeping the $x^0$ and $x^1$ coordinates fixed while rotating the remaining coordinates by $B$. That each $S^{\, 01}_{B}$ is a symmetry of $\approxsol$ can be seen as follows.  First, $S^{\, 01}_{B}$ commutes with $R_{2 \alpha + \tau}$  so that it is sufficient to check that one perturbed hypersphere, one neck region $\mathcal N^0_\eps$ and one pair of transition regions $\mathcal T_\eps^{0,\pm}$ is preserved by $S^{\, 01}_{B}$.  But this holds because the construction of the interpolating catenoidal neck guarantees that rotational symmetry around the $y^1$-axis in the stereographic coordinates is preserved.

One additional symmetry must be identified.  Define the reflection $T \in O(n+2)$ by $$T(x^0, x^1, x^2, \ldots, x^{n+1}) = (x^0, - x^1, x^2, \ldots, x^{n+1}) \, .$$ That $T$ is also a symmetry of $\approxsol$ can be seen as follows.  First, $T \circ R_{2 \alpha + \tau} = R_{2 \alpha + \tau}^{-1} \circ T$ so that $T$ permutes the perturbed hyperspheres and the interpolating catenoidal necks amongst each other, perhaps with opposite orientation.  The change of orientation doesn't matter for the perturbed hyperspheres because the graphing function $G$ satisfies $G(\mu) = G(\pi - \mu)$.  Morevoer, the catenoidal necks and transition regions were constructed so as to be preserved under the transformation $y^1 \longmapsto - y^1$ in the stereographic coordinates.  Hence the change of orientation doesn't matter for these either.  Thus $T$ is also a symmetry of $\approxsol$.

\subsection{The Analytic Set-Up}

\mylabel{sec:jacobi}

\paragraph*{Deforming the approximate solution.} The approximate solution $\approxsol$ constructed in the Section \ref{subsec:assembly} has mean curvature exactly equal to $H_\alpha$ everywhere except in the neck and transition regions where the mean curvature becomes to zero (but perhaps in non-uniform way due to the errors introduced by the cut-off functions).  The next task is to develop a means for deforming $\approxsol$ into an exactly constant mean curvature hypersurface, with mean curvature equal to $H_\alpha$.  

Since $\approxsol$ is a hypersurface, it is possible to parametrize deformations of $\approxsol$ in a very standard way via normal deformations.  These can be constructed by choosing a function $f : \approxsol \rightarrow \R$ and then considering the deformation $\phi_f : \approxsol \rightarrow \Sph^{n+1}$ given by $\phi_f (q) := \exp_q (f(q) \cdot N(q))$ where $\exp_q$ is the exponential map at the point $q$ and $N(q)$ is the outward unit normal vector field of $\approxsol$ at the point $q$.  For any given function $f$, the hypersurface $\phi_f(\approxsol)$ is a normal graph over $\approxsol$, provided $f$ is sufficiently small in a $C^1$ sense.  Finding an exactly CMC normal graph near $\approxsol$ therefore consists of finding a function $f$ satisfying the equation $H_{\phi_f(\approxsol)} = H_\alpha$, where $H_\Lambda$ denotes the mean curvature of a hypersurface $\Lambda$.  

\begin{defn} 
Let  $\difop$ be the operator $f \longmapsto H_{\phi_f(\approxsol)} - H_\alpha $.  
\end{defn}

\noindent This is a quasi-linear, second-order partial differential operator for the function $f$ whose zero gives the desired deformation of $\approxsol$.  

\paragraph*{The Banach space inverse function theorem.} Finding a solution of the equation $\difop(f) = 0$ when $\tau$ is sufficiently small will be accomplished by invoking the \emph{Banach space inverse function theorem}.   To provide a focus for the remainder of the proof of Main Theorem \ref{result1}, this fundamental result will be stated here in fairly general terms.  See \cite{amr} for a proof.

\begin{nonumthm}[IFT]
	Let $\Phi : X \rightarrow Z$ be a smooth map of Banach spaces, set $\Phi(0) := E$ and define the linearized operator $\mathcal L := \Dif \Phi(f) =  \left. \frac{\dif}{\dif s} \Phi(f+ s u)  \right|_{s=0}$.  Suppose $\mathcal L $ is bounded and either $\mathcal L$ is invertible and satisfies 
	\begin{subequations}
	\mylabel{eqn:iftestone}
	\begin{equation}
		\Vert \mathcal L^{-1} (z)\Vert \leq C \Vert z \Vert
	\end{equation}
	for all $z \in Z$; or else $\mathcal L$ is surjective and possesses a bounded right inverse $\mathcal R : Z \rightarrow X$ satisfying
	\begin{equation}
		\Vert \mathcal R (z)\Vert \leq C \Vert z \Vert
	\end{equation}
	\end{subequations}
	for all $z \in Z$.	 Choose $\rho$ so that if $y \in B_\rho(0) \subseteq X$, then 
	\begin{equation}
		\mylabel{eqn:iftesttwo}
		\Vert \mathcal L( x) - \Dif \Phi(y) (x) \Vert \leq \frac{1}{2C}  \Vert x \Vert
	\end{equation}
for all $x \in X$, where $C>0$ is a constant.  Then if $z \in Z$ is such that 
	\begin{equation}
		\mylabel{eqn:iftestthree}
		\Vert z - E \Vert \leq \frac{\rho}{2C} \, ,
	\end{equation}
there exists a unique $x \in B_\rho(0)$ so that $\Phi(x) = z$.  Moreover, $\Vert x \Vert \leq 2 C \Vert z - E \Vert$. 
\end{nonumthm}

As the statement of theorem makes clear, it must be the case that $\linop$ is surjective in order to solve the equation $\difop(f) = 0$ in a Banach subspace $X$ of at least $C^2$ functions on $\approxsol$. In addition, it is necessary to establish the three fundamental estimates \eqref{eqn:iftestone}, \eqref{eqn:iftesttwo}  and \eqref{eqn:iftestthree} of the theorem.  It will turn out that $\Dif \difop(0)$ is invertible in the setting of Main Theorem 1 and the first of these estimates, called the \emph{linear estimate} is proved in Section \ref{sec:linest} after some preliminary work in Sections \ref{sec:funcspace} and \ref{sec:jacobi} that helps identify the correct Banach subspace $X$.  The second estimate determines the variation of the operator $\Dif \difop(f)$ as $f$ varies.  The third estimate is of the size of $\difop(0)$ in the norm of $Z$.  The latter two estimates, called the \emph{non-linear estimates}, are proved in Section \ref{subsec:nonlinest}.  In all three estimates, it will be necessary to determine explicitly the $\tau$-dependence of the constant $C$ in order to apply them in a uniform way when $\tau$ is close to zero.

\paragraph*{The linearized operator and the Jacobi fields.}
In order to solve the equation $\difop (f) = 0$, the linearized operator $\Dif \difop (0)$ (hereinafter abbreviated $\linop$) must satisfy two fundamental requirements: it must be surjective; and it must possess a right inverse bounded above by a constant independent of $\tau$.  It is however a general phenomenon in singular perturbation problems that the linearized operator often has a kernel and a co-kernel, as well as small eigenvalues tending to zero as the singular parameter (in this case $\tau$) tends to zero.  These are the \emph{obstructions} preventing the deformation to an exactly CMC hypersurface.  

The origin of these obstructions in the present case is geometric.   To see this, recall the general fact that any one-parameter family of isometries of the ambient space in which a CMC hypersurface is situated gives rise to an element in the kernel of the linearized mean curvature operator as follows.

\begin{lemma}
	\mylabel{lemma:linmeancurv}
	Let $\Lambda$ be a closed hypersurface in a Riemannian manifold $X$ with mean curvature $H_\Lambda$, second fundamental form $B_\Lambda$ and unit normal vector field $N_\Lambda$.   If $R_t$ is a one-parameter family of isometries of $X$ with deformation vector field $V = \left. \frac{\dif }{\dif t} R_t \right|_{t=0} $, then the function $q_V := \langle V, N_\Lambda \rangle$ belongs to the kernel of $\Lambda$.
\end{lemma}

\begin{proof}
	Since $R_t$ is a family of isometries, then $H(R_t(\Lambda) ) = H(\Lambda)$ for all $t$ and $\left. \frac{\dif}{\dif t} \right|_{t=0}H(R_t(\Lambda)) = 0$.  The function $q_V := \langle V, N_\Lambda\rangle$ is thus in the kernel of $\Dif H_\Lambda(0)$ because $q_V$ generates a normal deformation of $\Lambda$ whose action, to first order, coincides with $R_t$.
\end{proof}

Both kinds of obstructions to the solvability of the equation $\difop(f) = 0$ on $\approxsol$ (namely those corresponding to elements in the kernel and those corresponding to small eigenvalues) can be explained using Lemma \ref{lemma:linmeancurv}.  First of all, since rotations of $\Sph^{n+1}$ are isometries, every one-parameter rotation of $\approxsol$ generates a function in the kernel of $\linop$ exactly as in the lemma.  These functions are known as \emph{Jacobi fields}.  Second of all, one can imagine transformations of $\approxsol$ which rotate exactly one of its constituent hyperspheres or catenoidal necks by a rotation in $SO(n+2)$ while leaving all the other constituent hyperspheres and necks fixed. The associated function $q_V$ coincides with an exact Jacobi field on the constituent being transformed and vanishes on the other constituents of $\approxsol$.   It is well known that the linear span of these functions approximates the small eigenspaces of $\linop$ \cite[Appendix B]{kapouleas7}.  These functions will be called  \emph{approximate Jacobi fields}. 

The existence of these two types of obstructions makes it impossible to meet the requirements for solving for the perturbation of $\approxsol$ into a CMC hypersurface.  One way to avoid these obstructions is to exploit the natural symmetry of $\approxsol$ and then to deform $\approxsol$ equivariantly; i.e.~by normal deformations corresponding to functions on $\approxsol$ that are invariant under these symmetries.  The key is to have enough symmetries so that there are no Jacobi fields --- either approximate or true --- that are invariant with respect to all of these symmetries at once.  That this situation holds for $\approxsol$ will be shown in Section \ref{sec:linest}.

To continue, it will be necessary to have an explicit representation of the Jacobi fields and approximate Jacobi fields of $\approxsol$.  First, it is clear what the true Jacobi fields of $\approxsol$ look like: they are of the form $\langle V, N \rangle$ where $V$ generates a rotation and $N$ is the unit normal vector field of $\approxsol$.  Second, according to \cite[Appendix B]{kapouleas7}, the approximate Jacobi fields can be constructed by multiplying the function $\langle V, N \rangle$ by a cut-off function whose support is exactly one of the constituent hyperspheres or catenoidal necks of $\approxsol$.  It is thus necessary to understand the functions $\langle V, N \rangle$, as $V$ ranges over all generators of $SO(n+2)$, restricted to a hypersphere or to a catenoidal neck.  What follows is a presentation of these function.

\bigskip \noindent \itshape  1. Jacobi fields of the  hyperspheres. \upshape \smallskip
	
The linearized mean curvature operator of $S_\alpha$ is easily computed to be 
$$\mathcal L_a := \sin^{-2}(\alpha) \big(  \Delta_{\Sph^n} + n \big) \, .$$
Therefore, the Jacobi fields of $\mathcal L_a$ are simply the eigenfunctions of the $n$-sphere with eigenvalue $n$.  In the context of Lemma \ref{lemma:linmeancurv}, these can be derived by considering all non-trivial rotations of $S_\alpha$, namely the rotations generated by the vector fields
$$V_k := x^k \frac{\partial}{\partial x^0} - x^0 \frac{\partial}{\partial x^k} \: \: \mbox{for $k= 1, \ldots, n+1$} \, .$$  Taking the inner product of $V_k$ with $N_\alpha$ and restricting the resulting function to $S_\alpha$ gives a Jacobi field for each $k$.  This procedure yields the coordinate functions $x^k$ restricted to $S_\alpha$.

\bigskip \noindent \itshape 2. Jacobi fields of the catenoidal necks. \upshape \smallskip

The catenoidal necks of $\approxsol$ are catenoids $\Sigma$ in $\R^{n+1}$ that have been
re-scaled and embedded in $\Sph^n$ by the inverse of the canonical stereographic projection.  It turns out that when the scale parameter $\eps$ is sufficiently small, it is sufficient to consider the Jacobi fields of $\Sigma$ treated as a hypersurface in $\R^{n+1}$, where the ambient metric is Euclidean (rather than the one induced by stereographic projection from the sphere) and the ambient isometries are the rigid motions of $\R^n$.  The linearized mean curvature operator of $\Sigma$ with respect to this background is easily computed to be
$$
\mathcal{L}_\Sigma  :=  \frac{1}{\phi^{n}} \, \frac{\partial}{\partial s} \left( \phi^{n-2}  \frac{\partial}{\partial s} \right) + \frac{1}{\phi^{2}}Ê\, \Delta_{\Sph^{n-1}} +  \frac{n(n-1)}{\phi^{2n}}
$$
in the standard parametrization of Definition \ref{defn:gencat}. The isometries generating the relevant Jacobi fields of $\Sigma$ are as follows.  First, the ambient space $\R^{n+1} = \R \times \R^{n}$ possesses $n$ translations along the $\R^{n}$ factor and one translation in the $\R$ direction, which are generated by the vector fields
	$$V^{\text{trans}}_k := \frac{\partial}{\partial y^k} \:\: \mbox{ for $k=1,\ldots, n+1$}  \, .$$ 
	Then there are $n$ rotations of $\R \times \R^{n} $ that do not preserve the $\R$-direction, which are generated by the vector fields
	$$ V^{\text{rot}}_{1k} := y^{1} \frac{\partial}{\partial y^k} - y^k \frac{\partial}{\partial y^{1}} \:\: \mbox{ for $k = 2, \ldots, n+1$} \, .$$
	Finally, the motion of dilation in $\R^{n+1}$, though not an isometry, does preserve the mean curvature zero condition and is thus a geometric motion to which Lemma \ref{lemma:linmeancurv} can be applied.  Dilation is generated by the vector field
	$$V^{\text{dil}} := \sum_{k=1}^{n+1} y^k \frac{\partial}{\partial y^k} \, . $$  
	
	The Jacobi fields of $\mathcal{L}_\Sigma$ arising from the three classes of motions above can be found by the procedure of Lemma \ref{lemma:linmeancurv}.   One obtains the following non-trivial functions:
\begin{equation}
	\mylabel{eqn:catjacobi}
	\begin{aligned}
		J_1 (s) &:=  \langle N_\Sigma, V^{\text{trans}}_1 \rangle = \frac{\dot \phi(s)}{\phi(s)} \\
		J_k(s, \Theta) &:= \langle N_\Sigma, V^{\text{trans}}_k \rangle = - \frac{\Theta^k}{\phi^{n-1}(s)} \qquad k= 2, \ldots, n+1 \\
		J_{1k} (s, \Theta) &:= \langle  N_\Sigma, V^{\text{rot}}_{1k} \rangle =  \Theta^k \! \left( \frac{\psi(s)}{\phi^{n-1}(s)} + \dot \phi(s) \right) \qquad k = 2, \ldots, n+1 \\
		J_0(s) &:= \langle N_\Sigma, V^{\text{dil}} \rangle = \frac{\psi(s) \dot \phi(s)}{\phi(s)} - \frac{1}{\phi^{n-2}(s)} \, .
	\end{aligned}
\end{equation}
Note that the functions $J_k$ with $k \neq 0$ have odd symmetry with respect to the central sphere of $\Sigma$, i.e.~with respect to the transformation $s \mapsto -s$; while $J_{1k}$ and $J_0$ have even symmetry.  Note also that $J_1$ is bounded while $J_0$ has linear growth when the dimension is $n=2$ and is bounded in higher dimensions; $J_k$ decays like $\exp(-(n-1)|s|)$ for large $|s|$; and $J_{1k}$ grows like $\exp(|s|)$ for large $|s|$.

\subsection{Function Spaces and Norms}
\mylabel{sec:funcspace}

It does not seem possible to obtain a `good' linear estimate of the form $\Vert \linop (u) \Vert \geq C \Vert u \Vert$ with any straightforward choice of Banach subspaces and norms, where `good' in this case means with a constant $C$ independent of $\tau$, due to the presence of Jacobi fields and approximate Jacobi fields.  In order to remedy this situation, one must implement two ideas.  As mentioned above, the first of these ideas is to impose symmetry conditions on the functions chosen for deforming $\approxsol$ which have the effect of ruling out the Jacobi fields.  The second idea is to introduce a special norm to measure the `size' of functions $f \in C^{2, \beta}(\approxsol)$ in order to properly determine the dependence on the parameter $\tau$ of the various estimates needed for the application of the inverse function theorem. 

The norm in question is a so-called \emph{weighted Schauder norm}.  To define this norm, one must first define a \emph{weight function} on $\approxsol$.   Let $\eps = \eps(\tau)$ be the scale parameter of $\approxsol$, recall the definitions from Section \ref{subsec:assembly} and fix some $\rho_0$ independent of $\tau$ satisfying $\rho_o \gg 2 \rho_\eps$ such that the balls of radii $2 \rho_0$ centered on two different neck regions do not intersect.

\begin{defn}  
	\mylabel{defn:weight}
	The \emph{weight function} $\zeta_\eps : \approxsol \rightarrow \R$ is defined by
\begin{equation*}
	\zeta_\eps (q) = 
	\begin{cases}
		\eps \cosh(s) &\quad q = R_{\alpha + \tau/2}^{2k+1} \circ K^{-1} ( \eps \psi(s), \eps \phi(s) \Theta) \in \mathcal N_\eps^k \\
		\mathit{Interpolation} &\quad q \in \mathcal T_\eps \\
		\mathrm{dist} (q, \gamma) &\quad q \in \mathcal E_\eps \cap \left[ \bigcup_{k=0}^{N-1} B_{\rho_0} \big( R_{\alpha + \tau/2}^{2k+1}(p) \big) \right]\\[1ex]
		\mathit{Interpolation} &\quad q \in \mathcal E_\eps \cap \left[ \bigcup_{k=0}^{N-1} \left( B_{2\rho_0} \big( R_{\alpha + \tau/2}^{2k+1}(p) \big) \setminus B_{\rho_0} \big( R_{\alpha + \tau/2}^{2k+1}(p) \big) \right)   \right] \\[1ex]
		\rho_0 &\quad q \in \mathcal E_\eps \setminus \left[ \bigcup_{k=0}^{N-1} B_{2\rho_0} \big( R_{\alpha + \tau/2}^{2k+1}(p) \big) \right] \, .
	\end{cases}
\end{equation*}	
\end{defn}

The interpolation is such that $\zeta_\eps$ is smooth and monotone, and is such that $\zeta_\eps$ is invariant under the symmetries of $\approxsol$.  The weighted Schauder norm is then defined as follows. (The weighted norm below is essentially an adaptation of the weighted norm used in \cite{pacardriviere}.)  Let $T$ be any tensor on $\approxsol$ and let $\mathcal U \subseteq \approxsol$ be any open subset.  Recall the notation 
\begin{equation*}
	| T |_{0, \mathcal U} = \sup_{q \in \,  \mathcal U} \Vert T(q) \Vert \\
	\qquad \mbox{and} \qquad [T]_{\beta, \, \mathcal U} = \sup_{q, q' \in \, \mathcal U} \frac{\Vert T(q') - \Xi_{q,q'} (T(q)) \Vert}{\mathrm{dist}(q,q')^\beta} \, ,
\end{equation*}
where the norms and the distance function that appear are taken with respect to the induced metric of $\approxsol$, while $\Xi_{q,q'}$ is the corresponding parallel transport operator from $q$ to $q'$.  Now let $\mathit{Tub}_\rho(\gamma)$ be the tubular neighbourhood of $\gamma$ having width $\rho$; and for any $\arctan(\eps/2) < \rho < \rho_0$ define the annular region $A_\rho = \approxsol \cap \big[\mathit{Tub}_\rho(\gamma) \setminus \mathit{Tub}_{\rho/2}(\gamma) \big]$. Then the norm on any subset $\mathcal U \subseteq A_\rho$ is
\begin{equation*}
	| f |_{l, \beta, \delta, \, \mathcal U \cap A_\rho} := \zeta_\eps^{-\delta}(\rho) \vert f \vert_{0, \, \mathcal U \cap A_\rho} + \ \cdots + \zeta_\eps^{-\delta + l}(\rho) | \nabla^l f |_{0, \, \mathcal U \cap A_\rho} + \zeta_\eps^{-\delta + l + \beta}(\rho) [ \nabla^l f ]_{\beta, \, \mathcal U \cap A_\rho} \, .
\end{equation*}
Now make the following definition.

\begin{defn} 
	Let $\mathcal U \subset \approxsol$. The $C^{l, \beta}_\delta$ norm of a function defined on $\mathcal U$ is given by
\begin{equation}
	\mylabel{eqn:weightnorm}
	| f |_{C^{l, \beta}_\delta (\mathcal U)} := \sum_{i=0}^l | \nabla^i f |_{0, \,\mathcal U \cap \left[ \approxsol \setminus N_{\rho_0}(\gamma) \right]} + [ \nabla^l f ]_{\beta, \, \mathcal U \cap \left[ \approxsol \setminus N_{\rho_0} \right]} + \sup_{\rho \in (2 \eps,\rho_0]} | f |_{l, \beta, \delta, \, \mathcal U \cap A_\rho} \, . 
\end{equation}  
\end{defn}

\noindent The notation for the norm $| \cdot |_{C^{l, \beta}_\delta(\approxsol)}$ will be abbreviated $| \cdot |_{C^{l, \beta}_\delta}$ but $C^{l,\beta}_\delta$ norms over other subsets will be written out in full.  Now let  $X$ be any space of tensor fields.  The Banach space $\clbg (X)$ denotes the $C^{l,\beta}$ tensor fields in $X$ measured with respect to the norm \eqref{eqn:weightnorm}.  Finally, it is well known that all the `usual' properties that one would expect from a Schauder norm (multiplicative properties, interpolation inequalities, and so on) are satisfied by the weighted $\clbg$ norms. 

The choice of the Banach space in which a solution of the deformation problem will be found is governed by the following considerations.  The functions must be $C^{2,\beta}$ and they must be measured with respect to the $\ctbg$ norm.  Also, they must be invariant with respect to all the symmetries of $\approxsol$ since this will have the effect of ruling out the presence of non-invariant Jacobi fields.  The following space meets all these needs. 

\begin{defn}
	Let $X :=  \{ f : \approxsol \rightarrow \R \: : \: f \circ R_{2\alpha + \tau} = f \circ T = f \: \mbox{and} \: f \circ S^{\, 01}_{B} = f  \: \:  \forall \: B \in SO(n) \}$.
\end{defn}

\noindent Since $R_{2 \alpha + \tau}$, $T$ and all $S^{\, 01}_{B}$ are isometries of $\Sph^n$, then $\difop (f) \circ R_{2 \alpha + \tau} = \difop ( f \circ R_{2 \alpha + \tau}) = \difop (f)$ and similarly $\difop (f) \circ S^{\, 01}_{B} =\difop(f) \circ T = \difop (f)$ whenever $f \in X$.  Hence $\difop$ can by symmetrized to yield a new operator (which will be given the same name) $\difop : X \rightarrow X$.  In addition, it is easy to verify that  $\linop : C^{2,\beta}_\delta (X) \rightarrow C^{0, \beta}_{\delta - 2}(X)$ is a bounded operator whose operator norm is bounded above by a constant independent of $\tau$.

\subsection{The Linear Estimate}
\mylabel{sec:linest}

The most involved estimate necessary for invoking the inverse function theorem in the proof of Main Theorem \ref{result1} is the \emph{linear estimate}, i.e.~the estimate \eqref{eqn:iftestone} of the linearization $\linop$ from below.  The purpose of this section is to prove this estimate in the Banach space $\ctbg(X)$.  Because the second-order part of $\linop$, namely the Laplace operator, has a double indicial root when the dimension of $\approxsol$ is $n = 2$, the proof of the estimate is more complicated in this dimension.  Thus the easier case of dimension $n\geq 3$ will be given first.   The method will be to construct an explicit solution of the equation $\linop(u) = f$ by patching together local solutions on the neck region and away from the neck region.   

\begin{prop}
		\mylabel{thm:linesthigh}
	Suppose $n\geq 3$ and choose $\delta \in (2-n,0)$ and $\tau >0$ sufficiently small. The linearized operator $\linop : \ctbg (X) \rightarrow \cobg (X)$ satisfies the estimate
	$$| \linop (u) |_{\cobg} \geq C | u |_{\ctbg}$$
where $C$ is a constant independent of $\tau$.
\end{prop}

\begin{proof}
The patching argument requires two carefully defined partitions of unity for the various pieces of $\approxsol$.  First, for any $\rho \in (\tilde \rho_\eps, \rho_0)$ define $\mathcal E^k_\eps(\rho) := \mathcal E_\eps^k \setminus \mathit{Tub}_{\rho}(\gamma)$, and define $\mathcal N^k_\eps(\rho)$ to be the disjoint component of $\mathit{Tub}_\rho(\gamma)$ containing $\mathcal N^k_\eps$.  Next, define the smooth, monotone cut-off functions
\begin{align*}
	\chi^k_{\mathit{ext}, \rho} (q) &:= 
	\begin{cases}
		1 &\qquad q \in \mathcal E_\eps^k(2\rho) \\
		0 &\qquad q \in \approxsol \cap \mathit{Tub}_\rho(\gamma)
	\end{cases} \\[1ex]
	\chi^k_{\mathit{neck}, \rho} (q) &:= 
	\begin{cases}
		1 &\qquad q \in \mathcal N_\eps^k(2\rho) \\
		0 &\qquad q \in \approxsol \setminus \mathit{Tub}_\rho(\gamma)
	\end{cases}	
\end{align*}
so that $\sum_{k=0}^{N-1} \chi^k_{\mathit{ext}, \rho} + \sum_{k=0}^{N-1} \chi^k_{\mathit{neck}, \rho} = 1$.  Finally, define another set of cut-off functions 
\begin{align*}
	\eta^k_{\mathit{ext}} (q) &:= 
	\begin{cases}
		1 &\qquad q \in \mathcal E_\eps^k \\
		0 &\qquad q \in \approxsol \setminus \big[ \mathcal T_\eps^{k-1, +} \cup \mathcal E_\eps^k \cup \mathcal T_\eps^{k,-} \big]
	\end{cases} \\[1ex]
	\eta^k_{\mathit{neck}} (q) &:= 
	\begin{cases}
		1 &\qquad q \in \mathcal N_\eps^k \\
		0 &\qquad q \in \approxsol \setminus \big[ \mathcal T_\eps^{k, -} \cup \mathcal N_\eps^k \cup \mathcal T_\eps^{k,+} \big]
	\end{cases}	
\end{align*}
so that once again $\sum_{k=0}^{N-1} \eta^k_{\mathit{ext}} + \sum_{k=0}^{N-1} \eta^k_{\mathit{neck}} = 1$.  In addition, one can assume that these cut-off functions are invariant under the all the desired symmetries.
	
Suppose now that $f \in \cobg(X)$ is given.  The solution of the equation $\linop (u) = f$ will be constructed in three stages: local solutions on the neck regions will be found; then local solutions on the exterior regions will be found; and finally these solutions will be patched together to form an approximate solution which can be perturbed to a solution by iteration.   To begin this process, fix a small $\rho \in (\tilde \rho_\eps, \rho_0)$ and write $f = \sum_{k=0}^{N-1} f_{\mathit{ext}}^k +  \sum_{k=0}^{N-1} f_{\mathit{neck}}^k$ where $f_{\mathit{ext}}^k := f \cdot \chi_{\mathit{ext}, \rho}^k$ and $f_{\mathit{neck}}^k := f \cdot \chi_{\mathit{neck}, \rho}^k$.  The various symmetries satisfied by $f$ ensure that $f^k_{\mathit{ext}, \rho} := f_{\mathit{ext}, \rho}$ and $f^k_{\mathit{neck}, \rho} := f_{\mathit{neck}, \rho}$ for all $k$.  Also, $f_{\ast, \rho} \circ S^{01}_B = f_{\ast, \rho}$ for all $B \in O(n)$ and $f_{\ast, \rho} \circ T = f_{\ast, \rho}$.

\paragraph{Step 1.  Local solutions on the neck regions.}
For each $k$, the subset $K \circ R_{\alpha + \tau/2}^{-(2k+1)}\big(\mathcal N_\eps^k(\rho) \big) \subseteq \R \times \R^n$ is the union of two graphs over the annulus $\{ \hat y : \eps \leq \| \hat y \| \leq \rho' \} \subseteq \R^n$ for some $\rho'$ which is small but independent of $\eps$.  The graphing functions are $y^1 = \pm \tilde F_\eps (\| \hat y \|)$ as defined in Section \ref{subsec:assembly}.  Furthermore, $K \circ R_{\alpha + \tau/2}^{-(2k+1)}\big(\mathcal N_\eps^k(\rho) \big)$ is a perturbation of the $\eps$-scaled catenoid $\eps \Sigma$.  Consequently, the function $f_{\mathit{neck}}$ and the equation $\linop (u) = f_{\mathit{neck}}$ can be pulled back to a compact subset of $\eps \Sigma$.  In this formulation, this compact subset carries a perturbation of the catenoid metric $4 \eps^2 g_\Sigma$.  However,  $f_{\mathit{neck}}$ will be viewed as a function of compact support on the complete catenoid $\eps \Sigma$ carrying exactly the metric $4\eps^2 g_\Sigma$, and the equation that will be solved in this step is $\frac{1}{4}\mathcal L_{\eps \Sigma} (u) = f_{\mathit{neck}}$ where  $\frac{1}{4}\mathcal L_{\eps \Sigma}$ is the linearized mean curvature operator of $\eps \Sigma$ with this metric.   In addition, $f_{\mathit{neck}, \rho}$ is invariant under the symetries of $\eps \Sigma$ induced by the transformations $(y^1 , \hat y) \longmapsto (- y^1, \hat y )$ and $(y^1, \hat y) \longmapsto (y^1, B \hat y)$ for all $B \in O(n)$.

Parametrize the catenoid by $(s, \Theta)$ as in Definition \ref{defn:gencat} and let $| \cdot |_{C^{k, \alpha}_\delta (\eps_a K)}$ denote the standard weighted $C^{l, \beta}_\delta$ norm on $\eps \Sigma$, so that 
$$
| u |_{{C}^{l, \beta}_\delta (\eps \Sigma)} : = |(\eps \cosh(s))^{-\delta} u|_{0, \eps \Sigma} + \cdots + | (\eps  \cosh(s))^{-\delta+ l } \nabla^l u |_{0, \eps \Sigma} + [ (\eps  \cosh(s))^{-\delta+ l + \beta}\nabla^l u ]_{\beta, \eps \Sigma}
$$ 
where the norms and derivatives correspond to the metric on $\eps \Sigma$.  In this parametrization, the symmetries induced on functions of $\eps \Sigma$ are $u(s, \Theta) = u(-s, \Theta)$ and $u(s, \Theta) = u(s, B( \Theta))$ for all $B \in O(n)$ (which just says that $u$ is independent of $\Theta$).  The corresponding spaces of functions invariant under these symmetries will be denoted by ${C}^{l, \beta}_{\delta , \mathit{sym}} (\eps \Sigma)$. 

A standard separation of variables argument shows that when $\delta <0$, the kernel of the operator $\frac{1}{4}\mathcal L_{\eps \Sigma} : C^{2, \beta}_\delta (\eps \Sigma) \rightarrow C^{0, \beta}_{\delta-2} (\eps \Sigma)$ consists of the linear span of the Jacobi fields $\{ J_{1k} : k = 2, \ldots, n+1\}$.  However, none of these Jacobi fields is invariant under the symmetries $(s, \Theta) \longmapsto (s, B (\Theta))$ for $B \in O(n)$.  Hence the operator $\frac{1}{4}\mathcal L_{\eps \Sigma} : C^{2, \beta}_{\delta, \mathit{sym}} (\eps \Sigma) \rightarrow C^{0, \beta}_{\delta-2, \mathit{sym}} (\eps \Sigma)$ is injective for $\delta < 0$; and by duality it is surjective for $\delta > 2-n$.  Let $u_{\mathit{neck}} \in C^{2, \beta}(\eps \Sigma)$ be the unique solution of the equation $\frac{1}{4}\mathcal L_{\eps \Sigma}(u_{\mathit{neck}}) = f_{\mathit{neck}}$.  The estimate $|u_{\mathit{neck}} |_{\ctbg(\eps \Sigma)} \leq C |f_{\mathit{neck}}|_{\cobg(\eps \Sigma)} \leq C |f|_{\cobg}$ is valid, where $C$ is a constant independent of $\eps$.  With slight abuse of notation, extend this function to all of $\approxsol$ by defining $\bar u_{\mathit{neck}} := \sum_{k=0}^{N-1} \chi_{\mathit{neck}, \rho}^k \!\cdot u_{\mathit{neck}}$.  One has the estimate $|\bar u_{\mathit{neck}} |_{\ctbg} \leq C  |f|_{\cobg}$.

\paragraph{Step 2. Local solutions on the exterior regions.}  Given the local solution $\bar u_{\mathit{neck}}$ constructed in the previous step, choose a small $\kappa \in (0,1)$ and define
$$
\hat  f_{\mathit{ext}}^k  : =  \chi_{\mathit{ext}, \kappa \rho}^k  \big( f -  \linop( \bar u_{neck})  \big) \, .
$$
This function vanishes within an $\eps$-independent distance from the union of all the neck regions of $\approxsol$.  Therefore one can determine without difficulty $
|  \hat  f_{\mathit{ext}}^k   |_{C^{0, \beta}} \leq C_\kappa   |  f |_{\cobg}
$ for some constant $C_\kappa$ that depends on $\kappa$ and $\delta$.  Here, $| \cdot |_{C^{0,\beta}}$ is the un-weighted Schauder norm.

By symmetry, one can say that $\hat f_{\mathit{ext}}^k := \hat f_{\mathit{ext}}$ for each $k$. In addition, $\hat f_{\mathit{ext}}$ is invariant under the symmetries $T$ and $S^{01}_B$ for all $B \in O(n)$.  Denote the space of $C^{l,\beta}$ functions on $S_\alpha$ possessing these symmetries by $C^{l,\beta}_{\mathit{sym}}(S_\alpha)$.  The function $\hat f_{\mathit{ext}}$ can be viewed as a function of compact support on the perturbed hypersphere $\tilde S_\alpha$.  Since $\tilde S_\alpha$ is a normal graph over the hypersphere $S_\alpha$, this function can be pulled back to the hypersphere $S_\alpha$.  The metric  carried by $S_\alpha$ in this identification is a perturbation of the standard induced metric $\sin^2(\alpha) g_{\Sph^{n}}$.  However, the equation that will be solved here is $\mathcal L_{\alpha} (u) = \hat f_{\mathit{ext}}$ where $\mathcal L_{\alpha}$ is the linearized mean curvature operator of $S_\alpha$ when it carries the un-perturbed metric $\sin^2(\alpha) g_{\Sph^{n}}$. 

By another standard argument, one can show that the kernel of $\mathcal L_\alpha$ on $S_\alpha$ consists of the linear span of the coordinate functions $x^k$ restricted to $S_\alpha$.  None of these functions satisfies all the required symmetries; hence the operator $L_\alpha $ has no kernel when the symmetry conditions are imposed and thus $\mathcal L_\alpha : C^{2,\beta}_{\mathit{sym}}(S_\alpha) \rightarrow C^{0,\beta}_{\mathit{sym}}(S_\alpha)$ is bijective.  Let $u_{\mathit{ext}} \in C^{2,\beta}_{\mathit{sym}}(S_\alpha)$ denote the unique solution of the equation $\mathcal L_{\alpha} (u_{\mathit{ext}}) = \hat f_{\mathit{ext}}$.  One has the estimate  $|u_{\mathit{ext}} |_{C^{2,\beta}(S_\alpha)} \leq C |\hat f_{\mathit{ext}}|_{C^{0,\beta}} \leq C_\kappa |f|_{\cobg}$.  

To proceed, recall the $(\mu, \Theta) \in (0, \pi) \times \Sph^{n-1}$ parametrization for $S_\alpha \setminus \{ p^+, p^- \}$.  By examining the Taylor expansion of $u_{\mathit{ext}}$ at the points $\mu _0 = 0$ or $\mu_0= \pi$ and invoking its symmetries, one finds that $u_{\mathit{ext}} = a + v_{\mathit{ext}}$ where $a := u_{\mathit{ext}}(\mu_0)$, while $| v_{\mathit{ext}} (\mu)|  \leq C_\kappa   |\mu - \mu_0 |^2 |f|_{\cobg}$.  One can now extend $u_{\mathit{ext}}$ to all of $\approxsol$ as follows.  Recall that  $J_0$ is the bounded Jacobi field on the neck having even symmetry, which can be normalized to have the asymptotic expansion $J_0(s) = 1 + \tilde J_0(s)$ where $\tilde J_0(s) = \mathcal O(\cosh^{2-n}(s))$ for large $|s|$.  Now define, again with slight abuse of notation, the function $$\bar u_{\mathit{ext}} := \sum_{k=0}^{N-1} \eta_{\mathit{ext}}^k \cdot u_{\mathit{ext}}  + \sum_{k=0}^{N-1} \eta_{\mathit{neck}}^k \cdot  a J_0 \, .$$ 

\paragraph*{Step 3. Estimates and convergence.}  In the preceding two steps, local solutions $\bar u_{\mathit{neck}}$ and $\bar u_{\mathit{ext}}$ on the neck regions and on the exterior regions, respectively, were found and extended to all of $\approxsol$.  Define $\bar u := \bar u_{\mathit{neck}} + \bar u_{\mathit{ext}}$.  Then 
\begin{equation}
	\mylabel{eqn:cvgce}
	\begin{aligned}
		\linop (\bar u) - f &= \sum_{k=0}^{N-1} \left[ \eta_{\mathit{ext}}^k (\linop - \mathcal L_\alpha) ( u_{\mathit{ext}} ) + \chi_{\mathit{neck}, \kappa \rho}^k  ( \linop - \tfrac{1}{4} \mathcal L_{\eps \Sigma} ) (u_{\mathit{neck}}) \right. \\
		&\qquad \quad \left. +  \eta_{\mathit{neck}}^k ( \linop - \tfrac{1}{4} \mathcal L_{\eps \Sigma}) ( a J_0) +  [ \linop , \eta_{\mathit{neck}}^k ] ( a \tilde J_0 ) + [ \linop , \eta_{\mathit{ext}}^k ] ( v_{\mathit{ext}} ) \right]
	\end{aligned}
\end{equation}
where $[\mathcal L, \eta] (u) :=\mathcal L (\eta u) - \eta \mathcal L(u)$.  Each term in \eqref{eqn:cvgce} will be shown to be small in the $\cobg$ norm.

Begin with the first term in \eqref{eqn:cvgce}.  Note that this term is supported in $\mathcal E_\eps$ meaning that it vanishes within  a distance of $\rho_\eps = \eps^{(3n-3)/(3n-2)}$ from the centers of the neck regions.  The fact that the operator $\linop$ differs only very slightly from the operator $\mathcal L_\alpha$ away from the neck regions is the key to deriving the estimate.  That is, because $\approxsol$ is the normal graph over $S_\alpha$ of the function  $\alpha + \eps^{n-1} G$ defined in Section \ref{subsec:assembly}, one can use Lemma \ref{lemma:normgrgeom} to find the form of $\linop$.  This yields a rather unwieldy expression; however, with some work one finds
\begin{align}
	\mylabel{eqn:difference}
	( \linop - \mathcal L_{\alpha})(u) &= (\tilde \Delta - \Delta_{\alpha} )(u) + ( \| \tilde B \|^2 - \|B_{\alpha} \|^2 ) u \notag \\[1ex]
	&= \frac{\Delta u + \mathcal O ( \eps^{2n-2} \| \nabla^2 G \| \| \nabla G \|)  \cdot  \nabla u  + \mathcal O(\eps^{2n-2} \| \nabla G \|^2) \cdot  \nabla^2 u }{ \sin^2( \alpha + \eps^{n-1} G)} - \frac{\Delta u}{\sin^2(\alpha)} \notag\\
	&\qquad  + \left( \frac{n \sin^4(\alpha + \eps^{n-1}G) \cos^2(\alpha + \eps^{n-1}G) + \mathcal O(\eps^{2n-2} \| \nabla^2 G \|^2)}{( \sin^2(\alpha + \eps^{n-1} G) + \eps^{2n-2} \| \nabla G\|^2)^3} - n \cot^2(\alpha) \right) u \notag\\[1ex]
	&
	\begin{aligned}
		 &=\mathcal O(\eps^{2n-2} \| \nabla ^2 G \|  \| \nabla G \|) \cdot \nabla u + \big( \mathcal O( \eps^{n-1} G) + \mathcal O (\eps^{2n-2} \| \nabla G \|^2) \big) \cdot \nabla^2 u
	\\
	&\qquad + \big( \mathcal O ( \eps^{n-1} G)  + \mathcal O( \eps^{2n-2} \| \nabla ^2 G \|^2) \big) u
	\end{aligned}
\end{align}
evaluated on an arbitrary function $u \in \ctbg(\approxsol)$, where $\tilde \Delta$ is the Laplacian of $\approxsol$ and $\tilde B$ is its mean curvature.  The notation $\mathcal O(\ast) \cdot \nabla^k u$ means a linear quantity in $\nabla^ku$ with coefficients satisfying $\mathcal O(\ast)$.  The reason one reaches \eqref{eqn:difference} is because the largest terms of $\tilde \Delta - \Delta_{\alpha}$ and $\| \tilde B \|^2 - \|B_{\alpha} \|^2$, as computed from Lemma \ref{lemma:normgrgeom}, have been retained at the expense of all the smaller terms.  By recalling the expansion of $G$ given in equation \eqref{eqn:normgraphasym}, one can now compute
$$| ( \linop - \mathcal L_{\alpha})(u_{\mathit{ext}}) |_{C^0_{\delta - 2}(\mathcal E_\eps)} \leq  C \eps^{(4n-3)/(3n-2)} |u _{\mathit{ext}} |_{\ctbg(\mathcal E_\eps)} \leq C \eps^{(4n-3)/(3n-2)} | f |_{\cobg}  \, .$$
An estimate for the H\"older coefficient can be found by similar computations.

Next, using the fact that the operator  $\linop$ differs only very slightly from the operator  $\mathcal L_{\eps \Sigma}$ near the neck regions, one can estimate the second term in \eqref{eqn:cvgce}.  Begin with the computation in $\mathcal N_\eps^k$.   Here, one employs Lemma \ref{lemma:stereomeancurv} and the parametrization for the catenoid used in Definition \ref{defn:gencat} to calculate
\begin{align}
	\mylabel{eqn:differencetwo}
	\linop(u) &= A^2 \mathcal L_{\eps \Sigma} (u) + (2-n) A \left(\frac{\psi \dot \psi + \phi \dot \phi}{\phi^2} \right) \frac{\partial u}{\partial s} + n \left( 1 + \eps^2 \left( \frac{\dot \phi \dot \psi}{\phi} - \frac{1}{\phi^{n-2}} \right) \right) u
\end{align}
where $A = \frac{1}{2} \big( 1 + \eps^2 ( \psi^2 + \phi^2) \big)$.  Therefore
\begin{align*}
	| ( \linop  - \tfrac{1}{4} \mathcal L_{\eps \Sigma} ) (u_{\mathit{neck}}) |_{C^0_{\delta-2}( \mathcal N_\eps^k ) }
	&= \mathcal O(\eps^2 \cosh^2(s)) \nabla^2 u_{\mathit{neck}} + \mathcal O (\eps \cosh(s)) \nabla u_{\mathit{neck}} + \mathcal O(1) u_{\mathit{neck}} \\
 	&\leq C \eps^{2(3n-3)/(3n-2)} |u_{\mathit{neck}}|_{\ctbg(\mathcal N_\eps^k)} \\
	&\leq C \eps^{2(3n-3)/(3n-2)} |f|_{\cobg}   \, .
\end{align*}
Now consider the same quantity further away from the neck region, in the region $\mathcal N_\eps^k( 2 \kappa \rho)$.  Similar calculations can be performed, this time computing $\linop$ from Lemma \ref{lemma:stereomeancurv} with a parametrization of a small perturbation of the catenoid.  The perturbation introduces additional error terms, but the largest error still comes from the $\mathcal O(1)$ terms.  Hence one obtains the estimate
$$| ( \linop  - \tfrac{1}{4} \mathcal L_{\eps \Sigma} ) (u_{\mathit{neck}}) |_{C^0_{2-\delta} (\mathcal N^k_\eps ( 2 \kappa \rho) \setminus \mathcal N_\eps^k)} \leq C \kappa^2 |u_{\mathit{neck}}|_{\ctbg(\mathcal N^k_\eps ( 2 \kappa \rho) \setminus \mathcal N_\eps^k)} \leq C \kappa^2 |f|_{\cobg}   \, .$$
Analogous H\"older coefficient estimates can be found in the same way.  In the end, one finds
\begin{align}
	\mylabel{eqn:finalneckest}
	|\chi_{\mathit{neck}, \kappa \rho}^k ( \linop  - \tfrac{1}{4} \mathcal L_{\eps \Sigma} )(u_{\mathit{neck}}) |_{\cobg} &\leq C \big( \kappa^2 + \eps^{2(3n-3)/(3n-2)} \big) |f|_{\cobg} \, .
\end{align}

The remaining terms in \eqref{eqn:cvgce} are considerably easier to estimate.  The third term can be estimated in much the same way as the first.  In fact, formula \eqref{eqn:differencetwo} applied to the function $u = J_0$ and the estimate $| a | \leq C |u_{\mathit{ext}}|_{C^{2,\beta}} \leq C_\kappa |f|_{\cobg}$ gives 
$$| a (\linop - \tfrac{1}{4} \mathcal L_{\eps \Sigma})(\eta_{\mathit{neck}}^k J_0)|_{\ctbg} \leq C_\kappa \eps^{2(3n-2)/(3n-2)} | f |_{\cobg} \, .$$
Finally, the estimate of the last term follows from the fact that in the support of the $[ \linop , \eta]$ quantities, the functions $v_{\mathit{ext}} $ and $a \tilde J_0$ are of size  $\mathcal O ( [\mathrm{dist} ( \gamma, \cdot) ]^2 |f|_{\cobg} ) $ while the $\cobg$ norm of the coefficients of $[ \linop , \eta ]$ are $\mathcal O(1)$ by definition of the cut-off functions.  Hence an estimate proportional to $|f|_{\cobg}$ is again valid, where the proportionality constant is of size $\mathcal O(\eps^{2(3n-2)/(3n-3)})$.

The outcome of the above analysis is that one has a function $\bar u \in \ctbg(\approxsol)$ satisfying the estimates $|\bar u | _{\ctbg} \leq C_\kappa | f|_{\cobg}$ as well as 
$$| \linop (\bar u ) - f |_{\cobg} \leq C \big( \kappa^2 + \eps^{2(3n-3)/(3n-2)}  + \eps^{(4n-3)/(3n-2)} \big) |f|_{\cobg}$$  
where $C$ is bounded above independently of $\kappa$.   Since the constant in front of  $|f|_{\cobg}$ can be made as small as desired by appropriate choices of $\eps $ and $\kappa$, a straightforward iteration argument now proves the existence of an actual solution $u \in \ctbg(\approxsol)$ of the equation $\linop (u) = f$ satisfying the estimate $|u | _{\ctbg} \leq C | f|_{\cobg}$.
\end{proof}

A modification of the preceding proof, which takes into account the double indicial root of $\linop$, can be found to handle the two-dimensional case. 
\begin{prop}
	\mylabel{thm:linest}
	Suppose $n=2$ and choose $\delta \in (-1,0)$ and $\tau >0$ sufficiently small. The linearized operator $\linop : \ctbg (X) \rightarrow \cobg (X)$ satisfies the estimate
	$$| \linop (u) |_{\cobg} \geq C \eps^{-\delta} | u |_{\ctbg}$$	%
where $\eps = \eps(\tau)$ is the scale parameter of $\approxsol$ and $C$ is a constant independent of $\tau$.
\end{prop}

\begin{proof}
The proof of the result in dimension $n=2$ follows broadly the same plan as the proof in higher dimensions.  For given $f \in \cobg(\approxsol)$, a solution of the equation $\linop(u) = f$ must be found and estimated.  Begin by defining the different regions of $\approxsol$ and their corresponding cut-off functions exactly as before, along with the decomposition $f = \sum_{k=0}^{N-1} (f_{\mathit{neck}}^k + f_{\mathit{ext}}^k)$.

The significant difference between this proof and the preceding one is how the local solutions are found on the neck.  The reason is that due to the indicial roots, there is no range of $\delta$ in which the operator $\mathcal L_{\eps \Sigma} : C^{2,\beta}_{\delta, \mathit{sym}} (\eps \Sigma) \rightarrow  C^{0,\beta}_{\delta-2, \mathit{sym}}  (\eps \Sigma)$ is both injective and surjective.  One must thus proceed differently.

First, recall that in dimension $n=2$, the function  $J_0$ is the linearly growing Jacobi field on the neck having even symmetry, which can be normalized to have the asymptotic expansion $J_0(s) = -1 + |s| + \tilde J_0(s)$ where $\tilde J_0(s) = \mathcal O(\cosh^{-2}(s))$ for large $|s|$.  The bounded Jacobi field having odd symmetry is $J_1$ and can be normalized to have the asymptotic expansion $J_1(s) = 1 + \tilde J_1(s)$ where $\tilde J_1(s) = \mathcal O(\cosh^{-2}(s))$ for large positive $s$.  Define a smooth, odd function $\chi : \eps \Sigma \rightarrow \R$ with the property $\chi(s) = 1$ for $s \geq 1$ and $\chi(s) = -1$ for $s \leq -1$.  Now set $K_1 := \chi J_1$ and define the linear subspaces $\mathcal D := \mathrm{span}_{\R} \{ J_0, K_1 \}$.  Together, $J_0$ and $K_1$ span all possible asymptotic behaviours corresponding to the indicial root at zero of $\mathcal L_{\eps \Sigma}$.  The following two facts hold, and are well known in the theory of elliptic operators on asymptotically flat manifolds.  The operator $\mathcal L_{\eps \Sigma} : C^{2, \beta}_{\delta , \mathit{sym}} (\eps \Sigma) \oplus \mathcal D \rightarrow C^{0,\beta}_{\delta - 2, \mathit{sym}}$ is surjective in the range $\delta \in (-1,0)$ with one-dimensional kernel spanned by $J_0$.  Furthermore, there is a bounded right inverse mapping into $C^{2, \beta}_{\delta , \mathit{sym}} (\eps \Sigma) \oplus \mathcal D_0$ where $\mathcal D_0 := \mathrm{span}_{\R} \{ K_1 \}$.  The proof of these facts follows from a similar proof that can be found for example in \cite{kusnermazzeopollack} in the context of constant mean curvature surfaces.  Consequently there is a solution $u_{\mathit{neck}} \in C^{2, \beta}_{\delta , \mathit{sym}} (\eps \Sigma) \oplus \mathcal D_0$ of the equation $\mathcal L_{\eps \Sigma} (u_{\mathit{neck}}) = f_{\mathit{neck}})$ which can be decomposed as $u_{\mathit{neck}} = v_{\mathit{neck}} + a_1 K_1$ where $v_{\mathit{neck}} \in C^{2, \beta}_{\delta , \mathit{sym}} (\eps \Sigma)$ and $a_1 \in \R$ satisfy $$|v_{\mathit{neck}} |_{\ctbg(\eps \Sigma)} + \eps^{-\delta} |a_1| \leq C |f_{\mathit{neck}}|_{C^{0,\beta}_{\delta - 2}(\eps \Sigma)} \leq C |f|_{\cobg}$$
for some constant $C$ independent of $\eps$.
 	
	Recall now that $G$ is the singular solution of $\mathcal L_\alpha (G) = 0$ on $S_\alpha \setminus \{ p^+, p^-\}$, for which one has $G(s) = \gamma_0 + \gamma_1 |s| + \tilde G(s)$ where $\tilde G(s) = \mathcal O(\eps^2 \cosh^2(s))$ and $\gamma_0$ and $\gamma_1$ are non-zero constants.  Note that $\gamma_0 = O(| \log(\eps)|)$.  This follows from the expansion for $G$ and a formula for changing between the $\mu$-coordinate of $S_\alpha$ and the $s$-coordinate of $\eps \Sigma$ on the regions of overlap. The function $u_{\mathit{neck}}$ can thus be extended to $\approxsol$ by prescribing
	$$\bar u_{\mathit{neck}} := \sum_{k=0}^{N-1}\chi_{\mathit{neck}, \rho}^k v_{\mathit{neck}} +  \sum_{k=0}^{N-1} \left( \eta_{\mathit{neck}}^k \big( b_1 J_0 + a_1 K_1) + \eta_{\mathit{ext}}^k c_1 G \right)$$
	where the constants $b_1$ and $c_1$ are chosen to ensure matching in the constant and linear terms of asymptotic expansions of $u_{\mathit{neck}}$ and $G$, i.e.~$b_1 = - a_1 \gamma_1 / (\gamma_1 + \gamma_2)$ and $c_1 =-a_1 /(\gamma_0 +\gamma_1)$.  In the end, this solution satisfies the estimates $|\bar u_{\mathit{neck}}|_{\ctbg} \leq C \eps^\delta |f|_{\cobg}$ as well as
\begin{align*}
	|\chi_{\mathit{neck}, \kappa \rho}^k (\linop(\bar u_{\mathit{neck}}) - f_{\mathit{neck}}) |_{\cobg} &\leq C \big( \kappa^2 + \eps^{3/2 + \delta} |\log(\eps) | \big) |f|_{\cobg} \, .
\end{align*}
Note that $\delta > -3/2$ is required to ensure that the quantity in front of $|f|_{\cobg}$ small as $\eps \rightarrow 0$.

Now that one has the solution $\bar u_{\mathit{neck}}$ on the neck region and extended to all of $\approxsol$, the next step is to define $\hat f_{\mathit{ext}}$ in the same way as before and find the local solution on the exterior region of $\approxsol$.  The extension to all of $\approxsol$ is again complicated by the fact that one must match constant and linear terms in the Taylor series at the points $p^+$ and $p^-$.  This can be accomplished by using $G$ and $J_0$.  That is, define the extended function 
$$\bar u_{\mathit{ext}} := \sum_{k=0}^{N-1} \eta_{\mathit{ext}}^k \cdot (u_{\mathit{ext}} + c G) + \sum_{k=0}^{N-1} \eta_{\mathit{neck}}^k \cdot  bJ_0$$  where $b$ and $c$ are chosen to ensure matching in the constant and linear terms of the Taylor expansion; i.e.~$b = - a \gamma_1 / (\gamma_1 + \gamma_2)$ and $c =-a /(\gamma_0 +\gamma_1)$, where $a = u_{\mathit{ext}}(0)$.

The remainder of the analysis is the same as before, and leads to a good approximate solution of the equation $\linop (u) = f$.  The solution procedure can be iterated to yield an exact solution.  The result is a solution satisfying the bound $|u|_{\ctbg} \leq C \eps^\delta |f|_{\cobg}$.	
\end{proof}

\subsection{The Non-Linear Estimates and the Conclusion of the Proof}
\mylabel{subsec:nonlinest}

\paragraph*{The non-linear estimates.}  According to the discussion in Section \ref{sec:jacobi}, the proof of Main Theorem \ref{result1} requires two more estimates in addition to the linear estimate from the previous section: it is necessary to show that $\difop(0) - H_\alpha$ is small in the $\cobg$ norm; and it is necessary to show that $\Dif \difop(f)  - \linop$ can be made to have small $\ctbg$-operator norm if $f$ is chosen sufficiently small in the $\ctbg$ norm.  These two estimates are the subject of this section.  Once these estimates are given, it will be possible to conclude the proof of Main Theorem \ref{result1} by invoking the inverse function theorem.   Begin with the first non-linear estimate.

\begin{prop}
	\mylabel{prop:error}
	The quantity $\difop(0)$, which is the mean curvature of $\approxsol$, satisfies the following estimate.  If $\tau > 0$ is sufficiently small, then there exists a constant $C$ independent of $\tau$ so that 
	\begin{equation}
		\mylabel{eqn:error}
		\vert \difop(0) - H_\alpha \vert_{\cobg} \leq C \eps^{(2 - \delta)(3n-3)/(3n-2)}
	\end{equation}
	where $\eps =  \eps(\tau)$ is the scale parameter of $\approxsol$ as in Section \ref{subsec:assembly}.
\end{prop}

\begin{proof} 

The estimate \eqref{eqn:error} is the consequence of three separate calculations: in $\mathcal E_\eps$ where $|\zeta_\eps| \geq 2\rho_\eps$, in $\mathcal T_\eps$ where $|\zeta_\eps| = \mathcal O(\rho_\eps)$, and in $\mathcal N_\eps$ where $|\zeta_\eps| \leq \rho_\eps /2$.  For the first of these, one uses the expression for the mean curvature of a normal graph over the  hypersphere $S_\alpha$ from Lemma \ref{lemma:normgrgeom} for the  graphing function $F = \eps^{n-1} G$.  One obtains the expression
\begin{equation}
	\begin{aligned}
		\mylabel{eqn:normgrmeancurv}
			H \big( \exp(\eps^{n-1} G N_\alpha)(S_\alpha) \big) &=  \frac{-\eps^{n-1}\Delta G +  n \sin(\alpha + \eps^{n-1} G) \cos(\alpha+\eps^{n-1} G)}{A \sin(\alpha + \eps^{n-1} G)}\\
			&\qquad  - \frac{\eps^{2n-2} \big[ \eps^{n-1} \ddot G  + \cos(\alpha + \eps^{n-1} G) \sin(\alpha+\eps^{n-1} G) \big] \dot G^2}{A^3 \sin(\alpha + \eps^{n-1} G)} 
	\end{aligned}
\end{equation}
where $\Delta$ is the Laplacian of the standard metric of $\Sph^n$ and $A = \big( \sin^2(\alpha + \eps^{n-1} G) + \eps^{2n-2} \dot G^2  \big)^{1/2}$.  Here $G$ is the singular solution of $(\Delta + n)(G) = 0$ and $\eps$ is very small.  Thus expanding \eqref{eqn:normgrmeancurv} yields
\begin{align*}
	H \big( \exp(\eps^{n-1} G)(S_\alpha) \big) &= n \cot (\alpha) - \eps^{2n-2} \big( (n-2) \dot G^2 + 2 n G^2 \big) \cot(\alpha)\csc^2(\alpha) \\
	&\qquad + \eps^{3n-3}\mathcal O \big( |\ddot G | \, \dot  G^2 + |G| \, \dot G^2 + |G|^3 \big) \, .
\end{align*}
Since $H_\alpha = n \cot (\alpha)$ and using the expansion of $G$ given in \eqref{eqn:normgraphasym} as well as  the fact that  $\mathrm{dist}(\gamma, q)  \geq \rho_\eps$ in $\mathcal E_\eps$,  then one finds 
$$| H(\approxsol) - H_\alpha|_{C^0_{2-\delta} (\mathcal E_\eps )} \leq C \eps^{(2 - \delta)(3n-3)/(3n-2) } \, .$$
This completes the supremum estimate needed for the $\cobg$ estimate in $\mathcal E_\eps$.  The H\"older coefficient estimate needed for the $\cobg$ estimate can be found in a similar (though more involved) way.

The remaining two calculations use Lemma \ref{lemma:stereomeancurv}, namely the formula for the mean curvature of an embedded hypersurface in $\R^{n+1}$ with respect to the metric induced from $\Sph^{n+1}$ by stereographic projection.  Consider first the neck region $\mathcal N_\eps$ whose image under stereographic projection is a catenoid with vanishing mean curvature with respect to the Euclidean metric.  Therefore by Lemma \ref{lemma:stereomeancurv}, the mean curvature at a point $(\eps \psi(s), \eps \phi(s) \Theta) \in \eps \Sigma \cap \{ (y^1, \hat y) : \| \hat y \| \leq \rho_\eps \}$ is
$$H(\approxsol)(s) = n \eps \frac{\dot \phi(s) \psi(s) - \dot \psi(s) \phi(s)}{(\dot \psi(s)^2 + \dot \phi(s)^2)^{1/2}} = 
\begin{cases}
	\mathcal O(\eps | \log(\eps)|) & \qquad n=2 \\
	\mathcal O(\eps) & \qquad n \geq 3 \, .
\end{cases}$$
This is because $\phi(s) \in [\eps, \rho_\eps]$, while $\dot \phi(s) = \mathcal O (\phi(s))$ and $\psi(s) = c(\eps, n) + \mathcal O(\phi(s)^{2-n})$ for large $s$, where $c(\eps, 2) = \mathcal O(| \log(\eps)|)$ and $c(\eps, n) = \mathcal O(1)$ for $n \geq 2$.  Thus 
$$| H(\approxsol) - H_\alpha|_{C^0_{2-\delta}( \mathcal N_\eps )} \leq C\eps^{(2 - \delta)(3n-3)/(3n-2)} \, .$$
Once again, the H\"older coefficient is similar, and yields the same sort of bound. 
Finally, it remains to estimate the mean curvature in the region $\mathcal T_\eps$ whose image under stereographic projection is the graph of the function  $F := \tilde F_{a, \tau}$ defined in equation \eqref{eqn:mergedfn}.  Let $y = \| \hat y \|$. Then at the point $(F(y), \hat y) \in \tilde \Sigma_\eps \cap \{ (y^1, \hat y) : \rho_\eps/2 \leq y \leq 2\rho_\eps \}$, the formula of Lemma \ref{lemma:stereomeancurv} yields
$$H(\approxsol)(y) = \frac{1+ y^2 + F(y)^2}{2} \left( \frac{\ddot F(y)}{(1+ \dot F(y)^2)^{3/2}} + \frac{(n-1) \dot F(y)}{y(1+ \dot F(y)^2)^{1/2}} \right) + \frac{n(F(y) - y \dot F(y))}{(1+ \dot F(y)^2)^{1/2}} \, .$$
The function $F(y)$ has the form
$$F(y) = \eps c_n + C_n \eps^{n-1}L_{ n}(y) + \eta \big( \mathcal O(y^2) + \mathcal O( \eps^{3n-3} y^{4 - 3n}) \big)$$
where $L_{ n}(y) = \log(y)$ when $n=2$ and $L_{ n}(y) =  y^{2-n}$ when $n\geq 3$, while $C_n$ and $c_n$ are constants and $\eta$ is the cut-off function whose $C^{2,\beta}_0$ norm is bounded by an $\eps$-independent constant.  Using the fact that $L_n$ satisfies the zero mean curvature equation $\ddot L_n + \frac{n-1}{y}(1+ \dot L_n^2) \dot L_n= 0$ then one sees that
$$ | H(\approxsol) - H_\alpha|_{C^0_{2-\delta}( \mathcal T_\eps)} \leq C\eps^{(2 - \delta)(3n-3)/(3n-2)}$$
in the region $y \in [\rho_\eps/2, 2 \rho_\eps]$.   Once again, the H\"older coefficient estimate in this region is similar, and yields the same sort of bound. 
\end{proof}

The second non-linear estimate is essentially the same as \cite[Prop.~26]{mepacard1} and \cite[Prop.~28]{mepacard2}.  The result is re-stated here for reference in the proof of Main Theorem 1 that will follow below.

\begin{prop}
    \mylabel{prop:nonlin}
    The linearized mean curvature operator satisfies the following general estimate.  If $\eps>0$ is sufficiently small and  $f \in \ctbg(\approxsol)$ has sufficiently small $\ctbg$ norm, then there exists a constant $C$ independent of $\eps$ so that
    \begin{equation}
        \mylabel{eqn:nonlin}
        \big\vert \Dif \Phi_{\alpha, \tau} (f) (u) - \mathcal{L}_{\alpha, \tau} ( u ) \big\vert_{\cobg} \leq C \eps^{\delta - 1} |f|_{\ctbg}Ê\, \vert u \vert_{\ctbg }
    \end{equation}
     for any function $u \in \ctbg(\approxsol)$.
\end{prop}

\begin{proof}
This holds because of the natural scaling property of the mean curvature operator.  Pick a point  $q \in \approxsol$ and work in geodesic normal coordinates in a neighbourhood $\mathcal U$ of $q$.  Let $\mathcal U_f$ denote the normal deformation of $\mathcal U$ determined by $f$.  Now set $s := \zeta_\eps(q)$ and consider re-scaling by the factor $1/s$.  The mean curvature at $q$ satisfies $H_1(\mathcal U_f) = s^{-1} H_{1/s} ( s^{-1} \mathcal U_{f/s} )$ where $H_r$ is the mean curvature of a hypersurface of the sphere of radius $r$.  Hence
\begin{align}
	\label{eqn:nonlinproof}
	s^{2-\delta} | \Dif H_1( \mathcal U_{f})(u) - \Dif H_1( \mathcal U_0)(u) | &= s^{-1-\delta} | \Dif H_{1/s}( s^{-1} \mathcal U_{f/s})(u) - \Dif H_{1/s}( s^{-1}\mathcal U_0)(u) | \notag \\
	&\leq C s^{-1-\delta} |f|_{C^{2,\beta}(s^{-1} \mathcal U_0)} |u|_{C^{2,\beta}(s^{-1} \mathcal U_0)} \notag  \\
	&\leq C s^{\delta-1} |f|_{\ctbg(\mathcal U_0)} |u|_{\ctbg( \mathcal U_0)}
\end{align}
by performing the scaling backwards.  Here $C$ is a universal constant independent of $s$.  Allowing $q$ to range over all points in $\approxsol$ now yields the desired estimate; the constant appearing in \eqref{eqn:nonlin} arises because the largest $s^{-1+\delta}$ can be in the right hand side of \eqref{eqn:nonlinproof} is $\mathcal O (\eps^{\delta -1})$.  
\end{proof}

\paragraph*{The proof of Main Theorem \ref{result1}.}  The estimates for the proof of Main Theorem \ref{result1} are now all in place and the conclusion of the theorem becomes a simple verification of the conditions of the inverse function theorem.  By Proposition \ref{thm:linesthigh} and Proposition \ref{thm:linest} the linearization $\linop$ is injective on $\ctbg(X)$.  But $\linop - \Delta$ is a compact operator so $\linop$ has the same index as $\Delta$ on $\ctbg(X)$.  By self-adjointness, this index is zero, so that $\linop$ must be surjective as well.  One has the estimate
$$| \linop^{-1} (f) |_{\ctbg} \leq C_L | f |_{\cobg} \, ,$$
where $C_L = \mathcal O(\eps^{\delta})$ in dimension $n=2$ and $C_L = \mathcal O(1)$ in higher dimensions.  Now in order to achieve the bound
$$| \Dif \difop (f)(u) - \linop(u) |_{\cobg} \leq \frac{1}{2 C_L} |u|_{\ctbg} \, ,$$
for any $u \in \ctbg(\approxsol)$, it is necessary to have $|f|_{\ctbg} \leq \rho$ where $\rho = \mathcal O(\eps^{1-2\delta})$ in dimension $n=2$ or $\rho = \mathcal O(\eps^{1-\delta})$ in higher dimensions.  Therefore by the inverse function theorem of Section \ref{sec:jacobi} along with Proposition \ref{prop:nonlin}, a solution $f := f_{\alpha, \tau}$  of the deformation problem can be found if $\approxsol$ has been constructed so that $| \difop(0) - H_\alpha |_{\cobg} \leq \frac{\rho}{2 C_L }$.  Since Proposition \ref{prop:error} asserts that $$| \difop(0) - H_\alpha |_{\cobg} = \mathcal O(\eps^{(2-\delta)(3n-3)/(3n-2)})\, ,$$ this is true so long as $\eps$ is made sufficiently small by a small enough choice of $\tau$ and $\delta$ is chosen appropriately.

As a further consequence of these estimates, the Banach space inverse function theorem asserts that the solution of the equation $\difop(f_{a, \tau}) = 0$ satisfies the estimate 
$$| f_{\alpha, \tau} |_{\ctbg} = \mathcal O \big(  C_L \eps^{(2-\delta)(3n-3)/(3n-2)}  \big)$$
which is much smaller that $\eps$.  Therefore the size of the perturbation of $\approxsol$ created by the normal deformation of magnitude $f_{\alpha, \tau}$ is much smaller than the width of $\approxsol$ at its narrowest points, i.e.~in the neck regions where the width is $\mathcal O(\eps)$.  Thus $\approxsol$ remains embedded under this normal deformation.  This concludes the proof of Main Theorem \ref{result1}.   \hfill \qedsymbol

\section{Gluing Two Delaunay-Like CMC Hypersurfaces Together}

\subsection{Introduction}

The first main theorem asserts that if one positions hyperspheres of radius $\cos(\alpha)$ in a highly symmetrical manner, namely equally spaced around an equator, then one can connect them to each other by small catenoidal necks and perturb the resulting hypersurface into one of constant mean curvature $H_\alpha$.  These hypersurfaces should be considered the analogues of the Delaunay surfaces in $\Sph^3$.

The celebrated constructions of Kapouleas \cite{kapouleas7} and Mazzeo, Pacard and Pollack \cite{mazzeopacard1,mazzeopacardpollack} show that it is possible to connect Delaunay surfaces in $\R^3$ together in a wide variety of ways and thereby deduce considerable information about the moduli space of constant mean curvature surfaces.  It is therefore tempting to try to find a similarly robust procedure in $\Sph^{n+1}$ using the gluing methods of the Section \ref{sec:delaunay}.  Indeed, one can readily imagine an enormous multitude of less symmetrical configurations of  hyperspheres, positioned this way and that in $\Sph^{n+1}$, only satisfying the condition that they meet each other tangentially (allowing catenoids to be inserted near the points of contact).  The question becomes when the approximate solution constructed from such a configuration can be perturbed into a CMC hypersurface.

The purpose of the present section of this paper is to give a construction of a CMC hypersurface based on the configuration of  hyperspheres corresponding to the connected sum of exactly two of the Delaunay-like hypersurfaces constructed in the previous section.   The initial configuration of hyperspheres in this case will be only partially symmetric and thus the analysis that needs to be carried out will become much more complicated.  One reason for this is the extra bookkeeping required to keep track of the various pieces of the initial configuration.  However, a more significant reason is the presence of non-trivial approximate Jacobi fields of the linearized operator.  A new technique, based on the Kapouleas \cite{kapouleas7} and Korevaar-Kusner-Solomon \cite{kks} balancing condition, will be introduced to deal with this difficulty.

\begin{rmk}
	The Delaunay-like CMC hypersurfaces of the previous section can be either immersed and non-compact, or immersed and compact, or embedded and compact, depending on the choices made during the construction.  Such freedom is also available here, though in this section of the paper only compact CMC hypersurfaces will be found.
\end{rmk}

\subsection{The New Approximate Solution}
\mylabel{subsec:newassembly}

The approximately CMC hypersurface that looks like the fusion of two Delaunay-like hypersurfaces in $\Sph^{n+1}$ will be constructed roughly as follows.  Consider \emph{two} orthogonally intersecting geodesics in $\Sph^{n+1}$.  Introduce the small spacing parameter $\tau$ so that $0 < \tau \ll 2\pi$ and $2 \alpha + \tau = 2 \pi m / N$ for a pair of positive and relatively prime integers $m$ and $N$, where $N$ is an even integer such that $N/2$ is odd.  Now position rotated copies of $S_\alpha$ in a symmetric manner around both geodesics so that each rotated hypersphere is separated from its two nearest neighbours by a distance $\tau$.  

Since this initial configuration of  hyperspheres is less symmetric than before, one expects there to be non-trivial approximate Jacobi fields on the compact hypersurface formed by joining the hyperspheres together using normal perturbations and small catenoidal necks as before.  Overcoming these Jacobi fields will require extra flexibility in the initial configuration.  This will be incorporated here by introducing small \emph{displacement parameters} for each hypersphere.  These will be denoted by $\sigma_k$ where $k = 0, \ldots, N-1$ and often abbreviated simply by $\sigma$.  The idea is to rotate the $k^{\mathit{th}}$  hypersphere on each of the geodesics by an angle of $k(2\alpha +  \tau) + \sigma_k$ which displaces it slightly from its `natural' position.  The details of the construction outlined above will now be given.

\paragraph*{Positioning the hyperspheres.} For $s=1,2$ let $\gamma_s$ denote the geodesic of $\Sph^{n+1}$ formed by intersecting $\Sph^{n+1}$ with the $(x^0, x^s)$-plane.  Define the following rotations
$$ R_{1,\theta} = \left(
\begin{array}{ccc|c}
	\cos(\theta) & \sin(\theta) & 0 & 0 \\
	-\sin(\theta) & \cos(\theta) & 0 & 0 \\
	0 & 0 & 1 & 0 \\
	\hline
	0 & 0 & 0 & I_{n-1}
\end{array}
\right) \qquad \mbox{and} \qquad R_{2,\theta} =  \left(
\begin{array}{ccc|c}
	\cos(\theta) & 0 & \sin(\theta)  & 0 \\
	0 & 1 & 0 & 0 \\
	-\sin(\theta) & 0 & \cos(\theta)  & 0 \\
	\hline
	0 & 0 & 0 & I_{n-1}
\end{array}
\right)$$
generating motions around $\gamma_1$ and $\gamma_2$ respectively.  Here $I_{n-1}$ is the $(n-1) \times (n-1)$ identity matrix.   

Introduce the displacement parameters $\sigma_0, \ldots, \sigma_{N-1}$ and consider the set of all hyperspheres of the form $R_{s, \sigma_k} \circ R_{s, 2 \alpha + \tau}^k (S_\alpha)$ for $s=1,2$ and $k = 0, \ldots, N-1$.  This is a configuration of $2N-2$  hyperspheres wrapping $m$ times around $\gamma_1$ and another $m$ times around $\gamma_2$.  The initial configuration of hyperspheres is then
\begin{equation}
	\mylabel{eqn:paraminitconfig}
	\Lambda_{\alpha, \tau, \sigma}^{\#} := \bigcup_{s=1}^2  \bigcup_{k=0}^{N-1} R_{s, \sigma_{k}} \circ R_{s, 2 \alpha + \tau}^k (S_\alpha) \, .
\end{equation}
The separation between $R_{s, \sigma_k} \circ R_{s, 2 \alpha + \tau}^k (S_\alpha)$ and its nearest neighbour $R_{s, \sigma_{k+1}} \circ R_{s, 2 \alpha + \tau}^{k+1} (S_\alpha)$ along the geodesic  $\gamma_s$ is $\tau_k := \tau + \sigma_{k+1} - \sigma_{k}$.   Note that precisely two pairs of hyperspheres amongst these coincide exactly: namely $S_\alpha = R_{1, 2 \alpha + \tau}^{0}(S_\alpha) = R_{2, 2 \alpha + \tau}^{0} (S_\alpha)$ as well as $- S_\alpha = R_{1, 2 \alpha + \tau}^{N/2}(S_\alpha) =  R_{2, 2 \alpha + \tau}^{N/2} (S_\alpha)$.  Note that if $m=1$ then no other hypersurfaces overlap and the initial configuration is embedded.     

It will be important to preserve as many symmetries as possible in the construction of the approximate solution.  Note first that the displacement parameters along one geodesic are equal to those along the other; this is meant to ensure symmetry with respect to the rotation which takes $\gamma_1$ into $\gamma_2$.  Though not strictly necessary for the success of the construction (for example, it is possible to consider a different number of windings around each geodesic), this condition does simplify the forthcoming analysis.  It will be necessary, however, to impose reflection symmetries with respect to the $x^0$, $x^1$ and $x^2$ axes.  For this to be true, the displacement parameters must satisfy
\begin{equation}
	\mylabel{eqn:displacement}
	\sigma_0 = \sigma_{N/2} = 0 \qquad \mbox{and} \qquad
	\begin{gathered}
		\sigma_k = - \sigma_{N/2 - k} = \sigma_{N/2+k} = - \sigma_{N-k} \\
		\mbox{for all } \: k = 1, \ldots, (N-2)/4  \, .
	\end{gathered}
\end{equation}
Thus it will be assumed that these conditions are in force from now on. 

\paragraph{Perturbing the hyperspheres.}  Each hypersphere in $\Lambda_{\alpha, \tau, \sigma}^{\#}$ must be perturbed in the normal direction in order to improve its fit with the catenoids that will be used for gluing.  However, there are two changes from the method used before which need to be incorporated into the construction: the normal perturbation must be allowed to be different on each hypersphere; and the $k=0$ and the $k=N/2$ hyperspheres must be connected to \emph{four} neighbours.

First, choose $\eps_k>0$ and $b_k \in \R$ for $k = 0, \ldots, N-1$, and let $\tilde S_{\alpha, \eps_k, b_k} := \exp(\eps_k^{n-1} ( G + b_k x^1)(S_\alpha)$ be the perturbed hypersphere obtained from the normal perturbation corresponding to the function $ G + b_k x^1 : S_\alpha \setminus \{ p^+, p^-\} \rightarrow \R$.  Here, $G$ is the singular solution of the equation $\mathcal L_\alpha (G) = 0$ that was used before, while $x^1 = x^1 \big|_{S_\alpha}$ is the first coordinate function  restricted to $S_\alpha$ which is also a solution of the equation $\mathcal L_\alpha(x^1) = 0$ on $S_\alpha$.  The purpose of the quantity $b_k f$ is to allow the constant term in the asymptotic expansion of $G + b_k x^1$ near $p^\pm$ to vary.  Then for $s=1,2$ and all $k$ except $k=0, N/2$ one obtains the family of perturbed hyperspheres $R_{s, \sigma_{k}} \circ R_{s, 2 \alpha + \tau}^k (\tilde S_{\alpha, \eps_k, b_k})$.  

Next, let  $p_s^\pm := R_{s,\alpha}^{\pm 1} (p)$ be the four points where the geodesics $\gamma_s$ intersect $S_\alpha$.  Now let $G^{\, \prime} : S_{\alpha} \setminus \{p_1^+, p_1^-, p_2^+, p_2^- \} \rightarrow \R$ denote the singular solution of the equation $\mathcal L_\alpha (G^{\, \prime}) = 0$ defined by $G^{\, \prime} := c(G + G \circ Q)$ where $c$ is  the constant whose purpose is to give $G^{\, \prime}$ the same asymptotic expansion, given in equation \eqref{eqn:normgraphasym}, at all four singular points.  Similarly, let $f' = x^1 \big|_{S_\alpha} + x^2 \big|_{S_\alpha}$. Then for $\eps_0 = \eps_{N/2}>0$ and $b_0 = b_{N/2} \in \R$ one obtains the two perturbed hyperspheres 
\begin{align*}
	\tilde S_{\alpha,\eps_0, b_0}^{\, \prime} &:= \exp (\eps_0^{n-1} (G^{\, \prime} + b_0 f') N_\alpha)(S_\alpha) \\ 
	\tilde S_{\alpha,\eps_{N/2}, b_{N/2}}^{\, \prime} &:= - \exp (\eps_ {N/2} ^{n-1} (G^{\, \prime} + b_ {N/2} f') N_\alpha)(S_\alpha) \, .
\end{align*}

\paragraph{Inserting catenoids and assembling the approximate solution.} The appropriate perturbed, truncated and re-scaled catenoidal necks must now be found so that optimal matching with the perturbed hyperspheres is achieved.  However, this is a much more difficult problem than before because the non-zero displacement parameters alter the scaling parameters for each neck separately.  

To begin, choose $\bar \eps_k >0$ and $\bar b_k \in \R$ for $k = 0, \ldots, N-1$ and let $\tilde \Sigma_{\bar \eps_k, \bar b_k} := \tilde \Sigma_{\bar \eps_k} + \bar b_k$ denote the catenoidal neck defined in Section \ref{subsec:assembly} but translated along its axis by an amount $\bar b_k$.  The purpose of the parameter $\bar b_k$ is again to allow the constant term in the asymptotic expansion of the graphing function for the ends of the catenoid to vary.   The idea is now to choose $\eps_k, b_k, \bar \eps_k, \bar b_k$ so that for every $k$, the neck $R_{s, (\sigma_k+\sigma_{k-1})/2} \circ R^{ 2k+1}_{s,\alpha + \tau/2} \circ K^{-1}  \big( \tilde \Sigma_{\bar \eps_k, \bar b_k} \big)$ matches optimally with the perturbed hyperspheres $R_{s, \sigma_{k}} \circ R_{s, 2 \alpha + \tau}^k (\tilde S_{\alpha, \eps_k, b_k})$ and $R_{s, \sigma_{k+1}} \circ R_{s, 2 \alpha + \tau}^{k+1} (\tilde S_{\alpha, \eps_{k+1}, b_{k+1}})$ in the regions of overlap.  To this end, one writes the image of the $k^{\mathit{th}}$ and $(k+1)^{\mathit{st}}$ perturbed hyperspheres under the stereographic coordinate chart for the $k^{\mathit{th}}$ neck as graphs over the $\hat y $-hyperplane and equates the leading terms of the asymptotic expansions of the graphing functions given in equation \eqref{eqn:greensurfasym} with those of the graphing function of the catenoid. 

For the positive end of the neck overlapping with the $(k+1)^{\mathit{st}}$ perturbed hypersphere this comparison yields the equations
\begin{subequations}
\label{eqn:asymcomparison}
\begin{equation}
	\begin{array}{rcll}
		\bar \eps_k \log(2 / \bar \eps_k) + \bar \eps_k \log(\| \hat y \|) + \bar b_k \! \! \! &=& \! \! \!
		\tan(\tau_k/4) - \eps_{k+1} \big( - c_{2k} - C_{2k}  \log(\| \hat y \|) + b_{k+1} \big)
		&\quad n=2 \\[1ex]
		\bar \eps_k c_n - \dfrac{\bar \eps_k^{n-1}}{(n-2) \| \hat y \|^{n-2}} + \bar b_k \! \! \! &=& \! \! \!
		\tan(\tau_k/4) - \eps_{k+1}^{n-1} \Biggl( \dfrac{C_{nk}} {(n-2) \| \hat y \|^{n-2}} + b_{k+1} \Biggr)
		&\quad n \geq 3
	\end{array}
\end{equation}
where $C_{nk}$ and $c_{2k}$ are the quantities in \eqref{eqn:expansionconst} with $\tau$ replaced by $\tau_k$ and $c_n$ is the constant appearing in \eqref{eqn:catenoidasym}.  For the negative end of the neck overlapping with the $k^{\mathit{th}}$ perturbed hypersphere one finds
\begin{equation}
	\begin{array}{rcll}	
		\bar \eps_k \log(2 / \bar \eps_k) + \bar \eps_k \log(\| \hat y \|) - \bar b_k  \! \! \! &=& \! \! \! \tan(\tau_k/4) - \eps_{k} \big( - c_{2k} - C_{2k}  \log(\| \hat y \|) - b_{k} \big)
		&\quad n=2 \\[1ex]
		\bar \eps_k c_n - \dfrac{\bar \eps_k^{n-1}}{(n-2) \| \hat y \|^{n-2}} - \bar b_k   \! \! \! &=& \! \! \! \tan(\tau_k/4) - \eps_{k}^{n-1} \Biggl( \dfrac{C_{nk}} {(n-2) \| \hat y \|^{n-2}} - b_{k} \Biggr)
		&\quad n \geq 3
	\end{array}
\end{equation}
\end{subequations}
To solve equations \eqref{eqn:asymcomparison}, it is clear that one must choose $\eps_k := \eps$ for all $k$, where $\eps>0$ is still an unknown quantity; and then one must  choose $\bar \eps_k = \eps C_{2k}$ for all $k$ in dimension $n=2$ and $\bar \eps_k = \eps  C_{nk}^{1/(n-1)}$ for all $k$ in dimension $n\geq 3$.  The remaining equations for $b_k$ and $\bar b_k$ then become
\begin{equation}
	\label{eqn:paramchoice}
	\begin{aligned}
		\bar b_k + \frac{\eps^{n-1}}{2} \big( b_{k+1} + b_k \big) &= 0 \\
		\frac{\eps^{n-1}}{2} \big( b_{k+1} - b_k \big) & = T_{nk}
	\end{aligned}
\end{equation}
where $T_{2k} = \tan(\tau_k / 4) + \eps C_{2k} \big( 1 - \log( 2/(\eps C_{2k}) \big)$ and $T_{nk} = \tan(\tau_k / 4) - \eps c_n C_{nk}^{1/(n-1)}$. One can check that this $2N \times 2N$ linear system has a one-dimensional kernel spanned by $b_k = \eps^{1-n}$ and $\bar b_k = -1$ for all $k$ and that the condition $\sum_{k=0}^{N-1} T_{nk} = 0$ guarantees that a solution of \eqref{eqn:paramchoice} exists.  This condition can be rearranged as
\begin{equation}
	\label{eqn:epschoice}
	\sum_{k=0}^{N-1} \tan(\tau_k/4) = 
	\begin{cases}
		\displaystyle \eps \sum_{k=0}^{N-1} C_{2k} \big( \log \bigl(2 / (\eps C_{2k}) \bigr) - 1 \big)   &\qquad n=2 \\[3ex]
		\displaystyle \eps \sum_{k=0}^{N-1} c_n  C_{nk}^{1/(n-1)}    & \qquad n\geq 3
	\end{cases}
\end{equation}
so that the parameter $\eps$ can be determined uniquely from $\tau$ and the displacement parameters in this way.  One can then find a solution $b_k$ and $\bar b_k$ that satisfies $b_0 = 0$ and thus $b_{N/2} = 0$ by the symmetry condition \eqref{eqn:displacement}.  It can be shown that for $k = 1, \ldots, \frac{1}{2}(N/2 - 1)$ quantities on the $k^{\mathit{th}}, (N/2 -k -1)^{\mathit{st}}, (N/2+k)^{\mathit{th}}$ and $(N-k-1)^{\mathit{st}}$ necks are equal and quantities on the $k^{\mathit{th}}, (N/2 -k )^{\mathit{th}}, (N/2+k)^{\mathit{th}}$ and $(N-k)^{\mathit{th}}$ perturbed hyperspheres are equal for the same reason.  This solution satisfies the estimates $\max_k \bigl( \eps^{n-1} |b_k| + | \bar b_k| \bigr) = \mathcal O(\eps)$.  

\begin{rmk}
	The estimate for $|b_k|$ will ensure that the errors introduced in the mean curvature of $S_\alpha$ by the perturbation corresponding to $\eps_k^{n-1} ( G + b_k x^1)$ are of the same magnitude as in Proposition \ref{prop:error} when the $b_k x^1$ term was absent. 
\end{rmk} 	

Once the parameters $\eps_k, b_k, \bar \eps_k$ and $\bar b_k$ have been chosen according to the procedure above, one can then assemble the approximate solution in the same way as before.  Note that the condition $b_0 = b_{N/2} = 0$ ensures proper matching of all four perturbed hyperspheres near $\tilde S_{\alpha,\eps_0}^{\, \prime}$ and $\tilde S_{\alpha,\eps_{N/2}}^{\, \prime}$.

\begin{defn}
	\mylabel{defn:secapproxregions}
	Recall the terminology from Section \ref{subsec:assembly}. 
	\begin{itemize}
		\item Define the neck regions
		\begin{align*}
			\mathcal N^{s,k, \pm}_{\eps} &:= R_{s, (\sigma_k+\sigma_{k+1})/2} \circ R^{ 2k+1}_{s,\alpha + \tau/2} \circ K^{-1}  \big( \tilde \Sigma_{\bar \eps_k, \bar b_k}^\pm \cap \mathit{Cyl}(\rho_{\eps}/2) \big) \, .
		\end{align*}
	Set $\mathcal N_{\eps}^{s,k} := \mathcal N_{\eps}^{s,k,+} \cup \mathcal N_{\eps}^{s,k,-} $.  These necks are oriented so that $\mathcal N_{\eps}^{s,k,+}$ lies ahead of $\mathcal N_{\eps}^{s,k,-}$ along the geodesic $\gamma_s$, and that the $k^{\mathit{th}}$ neck is centered on the point $R_{s,(\sigma_k+ \sigma_{k+1})/2} \circ R^{2k+1}_{s, \alpha + \tau/2}(p)$.  
	
		\item Define the transition regions
		\begin{align*}
			\mathcal T^{s,k, \pm}_{\eps} &:= R_{s, (\sigma_k + \sigma_{k+1})/2} \circ R^{ 2k+1}_{s,\alpha + \tau/2} \circ K^{-1}  \big( \tilde \Sigma_{\bar \eps_k, \bar b_k}^\pm \cap \big[ \mathit{Cyl}(2\rho_{\eps}) \setminus \mathit{Cyl}(\rho_{\eps}/2) \big] \big) \, .
		\end{align*}
	Set $\mathcal T_{\eps}^{s,k} := \mathcal T_{\eps}^{s,k,+} \cup \mathcal T_{\eps}^{s,k,-} $.  	
	
		\item Define the exterior regions
		\begin{align*}
			\mathcal E_{\eps}^0 &:= \tilde S_{\alpha, \eps_0, b_0}^{\, \prime} \setminus \big[ B_{\tilde \rho_{\eps}} (p_1^+) \cup B_{\tilde \rho_{\eps}} (p_1^-) \cup B_{\tilde \rho_{\eps}} (p_2^+) \cup  B_{\tilde \rho_{\eps}} (p_2^-) \big] \\[0.5ex]
			\mathcal E_{\eps}^{N/2} &:= - \mathcal E_{\eps}^{0} \\[0.5ex]
			\mathcal E_{\eps}^{s,k} &:= R_{s,\sigma_k} \circ R^{k}_{s,2 \alpha + \tau} \big( \tilde S_{\alpha, \eps_k, b_k} \setminus \big[ B_{\tilde \rho_{\eps}} (p_s^+) \cup B_{\tilde \rho_{\eps}} (p_1^-)\big] \big) \, .
		\end{align*} 
	\end{itemize}
\end{defn}
\noindent The approximate solution corresponding to the configuration $\Lambda_{\alpha,\tau,\sigma}^{\#}$ can now be defined as follows.
\begin{defn}
	\mylabel{defn:secapproxsol}
	The \emph{approximate solution} with parameters $\alpha, \tau$ and $\sigma$ is the hypersurface 
	\begin{equation}
		\mylabel{eqn:secondaproxsol}
		\secapprox := \mathcal E_{\eps}^0 \cup  \mathcal E_{\eps}^{N/2}  \cup \left[ \bigcup_{s=1}^2\bigcup_{k =1}^{ N/2 - 1} \!\! \Big(  \mathcal E_{\eps}^{s,k} \cup \mathcal E_{\eps}^{s,N-k}\cup \mathcal T_{\eps}^{s, k} \cup \mathcal N_{\eps}^{s,k} \cup \mathcal T_{\eps}^{s, N-k} \cup \mathcal N_{\eps}^{s,N-k}  \Big) \right] \, .
	\end{equation}	
	This definition contains no redundancy.  
\end{defn}

\subsection{Symmetries of the Approximate Solution}

The approximate solution $\secapprox$ is considerably less symmetric than before since symmetry with respect to each $R_{s,2\alpha+\tau}$ separately no longer holds.  However, it is important to realize that $\secapprox$ does still possess many symmetries that can be exploited to reduce the number of approximate Jacobi fields that will have to be dealt with during the perturbation of $\secapprox$ into a CMC hypersurface.  These symmetries are as follows.

\medskip
\noindent 1. \itshape Transverse rotational symmetries. \upshape
\medskip

Define the transformation $S^{\, 012}_{B} \in O(n+2)$ by choosing $B \in O(n-1)$ and then setting
	$$S^{\, 012}_B := \left( 
		\begin{array}{c|c} 
			\! \! \mbox{\scriptsize$\begin{array}{ccc}
				1 && \\[-0.5ex]
				& 1&  \\[-0.5ex]
				&&1
			\end{array}$} \! \!& 0\\ 
			\hline 0 & B 
		\end{array}
		\right)$$
	which is the transformation keeping the $x^0$, $x^1$ and $x^2$ coordinates fixed while transforming the remaining coordinates by $B$. Each $S^{\, 012}_{B}$ is a symmetry of $\secapprox$ by arguments similar to those for the analogous symmetries in the construction of the Delaunay-like hypersurfaces in the previous section of this paper.
	
\medskip
\noindent 2. \itshape Reflection symmetries. \upshape
\medskip

For $s \in 0, 1, 2$, define the reflections $T_s \in O(n+2)$ by 
	$$T_s(x^0, \ldots, x^n) := (x^0, \ldots, -x^s, \ldots, x^n) \, .$$  
	These reflections satisfy: $T_{s} \circ R_{1, \theta} \circ T_{s}  = R_{1,\theta}^{-1}$ for $s=0,1$ and $T_{s} \circ R_{2, \theta} \circ T_{s}  = R_{2,\theta}^{-1}$ for $s=0,2$; as well as $T_{1} \circ R_{2, \theta} \circ T_{1}  = R_{2,\theta}$ and $T_{2} \circ R_{1, \theta} \circ T_{2}  = R_{1,\theta}$.  The initial configuration of un-connected hyperspheres $\Lambda_{\alpha, \tau, \sigma}^{\#}$ is invariant under these symmetries because the displacement parameters satisfy \eqref{eqn:displacement}. Since the neck regions are exactly determined by the separation and displacement parameters (via the neck scale parameters) further consequences are:  $T_0$ reflects each neck region parallel to its axis and then permutes it with another neck region on the same geodesic;  $T_1$ reflects the neck regions on $\gamma_1$ parallel to their axis and permutes them amongst themselves and reflects the neck regions on $\gamma_2$ perpendicular to their axis and permutes them amongst themselves;  and $T_2$ reflects the neck regions on $\gamma_2$ parallel to their axis and permutes them amongst themselves and reflects the neck regions on $\gamma_1$ perpendicular to their axis and permutes them amongst themselves.  Hence $\secapprox$ is invariant under $T_0, T_1, T_2$.
	
\medskip
\noindent 3. \itshape Parallel rotational symmetry. \upshape
\medskip

Let $Q \in SO(n+1)$ denote the rotation 
	\begin{equation}
		\mylabel{eqn:x1x2rot}
		Q := \left(
		\begin{array}{c|c}
			\begin{array}{ccc}
				1 & 0 & 0  \\
				0 & 0 & -1  \\
				0 & 1 & 0
			\end{array} & 0 \\
			\hline
			0  & I
		\end{array} \right)
	\end{equation}
which rotates $\gamma_1$ into $\gamma_2$. Furthermore, $Q \circ R_{1,\theta} \circ Q^{-1} = R_{2,\theta}$ so that $\Lambda_{\alpha, \tau, \sigma}^{\#}$ is invariant under $Q$.  Also, $Q$ rotates each neck region on $\gamma_1$ to a neck region on $\gamma_2$.  Hence $\secapprox$ is invariant under $Q$.

\subsection{The Analytic Set-Up}
\mylabel{sec:secanalysis}

For each sufficiently small choice of displacement parameter $\sigma$, the approximate solution $\secapprox$ has mean curvature equal to $H_\alpha$ everywhere except in the neck regions.  Each of these hypersurfaces is thus a candidate for deformation into an exactly CMC hypersurface using the methods of the previous section.  However, it will turn out that because of the loss of symmetry, only a specific choice of $\sigma$ leads to a candidate for which the deformation can be carried out successfully, and identifying this $\sigma$ requires new techniques.

The analytic set-up for the present deformation problem, however, begins in the same way as before.  The approximate solution $\secapprox$ will be deformed equivariantly with respect to the symmetries $T_0$, $T_1$, $T_2$, $Q$ and $S^{\, 012}_{B}$ for all $B \in SO(n-2)$ by using normal deformations generated by functions invariant under these symmetries.  
\begin{defn}
	Let $X := \{ f : \secapprox \rightarrow \R \: : \:  f \circ T_s =  f \circ S^{012}_B = f \circ Q = f \:\:\: \forall \: s=0,1,2 \: \mbox{and} \: B \in SO(n-2) \}$.
\end{defn}
\noindent Next, let $\phi_f(y) := \exp_y (f(y) \cdot N(y))$ denote the deformation of $\secapprox$ obtained by taking the exponential of the vector field $f \cdot N$ where $N$ is the unit normal vector field of $\secapprox$ and $f$ is a function on $\secapprox$.  
\begin{defn} 
Let $\secdifop$ be the operator $f \longmapsto H_{\phi_f(\secapprox)} - H_\alpha$ with linearization at zero $\seclinop$.  
\end{defn}
\noindent Finally, let $C^{k,\beta}_\delta$ denote a weighted H\"older space defined as in Section \ref{sec:funcspace} except where the weight function is modified to take into account the additional perturbed hyperspheres and neck regions in $\secapprox$.  The operator $\secdifop : C^{2,\beta}_\delta (X) \rightarrow C^{0,\beta}_{\delta - 2} (X)$ is a well-defined, non-linear partial differential operator on $\secapprox$ whose zero, if one can be found, corresponds to a deformation of $\secapprox$ to an exactly CMC hypersurface.

The Banach space inverse function theorem will once again be used to find a zero of $\secdifop$ in $C^{2, \beta}_\delta( X )$ for the right choice of $\sigma$.  However, it turns out that not all the approximate Jacobi fields of $\seclinop$ are excluded by the symmetries.  To see which ones, proceed as follows.  First define smooth and monotone cut-off functions $\chi^{s,k}$ for $k \neq 0, N/2$ as well as  $\chi^0$ and $\chi^{N/2}$ on $\secapprox$ that are supported, respectively, in subsets of $\mathcal E_\eps^{s,k}$, $\mathcal E_\eps^0$ and $\mathcal E_\eps^{N/2}$ that will be made precise later.  Let $x^t $ be the $t^{\mathit{th}}$ coordinate function for $t = 1, \ldots , n$.  Define the following functions on $\secapprox$:
\begin{equation}
	\mylabel{eqn:approxjac}
	\begin{aligned}
		\tilde q_{s,k}^t  &:= 
			\chi^{s,k}  \cdot \left( x^t  \big\vert_{\tilde S_{\alpha, \eps_k, b_k} } \circ R_{s,2 \alpha + \tau}^{-k} \circ R^{-1}_{s,\sigma_k} \right) \\[0.5ex]
		\tilde q^t_0 &:= \chi^0 \cdot \left( x^t  \big\vert_{\tilde S_{\alpha, \eps_0}^{\, \prime}}\right)  \\[0.5ex]
			\tilde q^t_{ N/2} &:= \chi^{N/2}  \cdot \left( x^t  \big\vert_{- \tilde S_{\alpha, \eps_{N/2}}^{\, \prime}}\right) \, .
	\end{aligned}
\end{equation}
Then according to the discussion of Section \ref{sec:jacobi}, the vector space of \emph{all} approximate Jacobi fields in $C^{2, \beta}(\secapprox)$ is $\tilde{ \mathcal K} := \linspan_\R \{ \tilde q_{s,k}^t :\mbox{all} \: s, t, k \}$.  It remains to identify the subspace of these functions which have the symmetries needed to belong to  $C^{2, \beta}_\delta ( X)$.  By using the commutation relations between the various symmetries and the rotations $R_{s,\theta}$, one finds the following subspace.  

\begin{defn}
	\mylabel{defn:apkerbasis}
	The subspace of approximate Jacobi fields of $\seclinop$ in $C^{2,\beta}_\delta (X)$ is the subspace of functions  
	\begin{equation*}
		\tilde{\mathcal K}_{\mathit{sym}} := \linspan_\R \left\{ 
		\sum_{s=1}^2 \big( \tilde q^1_{s,k} - \tilde q^1_{s, N/2 - k} + \tilde q^1_{s,N/2 + k} - \tilde q^1_{s, N - k} \big) \: : \: k = 1, \ldots, \tfrac{N-2}{4}
		\right\} \, .
	\end{equation*}
	Denote the functions spanning $\tilde{\mathcal K}_{\mathit{sym}}$ by $\tilde q_k$.  Note that $\mathrm{dim} (\mathcal{\tilde K}_{\mathit{sym}}) = (N-2)/4$.
\end{defn}

\noindent Essentially, each function $\tilde q_k$ is identical on the $k^{\mathit{th}}$, $(N/2 - k)^{\mathit{th}}$, $(N/2+k)^{\mathit{th}}$ and $(N-k)^{\mathit{th}}$ perturbed hyperspheres, and they pull back to the coordinate function $x^1$ restricted to $\tilde S_{\alpha, \eps_k, b_k}$.

\subsection{The New Linear Estimate}

Because of the presence of the approximate Jacobi fields in $\tilde{\mathcal K}_{\mathit{sym}}$, the linearized operator $\seclinop$ on $C^{2,\beta}_\delta(X)$ no longer has an inverse whose norm is uniformly bounded with respect to the singular parameter $\tau$.  One way out of this difficulty is to project $\seclinop$ onto a subspace of $\cobg(X)$ which is transverse to the approximate co-kernel associated to $\tilde{\mathcal K}_{\mathit{sym}}$, and construct a bounded right inverse for the projected operator.  Since  $\seclinop$ is self-adjoint, an appropriate subspace to choose is the $L^2$-orthogonal complement of $\tilde{\mathcal K}_{\mathit{sym}}$, denoted $X^\perp$.  Let $\pi : X \rightarrow X^{\perp}$ denote the $L^2$ projection onto $X^{\perp}$ and set $\seclinop^\perp := \pi \circ \seclinop$. 

\newcommand{\secrinv}{\tilde{\mathcal R}_{\alpha, \tau, \sigma}}

\begin{prop}
	\mylabel{thm:seclinest}
	Suppose first that the dimension of $\secapprox$ is $n\geq3$.  Choose $\delta \in (n-2, 0)$ and  sufficiently small separation and displacement parameters $\tau$ and $\sigma$. Then the linearized operator $\seclinop^\perp : \ctbg (X) \rightarrow \cobg (X^{\perp})$ possesses a bounded right inverse $\secrinv$ satisfying the estimate
	$$| \secrinv(f) |_{\ctbg} \leq C | f |_{\cobg}$$
where $C$ is a constant independent of $\tau$ and $\sigma$.  If the dimension of $\secapprox$ is $n=2$ then one can choose $\delta \in (-1, 0)$ and find a right inverse satisfying the estimate
	$$| \secrinv(f) |_{\ctbg} \leq C \eps^{\delta} | f |_{\cobg}$$
where  $C$ is a constant independent of $\tau$ and $\sigma$.
\end{prop}

\begin{proof}
	
	The proof of this result is mostly the same as the proofs of Proposition \ref{thm:linesthigh} and Proposition \ref{thm:linest}, but with a few essential differences.  The idea is to solve the equation $\seclinop (u) = \pi(f)$ in $C^{2,\beta}_\delta(X)$ using the patching technique to obtain $u := \secrinv(f)$ and then to estimate $|u|_{\ctbg}$ in terms of $|f|_{\cobg}$.  
	
	For any $\rho \in (\tilde \rho_\eps, \rho_0)$ define the enlarged neck regions and exterior regions $\mathcal N^{s,k}_\eps(\rho)$ and $\mathcal E^{s,k}_\eps(\rho)$ for each $s = 1, 2$ and $k = 1, \ldots,N/2-1, N/2 +1, \ldots, N-1$ in the same way as before.  Define also the cut-off function $\chi_{\mathit{neck}, \rho}^{s,k}$ and $\chi_{\mathit{ext},\rho}^{s,k}$ corresponding to $\mathcal N^{s,k}_\eps(\rho)$ and $\mathcal E^{s,k}_\eps(\rho)$, respectively, as well as the cut-off functions $\eta_{\mathit{neck}}^{s,k}$ and $\eta_{\mathit{ext}}^{s,k}$ corresponding to $\mathcal N^{s,k}_\eps$ and $\mathcal E^{s,k}_\eps$.  Finally, define enlarged regions $\mathcal E_\eps^0(\rho)$ and $\mathcal E_\eps^{N/2}(\rho)$ and their corresponding cut-off functions in an analogous manner (taking into account that $\mathcal E_\eps^0$ and $\mathcal E_\eps^{N/2}$ are perturbed hyperspheres with four points removed).   Suppose now that $n \geq 3$.  Choose $\pi(f) \in \ctbg (X^{\perp})$ with $\delta \in (2-n,0)$ and decompose $\pi(f)$ as 
\begin{gather*}
	\pi(f) = f^{0}_{\mathit{ext}} + f^{N/2}_{\mathit{ext}} + \sum_{s=1}^2 \left( \sum_{k=1}^{N/2 -1} \big( f^{s,k}_{\mathit{ext}} + f^{s,N-k}_{\mathit{ext}} \big) + \sum_{k=0}^{N-1} f^{s,k}_{\mathit{neck}} \right) \\[0.5ex]
	\intertext{where} 
	\begin{aligned}
		f^{0}_{\mathit{ext}} &:= \chi^{0}_{\mathit{ext}} \cdot \pi(f) \\
		f^{N/2}_{\mathit{ext}} &:= \chi^{N/2}_{\mathit{ext}} \cdot \pi(f)
	\end{aligned} \qquad \mbox{and} \qquad
	\begin{aligned}
		f^{s,k}_{\mathit{ext}} &:= \chi^{s,k}_{\mathit{ext}}  \cdot \pi(f) \\ 
		f^{s,k}_{\mathit{neck}} &:= \chi^{s,k}_{\mathit{neck}}  \cdot \pi(f) \, .
	\end{aligned}
\end{gather*}

\begin{rmk}
	Use $\chi_{\mathit{ext}}^{s,k}$, $\chi^0_{\mathit{ext}}$ and $\chi^{N/2}_{\mathit{ext}}$ as the cut-off functions in Definition \ref{defn:apkerbasis}.
\end{rmk}
	
Observe that the following symmetries hold for these various functions.  First, the invariance with respect to $Q$ and $T_0, T_1, T_2$ implies  $f_{\mathit{ext}}^0 = f_{\mathit{ext}}^{N/2}$ and 
\begin{gather*}
	f^{1,k}_{\mathit{ext}} = f^{1, N/2-k}_{\mathit{ext}} = f^{1, N/2+k}_{\mathit{ext}} = f^{1,N-k}_{\mathit{ext}} = f^{2,k}_{\mathit{ext}} = f^{2, N/2-k}_{\mathit{ext}} = f^{2, N/2+k}_{\mathit{ext}} = f^{2,N-k}_{\mathit{ext}} \\
	f^{1,k}_{\mathit{neck}} = f^{1, N/2-k-1}_{\mathit{neck}} = f^{1, N/2+k}_{\mathit{neck}} = f^{1,N-k-1}_{\mathit{neck}} = f^{2,k}_{\mathit{neck}} = f^{2, N/2-k-1}_{\mathit{neck}} = f^{2, N/2+k}_{\mathit{neck}} = f^{2,N-k-1}_{\mathit{neck}} \, . 
\end{gather*}
In addition, $f^{1,(N-2)/4}_{\mathit{neck}}$ is an even function.  Next, choose $k = 1, \ldots, N/2 - 1$ and consider what effect the invariance with respect to $S^{012}_B$ for $B \in O(n-1)$ and $T_{3-s}$ has on these functions.  The function $f^{s, k}_{\mathit{neck}}$ can be pulled back to a function on $\tilde \Sigma_{\bar \eps_k}$.  Write the pull back function as $f^{s, k}_{\mathit{neck}} := f^{s, k}_{\mathit{neck}}( y^2, y^3, \ldots, y^{n-1})$ using a parametrization for $\tilde \Sigma_{\bar \eps_k}$ as a cylindrically symmetric graph in $\R \times \R^{n}$. Then one can check that $f^{s, k}_{\mathit{neck}}$ satisfies the symmetries 
\begin{equation}
	\mylabel{eqn:symneck}
	f^{s, k}_{\mathit{neck}} (y^2, B(y^3, \ldots, y^{n+1}) ) = f^{s, k}_{\mathit{neck}} (-y^2, y^3, \ldots, y^{n+1}) = f^{s, k}_{\mathit{neck}}(y^2, y^3,\ldots, y^{n+1})
\end{equation}
where $B$ acts on $(y^3, \ldots y^{n+1}) \in \R^{n-1}$ in the obvious way.   Similarly, the function $f^{s, k}_{\mathit{ext}}$ can be pulled back to a function on $\tilde S_{\alpha, \eps_k, b_k}$.  Write the pull-back function as $f^{s, k}_{\mathit{ext}} := f^{s, k}_{\mathit{ext}}(\mu, \Theta)$ using $(\mu, \Theta) \in (0, \pi) \times \Sph^{n-1}$ parametrization of $\tilde S_{\alpha, \eps_k, b_k}$. Then one can check that $f^{s, k}_{\mathit{ext}}$ satisfies the symmetries 
\begin{equation}
	\mylabel{eqn:symext1}
	f^{s, k}_{\mathit{ext}} (\mu, \hat S_B (\Theta) ) = f^{s, k}_{\mathit{ext}} (\mu, \hat T(\Theta)) = f^{s, k}_{\mathit{ext}}(\mu, \Theta)
\end{equation}
where $\hat S_B$ is the transformation of $\Sph^{n-1} \subseteq \R^n$ given by the matrix $\hat S_B := \left( \begin{smallmatrix} 1 & \\ & B\end{smallmatrix} \right)$, and $\hat T$ is the transformation of $\Sph^{n-1}$ given by the matrix $\hat T := \left( \begin{smallmatrix}-1 & \\ & I \end{smallmatrix} \right)$ where $I$ is the $(n-1) \times (n-1)$ identity matrix.   Finally, one can pull back $f^0_{\mathit{ext}}$ and $f^{N/2}_{\mathit{ext}}$ to  functions $f^\ast_{\mathit{ext}}:=f^\ast_{\mathit{ext}}(\mu, \Theta)$ on $\tilde S_{\alpha,\eps_0}^{\, \prime}$.  Then one can check that these functions satisfy the symmetries 
\begin{equation}
	\mylabel{eqn:symext2}
	f^\ast_{\mathit{ext}} (-\mu, \Theta) = f^\ast_{\mathit{ext}}(\mu, \hat T(\Theta)) = f^\ast_{\mathit{ext}}(\mu, \hat S_B (\Theta)) = f(\mu, \Theta) \, .
\end{equation}
This is because \emph{both} symmetries $T_1$ and $T_2$ fix $\mathcal E_{\eps}^0(\rho)$ and $\mathcal E^{N/2}_\eps(\rho)$.

The first step in the patching technique is to find local solutions in the neck regions.  That is, one solves the equation $\frac{1}{4}\mathcal L_{\eps_k \Sigma} (u_{\mathit{neck}}^{s,k}) = f^{s, k}_{\mathit{neck}}$ in $C^{2,\beta}_{\delta, \mathit{sym}}(\eps_k \Sigma)$ exactly as before (using the fact that the symmetries \eqref{eqn:symneck} ensure that $C^{2,\beta}_{\delta, \mathit{sym}}(\eps_k \Sigma)$ for $\delta \in (2-n, 0)$ contains none of the Jacobi fields of $\mathcal L_{\eps_k \Sigma}$ listed in \eqref{eqn:catjacobi} and thus $\frac{1}{4}\mathcal L_{\eps_k \Sigma}$ is bijective with $\delta$ in this range).  The solution can then be extended to all of $\secapprox$ as before and the extended function $\bar u_{\mathit{neck}}^{s,k}$ satisfies the estimate $|\bar u_{\mathit{neck}}^{s,k}|_{\ctbg} \leq C |\pi(f)|_{\cobg}$.   

The next step is to find local solutions on the exterior regions.  Choose a small $\kappa \in (0,1)$ and define the functions 
\begin{equation}
	\mylabel{eqn:secextfns}
	\begin{aligned}
		\hat f^{s,k}_{\mathit{ext}} &:= \chi_{\mathit{ext},\kappa \rho}^{s,k} \big( \pi(f) - \sum_{s',k'} \seclinop^\perp( \bar u_{\mathit{neck}}^{s',k'}) \big) \\
		\hat f^0_{\mathit{ext}} &:= \chi_{\mathit{ext},\kappa \rho}^{0} \big(\pi(f) - \sum_{s',k'} \seclinop^\perp( \bar u_{\mathit{neck}}^{s',k'}) \big) \\
		\hat f^{N/2}_{\mathit{ext}} &:= \chi_{\mathit{ext},\kappa \rho}^{N/2} \big( \pi(f) - \sum_{s',k'} \seclinop^\perp( \bar u_{\mathit{neck}}^{s',k'}) \big) \, .
	\end{aligned}
\end{equation}
One can view each of the functions \eqref{eqn:secextfns} as $C^{0,\beta}$ functions on $S_{\alpha}$ and attempt to solve the equation $\mathcal L_\alpha (u_{\mathit{ext}}^\ast) =  \hat f_{\mathit{neck}}^\ast$ up to projection onto the approximate co-kernel.

Begin with $\hat f_{\mathit{neck}}^{s,k}$.  The symmetries \eqref{eqn:symext1} imply that the only invariant Jacobi field of $\mathcal L_\alpha$ is the coordinate function $x^1$.   Compute $\lambda^{s,k} := \left. \int_{S_\alpha} \hat f^{s,k}_{\mathit{ext}} \cdot x^1 \right/ \int_{S_\alpha} \tilde q_k \cdot x^1$ so that $\big( \hat f^{s,k}_{\mathit{ext}} \big)^\perp := \hat f^{s,k}_{\mathit{ext}}  - \lambda^{s,k} \tilde q_k$ is orthogonal to $x^1$ where $\tilde q_k \in \tilde{ \mathcal K}_{\mathit{sym}}$ is defined in Definition \eqref{defn:apkerbasis}.   The equation $\mathcal L_\alpha (u^{s,k}_{\mathit{ext}}) = \big( \hat f^{s,k}_{\mathit{ext}} \big)^\perp$ can now be solved for $u^{s,k}_{\mathit{ext}}$ in $C^{2,\beta}(S_\alpha)$.  Moreover,  it is true that 
$$\int_{S_\alpha} \hat f^{s,k}_{\mathit{ext}} \cdot x^1  = \int_{S_\alpha \cap \mathit{Tub}_{2 \rho}(\gamma)} \hat f^{s,k}_{\mathit{ext}} \cdot x^1 + \int_{S_\alpha \setminus \mathit{Tub}_{2 \rho}(\gamma)} \hat f^{s,k}_{\mathit{ext}} \cdot x^1$$
and 
$$\int_{S_\alpha \setminus \mathit{Tub}_{2 \rho}(\gamma)} \hat f^{s,k}_{\mathit{ext}} \cdot x^1 = \int_{S_\alpha \setminus \mathit{Tub}_{2 \rho}(\gamma)} f^{s,k}_{\mathit{ext}} \cdot x^1  = \int_{\secapprox} f \cdot \tilde q_k = 0$$
using the fact that $\pi(f) \in \cobg(X^\perp)$.  Thus with the help of estimates similar to \eqref{eqn:finalneckest} one deduces 
$$|\lambda^{s,k}| \leq C \kappa \big(\kappa^2 +  \eps^{2(3n-3)/(3n-2)} \big) |\pi(f)|_{\cobg}$$ 
where $C_\kappa$ depends on $\kappa$ and $\delta$ but not $\eps$.   The solution therefore satisfies the estimate $|u^{s,k}_{\mathit{ext}}|_{C^{2,\beta}(S_\alpha)} \leq C_\kappa |\pi(f)|_{\cobg}$.  The remaining functions $f^0_{\mathit{ext}}$ and $f^{N/2}_{\mathit{ext}}$ are easier.  The symmetries \eqref{eqn:symext2} imply that there are no invariant Jacobi fields of $\mathcal L_\alpha$ to worry about.  Hence solutions of $\mathcal L_{\alpha} ( u^\ast_{\mathit{ext}} ) = f^\ast_{\mathit{ext}}$ can be found without further ado.  The same kind of estimate holds for these solutions.

The extension of the various local solutions on the perturbed hyperspheres to all of $\secapprox$ is slightly more involved than in the previous proof.  Suppose that $u^{s,k}_{\mathit{ext}}(p^\pm) := a^{s,k}_{\pm}$.  If $c_0^{s,k}$ and $c_1^{s,k}$ satisfy the equations $a^{s,k}_+ = c_0^{s,k} - c_1^{s,k}$ and $a^{s,k+1} = c_0^{s,k} + c_1^{s,k}$ then the leading term in the asymptotic expansion of the function $c_0^{s,k} J_0 + c_1^{s,k} J_1$ on the $k^{\mathit{th}}$ neck matches the leading term of the Taylor expansion of $u^{s,k}_{\mathit{ext}}$ on the overlap with the $k^{\mathit{th}}$ perturbed hypersphere and with the leading term of the Taylor expansion of $u^{s,k+1}_{\mathit{ext}}$ on the overlap with the $(k+1)^{\mathit{st}}$ perturbed hypersphere.  This can be done for $k=2, \ldots, N/2-1$ and $k=N/2+2, \ldots, N-1$.   By symmetry, the quantities  $u^{s, 1}_{\mathit{ext}} (p^-)$, $u^{s, N/2-1}_{\mathit{ext}} (p^+)$, $u^{s, N/2+1}_{\mathit{ext}} (p^-)$ and $u^{s, N-1}_{\mathit{ext}} (p^+)$ are all equal, as are $u^0_{\mathit{ext}}(p_s^\pm)$ and $u^{N/2}_{\mathit{ext}}(p_s^\pm)$.  Hence on can add the same multiples of $J_0$ and $J_1$ in the necks corresponding to $k=1, N/2, N/2+1$ and $N$ to match with the leading term in the Taylor series of $u^0_{\mathit{ext}}$ and $u^{N/2}_{\mathit{ext}}$ on the regions of overlap with these necks.  Finally, one arrives at the extended solution
$$\bar u_{\mathit{ext}} := \eta_{\mathit{ext}}^0 u^0_{\mathit{ext}} + \eta_{\mathit{ext}}^{N/2} u^ {N/2}_{\mathit{ext}} + \sum_{s=1}^2 \sum_{\substack{k=1 \\ k \neq N/2}}^{N-1} \eta_{\mathit{ext}}^{s,k} u_{\mathit{ext}}^{s,k} + \sum_{s=1}^2\sum_{k=1}^N \eta_{\mathit{neck}}^{s,k} \big( c_0^{s,k} J_0 + c_1^{s,k} J_1 \big) \, .$$
The extended function $\bar u_{\mathit{ext}}$ satisfies the estimate $|\bar u_{\mathit{ext}}|_{\ctbg} \leq C |\pi(f)|_{\cobg}$.

The proof now concludes as follows.  Set $\bar u := \bar u_{\mathit{ext}} + \sum_{s,k}  \bar u^{s,k}_{\mathit{neck}}$. Then estimates similar to those used before lead to the bounds  
\begin{equation}
	\mylabel{eqn:seccvgce}
	| \seclinop (\bar u ) - \pi(f) |_{\cobg} \leq C \big( \kappa^2  +  \eps^{2(3n-3)/(3n-2)} \big) |\pi(f)|_{\cobg}
\end{equation}
and $| \bar u |_{\ctbg} \leq C |\pi(f)|_{\cobg}$ for constants $C$ independent of $\eps$.  Since the factor in \eqref{eqn:seccvgce} that is in front of $|\pi(f)|_{\cobg}$ can be made as small as desired by choosing $\kappa$ and $\eps$ sufficiently small, an exact solution of $\seclinop^\perp (u) = f$ can be found by iteration.  It satisfies the bound $| u |_{\ctbg} \leq C |\pi(f)|_{\cobg}$.  This completes the proof in the case $n \geq 3$.	The case $n=2$ is similar, except that the double indicial root of $\seclinop$ must be dealt with by matching asymptotic expansions as in Proposition \ref{thm:linest} and one ends up with the estimate $| u |_{\ctbg} \leq C \eps^{\delta} |\pi(f)|_{\cobg}$.  
\end{proof}

\subsection{The Solution of the Non-Linear Problem up to Finite-Dimensional Error}

The new linear estimate derived in the previous section and non-linear estimates identical to those of Section \ref{subsec:nonlinest} can now be combined to solve the equation $\secdifop(f) = 0$ using the Banach space inverse function theorem up to a finite-dimensional error term in the following sense.  Let $\pi$ denote once again the $L^2$-projection onto $X^\perp$.

\begin{prop}
	\mylabel{prop:secperturb}
	If $\tau$ and $\Vert \sigma \Vert $ are sufficiently small, then there exists $f_{\alpha, \tau, \sigma} \in \ctbg (X)$ satisfying $\pi \circ \secdifop (f_{\alpha, \tau, \sigma}) = 0$ and there exists a constant $C$ independent of $\tau$ and $\sigma$ so that 
	$$| f_{\alpha, \tau, \sigma} |_{\ctbg} \leq C \cdot C_L(\eps) \eps^{(2-\delta)(3n-3)/(3n-2)}$$
	where $C_L = \mathcal O(1)$ in dimension $n\geq 3$ and $C_L = \mathcal O(\eps^{\delta})$ in dimension $n = 2$.   Therefore the hypersurface obtained by deforming $\secapprox$ in the normal direction by an amount determined by $f_{\alpha, \tau, \sigma}$ has constant mean curvature and is embedded if $\secapprox$ is.
\end{prop}

\begin{proof}
The linearization of $\pi \circ \secdifop$ at zero is $\Dif ( \pi \circ \secdifop )(0) = \seclinop^\perp$ and this operator possesses a bounded right inverse by Proposition \ref{thm:seclinest}.  The Banach space inverse function theorem can thus be applied to the equation $\pi \circ \secdifop (f) = 0$ provided that the three fundamental estimates \eqref{eqn:iftestone}, \eqref{eqn:iftesttwo} and \eqref{eqn:iftestthree} described in Section \ref{sec:jacobi} can be established.  The construction of the right inverse and its bound in Proposition \ref{thm:seclinest} constitutes the first of these estimates.  The second and third estimates can be proved by very straightforward modifications of the analogous estimates for Main Theorem 1, namely Proposition \ref{prop:error} and Proposition \ref{prop:nonlin}. The proof then follows as before.
\end{proof}

\subsection{The Balancing Problem and the Conclusion of the Proof} 
\mylabel{subsec:balmap}

\paragraph*{The balancing map.} Proposition \ref{prop:secperturb} shows that the equation $\secdifop (f) = 0$ can be solved up to a finite dimensional error term; i.e.~a function $f_{\alpha, \tau, \sigma} \in \ctbg(X^{\perp})$ can be found so that only the $L^2$-projection of $\secdifop( f_{\alpha, \tau, \sigma})$ to the subspace $\tilde{ \mathcal K}_{\mathit{sym}}$ fails to vanish identically as a function of $\sigma$.  It will now be shown, however, that there is a special choice of $\sigma$ that yields a so-called \emph{balanced} initial configuration $\Lambda_{\alpha, \tau, \sigma}^{\#}$ for which $\secdifop(f_{\alpha, \tau, \sigma})$ vanishes completely.  Therefore the solution $f_{\alpha, \tau, \sigma}$ for this choice of $\sigma$ yields the desired deformation of $\secapprox$ into an exactly CMC hypersurface.

In order to derive the \emph{balancing condition} that will be used determine this special value of $\sigma$, one must understand in greater detail the relationship between $\sigma$ and the quantity $(\mathrm{id} - \pi) \circ  \secdifop (f_{\alpha, \tau, \sigma})$ where $\pi$ is the $L^2$-projection onto $X^{\perp}$.  To this end, introduce the functions $\tilde q_k^{\, \ast} : \secapprox \rightarrow \R$ which are obtained from the functions in Definition \ref{defn:apkerbasis} by replacing the cut-off functions $\chi_{\mathit{ext}}^{s,k}$ by the cut-off functions $\eta_{\mathit{ext}}^{s,k}$ used in the proof of Proposition \ref{thm:seclinest}.

\begin{defn}
	\mylabel{defn:balancingmap}
	The \emph{balancing map} of $\secapprox$ is the function $B_{\alpha, \tau} : \R^{(N-2)/4} \rightarrow \R^{(N-2)/4}$ given by 
	\begin{equation}
		\mylabel{eqn:balmap}
		B_{\alpha, \tau} (\sigma): = \left( \int_{\secapprox} \secdifop (f_{\alpha, \tau, \sigma}) \cdot \tilde q_1^{\, \ast} \, , \ldots,  \int_{\secapprox} \secdifop (f_{\alpha, \tau, \sigma}) \cdot \tilde q_{(N-2)/4}^{\, \ast} \right), 
	\end{equation}
	where $f_{\alpha, \tau, \sigma} \in \ctbg(X^{\perp})$ is the solution found in Proposition \ref{prop:secperturb}. Note that the dimension of the domain of $B_{\alpha, \tau}$ is accurate because of the conditions \eqref{eqn:displacement} that are assumed to be valid for $\sigma$. 
\end{defn}

\noindent In terms of the balancing map, what remains to be done in order to prove Main Theorem \ref{result2} is to find a value of $\sigma$ for which $B_{\alpha, \tau}(\sigma) = 0$.  This is because the identity $\int_{\secapprox} \tilde q_j \cdot \tilde q_k^{\, \ast} =  \delta_{jk} + \mathcal O(\rho^n)$ implies that $(\mathrm{id} - \pi) \circ \secdifop (f_{\alpha, \tau, \sigma}) = 0$ as well.

The ordinary inverse function theorem for smooth functions on Euclidean space will be used to find this value of $\sigma$ in the following way.  First, $B_{\alpha, \tau}$ will be approximated by a simpler function $\mathring B_{\alpha, \tau} : \R^{(N-2)/4} \rightarrow \R^{(N-2)/4}$ depending only on the geometry of $\secapprox$ (i.e.~independent of $f_{\alpha, \tau, \sigma}$) using the ideas behind the Korevaar-Kusner-Solomon \cite{kks} and Kapouleas \cite{kapouleas7} balancing formula.   There will be a fairly obvious choice of $\sigma$ for which $\mathring B_{\alpha, \tau}(\sigma) = 0$; consequently it will be true that $B_{\alpha, \tau}(\sigma) \approx 0$ for this choice of $\sigma$.  Then by the inverse function theorem, it will be possible to find a small perturbation of $\sigma$ for which $B_{\alpha, \tau}(\sigma) = 0$ exactly since it will turn out that the derivative $\Dif B_{\alpha, \tau}(\sigma)$ is invertible with large norm.

\paragraph*{The proof of Main Theorem 2.}  The first task is to approximate the balancing map of $\secapprox$ by an easily-understood function that is independent of $f_{\alpha, \tau, \sigma}$.  Let $\bar \eps_k$ denote the scale parameter of the $k^{\mathit{th}}$ neck region of $\secapprox$.
\begin{prop}
	\mylabel{prop:balmapapprox}
	The balancing map of $\approxsol$ can be approximated as follows whenever $\tau$ is sufficiently small.  There exists  a function $\tilde B_{\alpha, \tau} : \R^{(N-2)/4} \times \R^{(N-2)/4}$ satisfying the estimate $$\Vert \tilde B_{\alpha, \tau} (\sigma) \Vert_{C^2} \leq C \rho_\eps^n | f_{\alpha, \tau, \sigma} |_{\ctbg}$$ for some constant $C$ independent of $\eps$ so that 
	$$B_{\alpha, \tau} (\sigma) = \mathring B_{\alpha, \tau} (\sigma)  + \tilde B_{\alpha, \tau}(\sigma)$$ 
	where
	\begin{equation}
		\mylabel{eqn:balmapapprox}
		\mathring B_{\alpha, \tau} (\sigma) = \omega \left(
		\begin{array}{cccccc}
			-1 & 1 &  &  & & \\
			 & -1 & 1 &  & & \\
			& & & \ddots \\
			& & & & -1 & 1 \\
		\end{array}
		\right) \cdot \left(
		\begin{matrix} 
			\bar \eps_0^{\, n-1} \\[1.5ex] \vdots \\[1.5ex] \bar \eps_{(N-2)/4}^{\, n-1}
		\end{matrix}
		\right) \, .
	\end{equation}
	Here $\omega$ is another constant independent of $\eps$.
\end{prop}

\begin{proof}
Let $[B_{\alpha, \tau}(\sigma)]_k$ denote the $k^{\mathit{th}}$ component of the balancing map.  By symmetry,  it is sufficient for the analysis of  $[B_{\alpha, \tau}(\sigma)]_k$ to examine only the $k^{\mathit{th}}$ spherical region along the geodesic $\gamma_1$ and its two adjoining neck and transition regions, namely the region
$$\mathcal S_{\eps}^k := \mathcal T_ {\eps}^{1,k,+} \cup{\mathcal N}_ {\eps}^{1,k,+} \cup \mathcal E_ {\eps}^k \cup \mathcal T_ {\eps}^{1,k+1,-} \cup{\mathcal N}_ {\eps}^{1,k+1,-}  \, .$$
Next, let
$$D_k := R_{1, (\sigma_k + \sigma_{k+1})/2} \circ R_{1, \alpha + \tau/2}^{2k+1} \circ K^{-1} \big( \{ (\bar b_k, \hat y) \: : \: \Vert\hat y \Vert \leq \eps_k \} \big)$$
be the central cross-sectional $n$-sphere of $\mathcal N_{\eps}^{1,k}$ and let $c_k := \partial D_k$ so that $\partial \mathcal S_{\eps}^k = c_k \cup c_{k+1}$.   Let $\phi_{\alpha, \tau, \sigma}$ denote the normal deformation of $\secapprox$ obtained from the solution $f_{\alpha, \tau, \sigma}$.  Using symmetry and simple estimates, one finds 
\begin{align*}
	[B_{\alpha, \tau}(\sigma)]_k &= 2 \int_{\mathcal \phi_{a, \tau, \sigma}(\secapprox)} \bar \eta^{1,k}_{\mathit{ext}} \cdot \big( H(\phi_{a,\tau, \sigma}(\secapprox)) - H_\alpha \big) \cdot \left( x^1 \circ R^{-1}_{1, k(2 \alpha + \tau) + \sigma_k} \right) \\
	&= 2 \int_{\phi_{a, \tau, \sigma}( \mathcal S_{\eps}^{1,k})} \big( H(\phi_{a,\tau, \sigma}(S_{\eps}^{1,k})) - H_\alpha \big) \cdot \left( x^1 \circ R^{-1}_{1, k(2 \alpha + \tau) + \sigma_k} \right) + \mathcal O(\rho_{\eps}^{n} |f_{\alpha, \tau, \sigma} |_{\ctbg} )
\end{align*}
since $\eta_{\mathit{ext}}^{1,k}$ is supported outside the neighbourhood $\mathit{Tub}_{\rho_{\eps}}(\gamma_1)$.  Here, $\mathcal O(\rho_{\eps}^{n} |f_{\alpha, \tau, \sigma} |_{\ctbg} )$ denotes a function bounded in $C^2$ norm by $\rho_{\eps}^{n} |f_{\alpha, \tau, \sigma} |_{\ctbg} $. Now recall the Korevaar-Kusner-Solomon identity given in equation \eqref{eqn:balancing} satisfied by the integral of the mean curvature of hypersurface against a Jacobi field.  That is, if $\nu_k$ denotes the outward unit normal vector field of $c_k$ parallel to $N_{\eps}^{1,k}$ and $V_k$ denotes the compatibly-oriented unit normal vector field of $D_k$ while $Y_k$ denotes the vector field 
$$Y_k := \big( R_{1, (\sigma_k + \sigma_{k+1})/2} \circ R_{\alpha+ \tau/2}^{2k+1} \big)_\ast \circ \left(x^0  \frac{\partial}{\partial x^1} - x^1 \frac{\partial}{\partial x^0} \right) \circ R_{\alpha+ \tau/2}^{-2k-1} \circ R_{1, (\sigma_k + \sigma_{k+1})/2}^{-1} \, ,$$
then one finds  
\begin{align}
	\mylabel{eqn:balmapone}
	\frac{1}{2} [B_{\alpha, \tau}(\sigma)]_k &=  \left(  \int_{\phi_{a, \tau, \sigma}(c_k)} \langle \nu_k , Y_k \rangle  -   H_\alpha \int_{\phi_{a, \tau, \sigma}( D_k) } \langle V_k, Y_k \rangle  \right) \notag \\
 	&\qquad -  \left(  \int_{\phi_{a, \tau, \sigma}(c_{k+1})} \langle \nu_{k+1} , Y_k \rangle  -   H_\alpha \int_{\phi_{a, \tau, \sigma}( D_{k+1}) } \langle V_{k+1}, Y_k \rangle  \right) + \mathcal O \big(\rho_{\eps}^{n} |f_{\alpha, \tau, \sigma} |_{\ctbg} \big)  \notag \\
	&= \int_{c_k} \langle \nu_k , Y_k \rangle  -   \int_{c_{k+1}} \langle \nu_{k+1} , Y_k \rangle + \mathcal O \big(\rho_{\eps}^{n} |f_{\alpha, \tau, \sigma} |_{\ctbg} \big) \, .
\end{align}
This is because the integrals over $\phi_{a, \tau, \sigma}(c_\ast)$ can be replaced by integrals over $c_\ast$ and the integrals over $\phi_{a, \tau, \sigma}(D_\ast)$ can be dropped altogether at the expense of further error terms whose size is strictly smaller than $\mathcal O \big(\rho_{\eps}^{n} |f_{\alpha, \tau, \sigma} |_{\ctbg} \big)$.  Finally, the calculation of the integrals $\int_{c_\ast} \langle \nu_\ast , Y_k \rangle$ in \eqref{eqn:balmapone}  can be carried out in the stereographic coordinate chart used to define the neck regions, in which case $Y_k$ simply becomes the translation vector field perpendicular to the axis of the neck.  This calculation is very straightforward, and yields a quantity proportional to the $(n-1)$ dimensional area of $c_k$, thus proportional to $\bar \eps_k^{\, n-1}$.
\end{proof}

Recall that $\bar \eps_k$ depends on the separation between the hyperspheres through the procedure culminating in formula \eqref{eqn:epschoice}.  Thus one can see that there are universal functions $E : \R \rightarrow \R$ and $C: \R \rightarrow \R$ so that $$\bar \eps_k(\sigma) := E \bigl( C(\tau + \sigma_{1} - \sigma_{0}) + \cdots + C(\tau + \sigma_{0} - \sigma_{N-1}) \bigr) \cdot C(\tau + \sigma_{k+1} - \sigma_{k}) \, .$$  It is thus clear how to choose $\sigma$ so that $\mathring B_{\alpha, \tau}$ vanishes: one chooses $\sigma_1 = \cdots = \sigma_{N-1} = 0$ which corresponds to an equal spacing of the hyperspheres making up $\secapprox$.  By Proposition \ref{prop:balmapapprox}, one has $B_{\alpha, \tau}(0) = \mathcal O \big(\rho_\eps^n |f_{\alpha, \tau, \sigma} |_{\ctbg} \big)$.  Note that when this choice of $\sigma$ is made, then all the various quantities on the different constituents of $\tilde \Lambda_{\alpha, \tau, 0}$ are equal and determined uniquely by the choice of $\tau$ through the scale parameter $\eps := \eps(\tau)$.

It remains to show that the derivative $\Dif B_{\alpha, \tau}(0)$ is invertible and is bounded below by a quantity much larger than $\mathcal O \big(\rho_\eps^n |f_{\alpha, \tau, \sigma}|_{\ctbg} \big)$.  Since $\Dif B_{\alpha, \tau}(0) = \Dif \mathring B_{\alpha, \tau}(0) + \mathcal O \big(\rho_\eps^n |f_{\alpha, \tau, \sigma}|_{\ctbg} \big)$, it is sufficient to establish this fact for $\Dif \mathring B_{\alpha, \tau}(0)$.

\begin{prop}
	\mylabel{prop:invert}
	The derivative matrix of $\mathring B_{\alpha, \tau}$ at $\sigma = 0$ is
	$$\Dif \mathring B_{\alpha, \tau}(0) =\omega(\eps) \left(
	\begin{array}{cccccccc}
		0 & 1  \\
		& 0 & 1  \\
		&  & 0 & 1 \\
		& & & & \ddots \\
		& & & & & & 0 & 1 \\
		1& & & & & & & 0
	\end{array}
	\right)$$
	where $\omega(\eps) = \mathcal O(\eps^{n-1})$.  This is an invertible linear transformation.  
\end{prop}

\begin{proof}
The matrix appearing in $\Dif \mathring B_{\alpha, \tau}(0)$ is easily found by applying the chain rule to $\mathring B_{\alpha, \tau}$ and the functions $\bar \eps_k^{\, n-1} := [\eps_k(\sigma)]^{n-1}$ at $\sigma =0$.
\end{proof}

The consequence of Proposition \ref{prop:invert} is that the derivative $\Dif B_{\alpha, \tau}(0)$ is invertible and much larger than $B_{\alpha, \tau}(0)$.  Consequently the inverse function theorem can be applied to find $\sigma$ near zero for which $B_{\alpha, \tau}(\sigma) = 0$.  The proof of Main Theorem 2 is now complete. \hfill \qedsymbol

\section{Attaching Delaunay-Like Handles to Clifford Tori}

\subsection{Introduction}

The techniques that have been developed so far allow the analogue of one more `classical' problem from the theory of constant mean curvature hypersurfaces in Euclidean space to be realized in $\Sph^{n+1}$.  That is, one can attempt to attach a number of Delaunay-like handles to a given CMC hypersurface and hope to perturb the resulting hypersurface to have constant mean curvature.  It should be possible to pose conditions, in very general terms, for when such a construction could be realized.  However it will be simplest to focus on one canonical example of the construction since this is easier for the reader to visualize and it is in keeping with the spirit of this paper.  Moreover, the generalization to more complicated settings will be readily believable once a simple example is shown to work.

The case that will be considered here is the attachment of a truncated Delaunay-like hypersurface in $\Sph^{n+1}$ to the generalized Clifford torus $\cliff$ of type $\Sph^{n_1} \times \Sph^{n_2}$ where $n_1 + n_2 = n$.  From what has been learned about the problem of gluing two Delaunay-like hypersurfaces together, it is intuitively clear that one must find a geodesic $\gamma$ that makes contact orthogonally with $\cliff$ in two places, and that hyperspheres of the same mean curvature as $\cliff$ must be positioned with equal spacing along $\gamma$.  The hyperspheres will then be in balanced position.  What is not clear is that this can be done with enough symmetry to rule out all the Jacobi fields of $\cliff$.  However, by analyzing these Jacobi fields carefully, it will be found that it is indeed possible rule out these Jacobi fields and thereby realize the handle-attachment procedure.

\subsection{The Approximate Handle-Attachment}
\label{subsec:handleassembly}

\paragraph*{The Generalized Clifford Torus in \boldmath $\Sph^{n+1}$.}  The generalized Clifford tori of type $\Sph^{n_1} \times \Sph^{n_2}$ with $n_1 + n_2 = n$ in $\Sph^{n+1}$ are defined by
$$\cliff := \big\{ (x_1, x_2) \in \R^{n_1+1} \times \R^{n_2+1} : \| x_1 \| = \cos(\alpha) \mbox{ and } \| x_2 \| = \sin(\alpha) \big\}$$
where  $\alpha \in (0, \pi/2)$.  Then $\cliff$ can be parametrized by $(\Theta^1, \Theta^2) \longmapsto ( \cos(\alpha) \Theta^1, \sin(\alpha) \Theta^2)$ where the map $\Theta^j : \Sph^{n_j} \rightarrow \R^{n_j+1}$ is an embedding of $\Sph^{n_j}$ as the unit sphere in $\R^{n_j+1}$.  The calculations of \cite{mepacard2} give the relevant facts about the geometry of the Clifford tori.  In particular, the mean curvature of $\cliff$ is $H^{\mathit{Cliff}}_\alpha := n_2 \cot(\alpha) - n_1 \tan (\alpha)$ and the linearized mean curvature operator is given by
$$\mathcal L_{n_1, n_2, \alpha} := \cos^{-2} (\alpha) \big(\Delta_{1} + n_1 \big) + \sin^{-2} (\alpha) \big(\Delta_{2} + n_2 \big) $$
where $\Delta_{1}$ and $\Delta_{2}$ are the Laplacians of $\Sph^{n_1}$ and $\Sph^{n_2}$, respectively. 
One deduces that the Jacobi fields of $\cliff$ are the products of the coordinate functions --- i.e.~$\mathcal L_{n_1, n_2, \alpha} \big( x_1^{t_1} x_2^{t_2}  \big|_{\cliff} \big) = 0$ for all $t_j = 0, \ldots, n_j$.  

\paragraph*{Assembly of the approximate solution.}  To begin assembling the approximate solution for the handle-attachment process, choose a canonical point $p \in \cliff$ defined by $p := (\cos(\alpha) P_1, \sin(\alpha) P_2)$ where $P_j \in \R^{n_j+1}$ is the point $P_j := (1, 0, \ldots, 0)$.   For any $\theta \in [0, 2 \pi]$, define the rotation 
\begin{equation}
	\label{eqn:torusgeo}
	R_\theta := \left(
	\begin{array}{c|ccc|c|ccc}
		\cos(\theta) & 0 & \cdots & 0 & -\sin(\theta) & 0 & \cdots & 0\\
		\hline 0 & & & & 0 & & &  \\
		\vdots & & I_{n_1} & & \vdots & & 0 & \\
		0 & & & & 0 & & &  \\
		\hline \sin(\theta) & 0 & \cdots & 0 & \cos(\theta) & 0 & \cdots & 0 \\
		\hline 0 & & & & 0 & & &  \\
		\vdots & & 0 & & \vdots & & I_{n_2} & \\
		0 & & & & 0 & & &
	\end{array} \right) \, .
\end{equation}
The normal vector of $\cliff$ at $p$ is $N_p := -\sin (\alpha) P_1 + \cos(\alpha) P_2$.  Let $\gamma_p$ be the normal geodesic emanating from $\cliff$ at the point $p$; then $\gamma_p$ can be parametrized by $t \longmapsto R_{t} (p)$ for $\gamma_p$.  Define the point $p' := (-\cos(\alpha) p_1, \sin(\alpha) p_2)$.  Then $$\gamma_p \cap \cliff = \{ p, p', -p, -p'\}$$ and these points occur at the parameter values $t = 0$, $t= \pi - 2 \alpha $, $t = \pi$ and $t= 2 \pi - 2 \alpha $, respectively.

Recall that the mean curvature of a hypersphere $S_{\hat \alpha}$ with parameter $\hat \alpha \in (0, \pi/2)$ is $H^{\mathit{Sph}}_{\hat \alpha} := n \cot (\hat \alpha)$.  Therefore the  hypersphere with the same mean curvature (up to sign) as $\cliff$ has $\hat \alpha  := \arccot \big| \frac{n_2}{n} \cot (\alpha) - \frac{n_1}{n} \tan (\alpha) \big|$.  Denote the common value of the mean curvature simply by $H_\alpha$.  Since the length of the geodesic segment from $p$ to $-p'$ is exactly $2 \sigma$, it is thus possible to position $N$ rotated  hyperspheres of parameter $\hat \alpha$ along $\gamma_p$ from $p$ to $-p'$, where each hypersphere is separated from its two nearest neighbours by a distance $\tau$, when $2 N \hat \alpha + (N+1) \tau = 2 \alpha + 2 m \pi$ for some integer $m \geq 0$.  These spheres wind $m$ times around $\gamma_p$ before ending up a distance of exactly $\tau$ from $-p'$.

\begin{rmk}
It is clear that for every fixed $\alpha \in (0, \pi/2)$ and for every $\tau_0>0$ there is always \emph{some} value of $\tau$ satisfying $0 < \tau < \tau_0$ and \emph{some} integers $N$ and $m$ so that the equation $2 N \hat \alpha + (N+1) \tau = 2 \alpha + 2 m \pi$ holds.     It is not clear if this equation can be made to hold in this way for fixed $m$ as well (e.g.~for $m=0$ which would correspond to embedded configurations).   Since $\lim_{\alpha \rightarrow \pi/2} \frac{\alpha}{ \hat \alpha(\alpha)} = \infty$ one can find $N$ and $\alpha_0$ so that $N \hat \alpha(\alpha_0)  = \alpha_0$. Then since the derivative $\frac{\dif \hat \alpha}{\dif \alpha}$ is bounded away from zero, one can determine $\alpha := \alpha(\tau)$ near $\alpha_0$ so that $2 N \hat \alpha(\alpha(\tau)) + (N+1) \tau = 2 \alpha(\tau)$ for all sufficiently small $\tau$. Moreover, if $|\alpha_0 - \pi/2| > T$ then one can arrange to have $|\alpha(\tau) - \pi/2| > T/2$.  This ensures that the mean curvature of the configuration, though large since $\alpha_0$ is near $\pi/2$, never becomes unbounded.
\end{rmk}

These considerations lead to the definition of the initial configuration of  hyperspheres and Clifford torus that will be glued together below.   For any $\theta \in [0, 2 \pi]$ first introduce the matrix $Q_\alpha \in SO(n+2)$ by
\begin{equation}
	\label{eqn:torusrotn}
	Q_\alpha := \left( 
	\begin{array}{cc|ccc|ccc}
		\cos(\alpha) & -\sin(\alpha) & 0 & \cdots & 0 & 0 & \cdots & 0 \\
		\hline 0 & 0 & & &&&& \\
		\vdots& \vdots & &I_1& & & 0 & \\
		0 & 0 & & &&&& \\
		\hline \sin(\alpha) & \cos(\alpha) & 0 & \cdots & 0 & 0 & \cdots & 0 \\ 
		\hline 0 & 0 & & &&&& \\
		\vdots& \vdots & &0& & & I_2 & \\
		0 & 0 & & &&&&
	\end{array} \right)
\end{equation}
where $I_j$ is the $n_j \times n_j$ identity matrix.  The matrix $Q_\alpha$ takes the point $p_0 := (1, 0, \ldots, 0) \in \R^{n+2}$ to the point $p$, the vector $N_0 := (0, 1, 0, \ldots 0) \in \R^{n+2}$ to the vector $N_p$, and the subspace $\{ 0, 0 \} \times \R^n$ to the tangent space of $\cliff$ at $p$.   Note that $R_\theta \circ Q_\alpha = Q_{\theta + \alpha}$.

The rotation $R_\theta \circ Q_ \alpha$ takes the  hypersphere $S_{\hat \alpha}$ to another  hypersphere centered on the point $R_\theta(p)$ along the geodesic $\gamma_p$. Let $N, m, \alpha, \tau$ and $\hat \alpha$ satisfy $2 N \hat \alpha + (N+1) \tau = 2 \alpha + 2 m \pi$ and introduce displacement parameters $\sigma_1, \ldots, \sigma_N$ as before.  Define the initial configuration
\begin{equation}
	\mylabel{eqn:handleattach}
	\Lambda^{\#}_{\alpha, \tau, \sigma} := \left[ \bigcup_{k=0}^{N-1} R^{-1}_{\sigma_k} \circ R^{-1}_{k(2 \alpha + \tau) + \alpha + \tau} \circ Q_\alpha (S_{\hat \alpha}) \right] \cup \cliff \, .
\end{equation}
The mean curvature of each component of $\Lambda^{\#}_{\alpha, \tau, \sigma}$ is $H_{\alpha}$.

The next step in the assembly of the approximate solution is the perturbation of each of the components of $\Lambda ^{\#}_{\alpha, \tau, \sigma}$ in order to improve their fit with the catenoidal necks that are needed to glue everything together.  By analogy with the procedure described in Section \ref{subsec:newassembly}, each hypersphere $R^{-1}_{k(2 \alpha + \tau) + \alpha + \tau + \sigma_k} \circ Q_\alpha (S_{\hat \alpha})$ should be replaced by $R^{-1}_{k(2 \alpha + \tau) + \alpha + \tau + \sigma_k} \circ Q_\alpha (\tilde S_{\hat \alpha, \eps_k, b_k})$ where $\tilde S_{\hat \alpha, \eps_k, b_k}$ is the perturbed hypersphere with scale parameter $\eps_k$ and translation parameter $b_k$.  A similar perturbation must also be found for the Clifford torus $\cliff$.  This is provided by the construction used in Butscher and Pacard's papers \cite{mepacard1,mepacard2} where it was shown how to perturb the Clifford torus in the normal direction using the singular solution of the equation $\mathcal L_{n_1, n_2, \alpha} (G) = 0$ on $\cliff \setminus \{ p, -p' \}$.   This perturbation is completely analogous to the perturbation of the hyperspheres developed in this paper, though it differs in such details as the numerical values of the coefficients of its asymptotic expansions near $p$ and $-p'$.  Denote the perturbed Clifford torus by $\tilde C_{\alpha, \eps}^{n_1, n_2}$.   Then as before, one can find the parameters $\eps_k, b_k$, and $\eps$ along with neck scale and translation parameters $\bar \eps_k$ and $\bar b_k$ that ensure optimal matching between these hypersurfaces and small catenoidal necks with the parameters $\bar \eps_k, \bar b_k$.   Moreover, these are all determined uniquely in terms of $\tau$.

\begin{defn}
	\mylabel{defn:handleregions}
	Recall the terminology from Sections \ref{subsec:assembly} and \ref{subsec:newassembly}. 
	\begin{itemize}
		\item Define the neck regions
		\begin{align*}
			\mathcal N^{k, \pm}_{\eps} &:= R^{-1}_{(\sigma_k+ \sigma_{k-1})/2} \circ R^{-1}_{(k+\frac{1}{2})(2\alpha + \tau) + \alpha + \tau} \circ Q_\alpha \circ K^{-1}  \big( \tilde \Sigma_{\bar \eps_k, \bar b_k}^\pm \cap \mathit{Cyl}(\rho_{\eps}/2) \big) \, .
		\end{align*}
	Set $\mathcal N_{\eps}^{k} := \mathcal N_{\eps}^{k,+} \cup \mathcal N_{\eps}^{k,-} $.  These necks are oriented so that $\mathcal N_{\eps}^{k,+}$ lies behind $\mathcal N_{\eps}^{k,-}$ along the geodesic $\gamma$ and that the $k^{\mathit{th}}$ neck is centered on the point $R^{-1}_{(\sigma_k + \sigma_{k+1})/2} \circ R^{-1}_{(k+\frac{1}{2})(2\alpha + \tau) + \alpha + \tau}(p)$.  
	
		\item Define the transition regions
		\begin{align*}
			\mathcal T^{k, \pm}_{\eps} &:= R^{-1}_{ \sigma_k} \circ R^{ -1}_{(k+\frac{1}{2})(2\alpha + \tau) + \alpha + \tau}\circ Q_\alpha \circ K^{-1}  \big( \tilde \Sigma_{\bar \eps_k, \bar b_k}^\pm \cap \big[ \mathit{Cyl}(2\rho_\eps) \setminus \mathit{Cyl}(\rho_\eps/2) \big] \big) \, .
		\end{align*}
	Set $\mathcal T_{\eps}^{k} := \mathcal T_{\eps}^{k,+} \cup \mathcal T_{\eps}^{k,-} $.  	
	
		\item Define the exterior regions
		\begin{align*}
			\mathcal E_{\alpha, \eps}^{\mathit{Cliff}} &:= \tilde C_{\alpha, {\eps_0}}^{n_1, n_2} \setminus \big[ B_{\tilde \rho_{\eps_0}} (p) \cup B_{\tilde \rho_{\eps_0}} (-p') \big] \\[0.5ex]
			\mathcal E_{\hat \alpha, \eps}^{k} &:= R^{-1}_{\sigma_k} \circ R^{-1}_{k(2 \alpha + \tau) + \alpha + \tau} \circ Q_\alpha \big( \tilde S_{\alpha, \eps_k, b_k} \setminus \big[ B_{\tilde \rho_{\eps}} (p^+) \cup B_{\tilde \rho_{\eps}} (p^-)\big] \big) \, .
		\end{align*} 
	\end{itemize}
\end{defn}

\noindent The approximate solution corresponding to the configuration $\Lambda _{\alpha,\tau,\sigma}^{\#}$ can now be defined as follows.
\begin{defn}
	\mylabel{defn:handleapproxsol}
	The \emph{approximate solution} with parameters $\alpha, \tau$ and $\sigma$ is the hypersurface 
	\begin{equation}
		\mylabel{eqn:handleapproxsol}
		\handleapprox := \mathcal E_{\alpha, \eps}^{\mathit{Cliff}} \cup \left[\bigcup_{k =0}^{ N-1} \Big(  \mathcal E_{\hat \alpha, \eps}^{k} \cup  \mathcal T_{\eps}^{ k} \cup \mathcal N_{\eps}^{k}\Big) \right] \cup  \big[ \mathcal T_{\eps}^{ N} \cup \mathcal N_{\eps}^{N} \big]\, .
	\end{equation}	
\end{defn}

\subsection{Symmetries and Jacobi Fields}

\paragraph*{Symmetries.} The generalized Clifford tori are highly symmetric objects in $\Sph^{n+1}$, being invariant under the diagonal $O(n_1+1) \times O(n_2+1)$ in $O(n+2)$.  Once the Delaunay-like handle is attached, there is considerably less symmetry.  The symmetries which preserve \emph{both} the perturbed  generalized Clifford torus and the perturbed hyperspheres are as follows.

\medskip 
\noindent 1. \itshape Rotational symmetries. \upshape 
\medskip

Define the transformation $S^{00}_{B_1, B_2} \in O(n+2)$ by choosing $B_j \in O(n_j+1)$ and then setting
\begin{equation}
	\label{eqn:fixsubgp}
	S_{B_1, B_2}^{00} := \left(
	\begin{array}{c|ccc|c|ccc}
		1& 0 & \cdots & 0 & 0 & 0 & \cdots & 0\\
		\hline 0 & & & & 0 & & &  \\
		\vdots & & B_{1} & & \vdots & & 0 & \\
		0 & & & & 0 & & &  \\
		\hline 0 & 0 & \cdots & 0 & 1 & 0 & \cdots & 0 \\
		\hline 0 & & & & 0 & & &  \\
		\vdots & & 0 & & \vdots & & B_{2} & \\
		0 & & & & 0 & & &
	\end{array}
\right)
\end{equation}
which is the transformation that preserves $\cliff$ and the geodesic $\gamma$ as well as the perturbed Clifford torus.   These transformations also preserve the perturbed hyperspheres positioned along $\gamma$ and the necks between them because of the commutation relations $S_{B_1, B_2}^{00} \circ R_\theta = R_\theta \circ S_{B_1, B_2}^{00}$ and  
$$Q_\alpha^{-1} \circ S_{B_1, B_2}^{00} \circ Q_\alpha = \left(
\begin{array}{cc|cc}
1 & & & \\
& 1 & & \\
\hline & & B_1 & \\
& & & B_2
\end{array}
\right)
$$
which is a transformation that preserves the perturbed hypersphere $\tilde S_{\alpha, \eps}$.

\medskip 
\noindent 2. \itshape Reflection symmetry. \upshape 
\medskip

One additional symmetry of $\handleapprox$ comes from the reflection that reverses the geodesic $\gamma_p$ with respect to the point half-way between $p$ and $-p'$, namely the point $R_{\pi - \alpha}(p)$.  This is the transformation $T \in O(n+2)$ given by 
$$T := \left(
\begin{array}{c|ccc|c|ccc}
	1& 0 & \cdots & 0 & 0 & 0 & \cdots & 0\\
	\hline 0 & & & & 0 & & &  \\
	\vdots & & I_{n_1} & & \vdots & & 0 & \\
	0 & & & & 0 & & &  \\
	\hline 0 & 0 & \cdots & 0 & -1 & 0 & \cdots & 0 \\
	\hline 0 & & & & 0 & & &  \\
	\vdots & & 0 & & \vdots & & I_{n_2} & \\
	0 & & & & 0 & & &
\end{array}
\right) \, .
$$
This reflection preserves the Clifford torus and its perturbation since it is an element of its symmetry group.  If the displacement parameters satisfy $\sigma_k = \sigma_{N-k-1}$, then $T$ also preserves the perturbed hyperspheres and the necks because of the commutation relations $T \circ R_\theta \circ T = R^{-1}_\theta$ and 
$$Q_{-\alpha}^{-1} \circ T \circ Q_\alpha = \left(
\begin{array}{cc|cc}
1 & & & \\
& -1 & & \\
\hline & & I_{n_1} & \\
& & & I_{n_2}
\end{array}
\right)
$$
which is a transformation that preserves the perturbed hypersphere $\tilde S_{\alpha, \eps}$.

\paragraph*{Jacobi fields.} It remains to find all the Jacobi fields of $\handleapprox$ that are invariant under the transformations $S^{00}_{B_1, B_2}$ and $T$.  As in Section \ref{sec:secanalysis}, one proceeds as follows.  First define smooth and monotone cut-off functions $\chi_k : \handleapprox \rightarrow \R$ for $k = 0, \ldots, N-1$ supported in $\mathcal E^k_{\eps}$ and another such function $\chi_{\mathit{Cliff}} : \handleapprox \rightarrow \R$ supported in $\mathcal E^{\mathit{Cliff}}_{\alpha, \eps}$.  Consider $\tilde S_{\hat \alpha, \eps}$ as a hypersurface in $\Sph^{n+1} \subseteq \R^{n+2}$ endowed with coordinates $(x^0, x^1 , \ldots, x^{n+1})$ and define the following functions:
\begin{subequations}
	\mylabel{eqn:handleapproxjac}
	\begin{equation}
		\tilde q_{k}^{\, t}  := 
			\chi_k \cdot \left( x^t  \big\vert_{\tilde S_{\hat \alpha, \eps} } \circ Q_\alpha^{-1} \circ R_{k(2 \alpha + \tau) + \alpha + \tau} \circ R_{\sigma_k} \right) 
	\end{equation}
	Also, in the standard coordinates $(x_1, x_2)$ for $\R^{n_1+1} \times \R^{n_2+1}$ that are in use at the moment, define the functions
	\begin{equation}
		\tilde q^{\, t_1, t_2}_{\mathit{Cliff}}  := \chi_{\mathit{Cliff}} \cdot \left( x_1^{t_1} x_2^{t_2}  \big\vert_{\tilde C^{n_1, n_2}_{\alpha, \eps}} \right)  \, .
	\end{equation}
\end{subequations}
Then the vector space of \emph{all} approximate Jacobi fields in $\handleapprox$ is the set of functions 
$$\tilde{ \mathcal K} := \linspan_\R \{ \tilde q_{k}^{\, t} : t = 1, \ldots, n \; \mbox{ and all } \; k \} \cup \{ \tilde q^{\, t_1, t_2}_{\mathit{Cliff}} : t_1 = 1, \ldots, n_1 \mbox{ and } t_2 = 1, \ldots, n_2  \} \, .$$
It is easy to see that the subspace of these functions which are invariant under the symmetries  $S_{B_1, B_2}^{00}$ for all $B_j \in O(n_j+1)$ is spanned by $\{\tilde q_{k}^{1} :  \mbox{ all }\; k \} \cup \{ \tilde q^{00}_{\mathit{Cliff}} \}$.   Under the further requirement of invariance with respect to $T$, the function $\tilde q^{\, 00}_{\mathit{Cliff}}$ is eliminated and one is left with the following subspace.

\begin{defn}
	The subspace of approximate Jacobi fields of  $C^{2,\beta}_\delta (\handleapprox)$ invariant under the group of symmetries generated by  $T$ and $S_{B_1, B_2}^{00}$ for all $B_j \in O(n_j+1)$  is the subspace of functions 
	\begin{equation}
		\mylabel{eqn:handleapker}
		\tilde{\mathcal K}_{\mathit{sym}} := \linspan_\R \left\{ \tilde q^1_k - \tilde q^1_{N-k-1} : \mbox{all } k
		\right\} \, .
	\end{equation}
	Note that the dimension of $\mathcal K_{\mathit{sym}}$ is $N$.
\end{defn}

\subsection{The Proof of Main Theorem 3}

The proof of Main Theorem 3 is now in all respects identical to the proof of Main Theorem 2, for the following reasons.  One can set up normal deformations of $\handleapprox$ using functions invariant under $T$ and $S^{00}_{B_1,B_2}$ for all $B_j \in O(n_j+1)$.  Then the construction of the right inverse of the linearized mean curvature operator and its estimate follow as in Proposition \ref{thm:seclinest} because the nature of the Jacobi fields of $\handleapprox$ is unchanged; namely, there is one invariant Jacobi field on each perturbed hypersphere and none on the necks and on the perturbed Clifford hypersphere by the nature of the symmetries and the choice of weighted function space.  This gives the linear estimate \eqref{eqn:iftestone} needed to invoke the inverse function theorem.  Then the estimates of the mean curvature of the Delaunay-like hypersurfaces and the estimates of the mean curvature of the generalized Clifford tori from \cite{mepacard1,mepacard2} can be combined to yield the first of the non-linear estimates \eqref{eqn:iftesttwo} needed to invoke the inverse function theorem.  Finally, the second of the non-linear estimates holds as before.  Hence when $\tau$ is sufficiently small the inverse function theorem can be applied and yields a deformation of $\handleapprox$ to a constant mean curvature hypersurface that is embedded if and only if $\handleapprox$ is embedded. \hfill \qedsymbol

\section{Doubling Clifford Tori --- Revisited}

\subsection{Introduction}

The methods developed in this paper allow one final example of a constant mean curvature hypersurface to be constructed.  It is the analogue of the Butscher-Pacard doubling construction for Clifford tori \cite{mepacard1,mepacard2}, except for Clifford tori of arbitrary mean curvature.  That is, one can take two Clifford tori with equal but opposite mean curvature of \emph{any} magnitude, choose a collection of symmetrically located gluing points on each, and place a sequence of hyperspheres of the same mean curvature end-to-end along the geodesics connecting the gluing points on one Clifford torus to those on the other.  The entire configuration can then be glued together by inserting catenoidal necks between the hyperspheres, and then perturbed to have exactly constant mean curvature.

\subsection{The Approximate Doubling Construction}

Let $\T^{n_1, n_2}_\alpha$ be the generalized Clifford torus of type $\Sph^{n_1} \times \Sph^{n_2}$ with  $n_1 + n_2 = n$ and parameter $\alpha \in (0, \pi/2)$ in $\Sph^{n+1}$.  This hypersurface has mean curvature $H_\alpha^{\mathit{Cliff}} := n_2 \cot(\alpha) - n_1 \tan(\alpha)$.  Let $\bar \alpha \in (0, \pi/2)$ be such that $n_2 \cot( \bar \alpha) - n_1 \tan(\bar \alpha) = - H_\alpha^{\mathit{Cliff}}$.  Note that $\alpha < \alpha_\ast < \bar \alpha$ where $\alpha_\ast = \arctan \bigl( \sqrt{n_2/n_1} \, \bigr)$ is such that $H_{\alpha_\ast}^{\mathit{Cliff}} =0$.  Then the Clifford torus $\T^{n_1, n_2}_{\bar \alpha}$ has equal but opposite mean curvature to $\T^{n_1, n_2}_\alpha$.  

Take a finite subgroup $G \subseteq SO(n+2)$, having order $|G|$, of the type considered in \cite{mepacard1,mepacard2}.  That is, suppose $G$ is a subgroup of the diagonal $O(n_1+1) \times O(n_2+1)$ in $O(n+2)$ and contains the element $T := (T_{1}, T_{2})$ where $T_j \in O(n_j+1)$ is the reflection symmetry across the $x_1$-axis defined by $T_j (x_1, x_2,\ldots, x_{n_j+1}) := (x_1, - x_2, \ldots, - x_{n_j+1})$.  Note that such a group preserves both $\T^{n_1, n_2}_\alpha$ and $\T^{n_1, n_2}_{\bar \alpha}$. Enumerate the elements of the subgroup $G$ by $G := \{ \omega_1 = \mathit{Id},\, \omega_2, \ldots, \omega_{|G|} \}$ and denote $\omega_k := (\omega_k^1, \omega_k^2)$ if needed, where $\omega_k^j \in O(n_j+1)$.   Now choose the canonical point $p \in \cliff$ defined by $p := (\cos(\alpha) P_1, \sin(\alpha) P_2)$ where $P_j \in \R^{n_j+1}$ is the point $P_j := (1, 0, \ldots, 0)$ and set $p_s := \omega_s(P)$ for every $s = 1, \ldots, |G|$.  Similarly, let $\bar p :=  (\cos(\bar \alpha) P_1, \sin(\bar \alpha) P_2) \in \T^{n_1,n_2}_{\bar \alpha}$ and set $\bar p_s := \omega_s (\bar P)$ for every $s = 1, \ldots, |G|$.  Assume that $\omega_1 = \mathit{Id}$ so that $p_1 = p$ and $\bar p_1 = \bar p$.  Let $M$ be the orbit of $p$ under $G$ and  let $\bar M$ be the orbit of $\bar p$ under $G$.   Finally, let $\gamma_s$ be the geodesic taking $p_s$ to $\bar p_s$.  Note that $\gamma_s$ is normal to $\cliff$ at $p_s$, normal to $\T^{n_1,n_2}_{\bar \alpha}$ at $\bar p_s$ and can be parametrized by  $t \longmapsto \omega_s \circ R_t (p)$ where $R_t \in SO(n+2)$ is the rotation introduced in \eqref{eqn:torusgeo}.

The next task is to position hyperspheres between $\cliff$ and $\T^{n_1, n_2}_{\bar \alpha}$ along each $\gamma_s$.  Recall that the mean curvature of a hypersphere $S_{\hat \alpha}$ with parameter $\hat \alpha \in (0, \pi/2)$ is $H^{\mathit{Sph}}_{\hat \alpha} := n \cot (\hat \alpha)$.  Therefore the  hypersphere with the same mean curvature (up to sign) as $\cliff$ is the one with $\hat \alpha := \arccot \big| \frac{n_2}{n} \cot (\alpha) - \frac{n_1}{n} \tan (\alpha) \big|$.  Suppose now that the parameters $\alpha, \bar \alpha$ and $\hat \alpha$ are such that there are positive integers $N$ and $m$ and a small separation parameter $0 < \tau \ll 2 \pi$ so that $2 N \hat \alpha + (N+1) \tau = \bar \alpha - \alpha + 2 m \pi$.  It is thus possible to position $N$ rotated  hyperspheres of parameter $\hat \alpha$ along each $\gamma_s$ from $p_s$ to $\bar p_s$, where each hypersphere is separated from its two nearest neighbours by a distance $\tau$.  These spheres wind $m$ times around $\gamma_s$ before ending up a distance of exactly $\tau$ from $\bar p_s$.

\begin{rmk}
It is clear that for every fixed $\alpha \in (0, \pi/2)$ and for every $\tau_0>0$ there is always \emph{some} value of $\tau$ satisfying $0 < \tau < \tau_0$ and \emph{some} integers $N$ and $m$ so that the equation $2 N \hat \alpha + (N+1) \tau = \bar \alpha - \alpha + 2 m \pi$ holds.     It is not clear if this equation can be made to hold in this way for fixed $m$ as well (e.g.~for $m=0$ which would correspond to embedded configurations).   But if $\alpha$ is sufficiently close to $0$ and thus $\bar \alpha$ is sufficiently close to $\pi/2$ then as before there is a sequence of examples having increasingly smaller $\tau$ satisfying the condition with $m=0$.
\end{rmk}

These considerations lead to the definition of the initial configuration of  hyperspheres and Clifford tori that will be glued together below.  Recall the rotation $Q_\alpha \in SO(n+2)$ introduced in \eqref{eqn:torusrotn} which takes the point $p_0 := (1, 0, \ldots, 0) \in \R^{n+2}$ to the point $p$, the vector $N_0 := (0, 1, 0, \ldots 0) \in \R^{n+2}$ to the vector $N_p$, and the subspace $\{ 0, 0 \} \times \R^n$ to the tangent space of $\cliff$ at $p$.  Let $\sigma_1, \ldots, \sigma_N$ be small displacement parameters and define the initial configuration
\begin{equation}
	\label{eqn:initdouble}
	\Lambda_{\alpha, \tau, \sigma}^{\#} := \cliff \cup \left[ \bigcup_{s=1}^{|G|} \bigcup_{k=1}^N \omega_s \circ R_{(2k-1) \hat \alpha + k \tau + \sigma_k} \circ Q_\alpha (S_{\hat \alpha})  \right] \cup \T^{n_1, n_2}_{\bar \alpha} \, .
\end{equation}
Note that there is redundancy in the labeling of the spheres in \eqref{eqn:initdouble} corresponding to the subgroup of $G$ fixing the point $p$.  When $\omega_s$ fixes $p$ then $\omega_s = S^{00}_{B_1, B_2}$ for some $B_j \in O(n_j+1)$ as defined in \eqref{eqn:fixsubgp}.  Thus $\omega_s \circ  R_{(2k-1) \hat \alpha + k \tau + \sigma_k} \circ Q_\alpha (S_{\hat \alpha}) = R_{(2k-1) \hat \alpha + k \tau + \sigma_k} \circ Q_\alpha \circ R(S_{\hat \alpha})$ where $R$ is a transformation preserving both the sphere $S_{\hat \alpha}$ and the $\{x^0, x^1\}$-equator.  Consequently, the entire initial configuration $\Lambda_{\alpha, \tau, \sigma}^{\#}$ is also invariant under $G$.

To glue the initial configuration together, one proceeds as follows. By analogy with the procedure described in Section \ref{subsec:handleassembly}, one first replaces each hypersphere $\omega_s \circ R_{(2k-1) \alpha + k \tau + \sigma_k} \circ Q_\alpha (S_{\hat \alpha})$ with $\omega_s \circ R_{(2k-1) \alpha + k \tau + \sigma_k} \circ Q_\alpha (\tilde S_{\hat \alpha, \eps_k, b_k})$ where $\tilde S_{\hat \alpha, \eps_k, b_k}$ is the perturbed hypersphere of Section \ref{subsec:newassembly} and the parameters $\eps_k, b_k$ are the usual scale and translation parameters.  Then on replaces the Clifford tori $\cliff$ and $\T^{n_1, n_2}_{\bar \alpha}$ with the perturbed hypersurfaces used in Butscher and Pacard's papers \cite{mepacard1,mepacard2}.  Denote these $\tilde C_{\alpha, \eps}^{n_1, n_2}$ and $\tilde C_{\bar \alpha, \eps}^{n_1, n_2}$  corresponding to  $\cliff$ and $\T^{n_1, n_2}_{\bar \alpha}$, respectively.  Then as before, one can find the parameters $\eps_k, b_k$, and $\eps$ along with neck scale and translation parameters $\bar \eps_k$ and $\bar b_k$ that ensure optimal matching between these hypersurfaces and small catenoidal necks with the parameters $\bar \eps_k, \bar b_k$.   Moreover, these are uniquely determined in terms of $\tau$.

\begin{defn}
	\mylabel{defn:doubleregions}
	Recall the terminology from Sections \ref{subsec:assembly}, \ref{subsec:newassembly} and \ref{subsec:handleassembly}. 
	\begin{itemize}
		\item Define the neck regions
		\begin{align*}
			\mathcal N^{sk, \pm}_{\eps} &:= \omega_s \circ R_{2 k \hat \alpha + (k+1/2)\tau +(\sigma_k+ \sigma_{k+1})/2} \circ Q_{\alpha} \circ K^{-1}  \big( \tilde \Sigma_{\bar \eps_k, \bar b_k}^\pm \cap \mathit{Cyl}(\rho_{\eps}/2) \big) \, .
		\end{align*}
	Set $\mathcal N_{\eps}^{sk} := \mathcal N_{\eps}^{sk,+} \cup \mathcal N_{\eps}^{sk,-} $.  These necks are oriented so that $\mathcal N_{\eps}^{sk,+}$ lies ahead of $\mathcal N_{\eps}^{sk,-}$ along the geodesic $\gamma_s$ and that the $k^{\mathit{th}}$ neck is centered on the point $\omega_s \circ R_{2 k \hat \alpha + (k+1/2)\tau +(\sigma_k+ \sigma_{k+1})/2}(p)$.  
	
		\item Define the transition regions
		\begin{align*}
			\mathcal T^{sk, \pm}_{\eps} &:= \omega_s \circ R_{2 k \hat \alpha + (k+1/2)\tau +(\sigma_k+ \sigma_{k+1})/2} \circ Q_\alpha \circ K^{-1}  \big( \tilde \Sigma_{\bar \eps_k, \bar b_k}^\pm \cap \big[ \mathit{Cyl}(2\rho_\eps) \setminus \mathit{Cyl}(\rho_\eps/2) \big] \big) \, .
		\end{align*}
	Set $\mathcal T_{\eps}^{sk} := \mathcal T_{\eps}^{sk,+} \cup \mathcal T_{\eps}^{sk,-} $.  	
	
		\item Define the exterior regions
		\begin{align*}
			\mathcal E_{\alpha, \eps}^{\mathit{Cliff}} &:= \tilde C_{\alpha, {\eps_0}}^{n_1, n_2} \setminus \big[ B_{\tilde \rho_{\eps_0}} (p) \cup B_{\tilde \rho_{\eps_0}} (-p') \big] \\[0.5ex]
			\mathcal E_{\bar \alpha, \eps}^{\mathit{Cliff}} &:= \tilde C_{\bar \alpha, {\eps_0}}^{n_1, n_2} \setminus \big[ B_{\tilde \rho_{\eps_0}} (p) \cup B_{\tilde \rho_{\eps_0}} (-p') \big] \\[0.5ex]
			\mathcal E_{\hat \alpha, \eps}^{sk} &:= \omega_s \circ R_{(2k-1) \hat \alpha + k \tau + \sigma_k} \circ Q_\alpha  \big( \tilde S_{\alpha, \eps_k, b_k} \setminus \big[ B_{\tilde \rho_{\eps}} (p^+) \cup B_{\tilde \rho_{\eps}} (p^-)\big] \big) \, .
		\end{align*} 
	\end{itemize}
\end{defn}

\noindent The approximate solution corresponding to the configuration $\Lambda _{\alpha,\tau,\sigma}^{\#}$ can now be defined as follows.
\begin{defn}
	\mylabel{defn:doubleapproxsol}
	The \emph{approximate solution} with parameters $\alpha, \tau$ and $\sigma$ is the hypersurface 
	\begin{equation}
		\mylabel{eqn:doubleapproxsol}
		\doubleapprox := \mathcal E_{\alpha, \eps}^{\mathit{Cliff}} \cup \big[\mathcal T_{\eps}^{ s0} \cup \mathcal N_{\eps}^{s0} \big] \cup \left[\bigcup_{s=1}^{|G|} \bigcup_{k =1}^{ N} \Big(  \mathcal E_{\hat \alpha, \eps}^{sk} \cup  \mathcal T_{\eps}^{ sk} \cup \mathcal N_{\eps}^{sk}\Big) \right] \cup \mathcal E_{\bar \alpha, \eps}^{\mathit{Cliff}}  \, .
	\end{equation}	
\end{defn}

\subsection{Symmetries and Jacobi Fields}

The initial configuration $\Lambda_{\alpha, \tau, \sigma}^{\#}$ is $G$-invariant and it is clear that the gluing procedure leading up to the definition of the approximate solution $\doubleapprox$ can also be done in a $G$-invariant manner.  It remains to see which of the Jacobi fields of $\doubleapprox$ are $G$-invariant. As before, one proceeds as follows.  First define a smooth and monotone cut-off function $\chi_{sk} : \doubleapprox \rightarrow \R$ supported in $\mathcal E^{sk}_{\eps}$ and two other such functions $\chi_{\alpha, \mathit{Cliff}} : \handleapprox \rightarrow \R$ supported in $\mathcal E^{\mathit{Cliff}}_{\alpha, \eps}$ and  $\chi_{\bar \alpha, \mathit{Cliff}} : \handleapprox \rightarrow \R$ supported in $\mathcal E^{\mathit{Cliff}}_{\bar \alpha, \eps}$.  Without loss of generality, the functions $\chi_{\alpha, \mathit{Cliff}}$ and $\chi_{\bar \alpha, \mathit{Cliff}}$ are $G$-invariant while $\chi_{sk} \circ \omega_s \circ \omega_t^{-1} = \chi_{tk}$ holds for all $s,t = 1, \ldots, |G|$.   Consider $\tilde S_{\hat \alpha, \eps}$ as a hypersurface in $\Sph^{n+1} \subseteq \R^{n+2}$ endowed with coordinates $(x^0, x^1 , \ldots, x^{n+1})$ and define the following functions.
\begin{subequations}
	\mylabel{eqn:doubleapproxjac}
	\begin{equation}
		\tilde q_{sk}^{\, t}  := 
			\chi_{sk} \cdot \left( x^t  \big\vert_{\tilde S_{\hat \alpha, \eps} } \circ Q_\alpha^{-1} \circ R^{-1}_{(2k-1) \alpha + k \tau+ \sigma_k} \circ \omega_s^{-1} \right) \, .
	\end{equation}
		Also, in the standard coordinates $(x_1, x_2)$ for $\R^{n_1+1} \times \R^{n_2+1}$ that are in use at the moment, define the functions
	\begin{equation}
		\begin{aligned}
			\tilde q^{\, t_1, t_2}_{\alpha, \mathit{Cliff}}  &:= \chi_{\alpha, \mathit{Cliff}} \cdot \left( x_1^{t_1} x_2^{t_2}  \big\vert_{\tilde C^{n_1, n_2}_{\alpha, \eps}} \right) \\
			\tilde q^{\, t_1, t_2}_{\bar\alpha, \mathit{Cliff}}  &:= \chi_{\bar\alpha, \mathit{Cliff}} \cdot \left( x_1^{t_1} x_2^{t_2}  \big\vert_{\tilde C^{n_1, n_2}_{\bar \alpha, \eps}} \right)  \, .
		\end{aligned}
	\end{equation}
\end{subequations}
Then the vector space of \emph{all} approximate Jacobi fields in $C^{2, \beta}(\handleapprox)$ is the set of functions
$$\tilde{ \mathcal K} := \linspan_\R  \bigcup_{s=1}^{|G|} \bigcup_{k=1}^N \{ \tilde q_{sk}^{\, t} : t = 1, \ldots, n \} \cup \{ \tilde q^{\, t_1, t_2}_{\alpha , \mathit{Cliff}},  \tilde q^{\, t_1, t_2}_{\bar \alpha , \mathit{Cliff}} :  t_1 = 1, \ldots, n_1 \mbox{ and } t_2 = 1, \ldots, n_2 \} \, .$$
At this point, impose the condition on the group $G$ used in Butscher and Pacard's doubling construction \cite{mepacard1,mepacard2} that there is no linear combination of the $\tilde q^{\, t_1, t_2}_{\alpha , \mathit{Cliff}}$ and $\tilde q^{\, t_1, t_2}_{\bar \alpha , \mathit{Cliff}}$ that is invariant under every element of $G$.  This eliminates the Jacobi fields on the perturbed Clifford tori from consideration when only equivariant perturbations are allowed.  Furthermore, invariance with respect to $G$ implies first of all that the Jacobi fields on the $k^{\mathit{th}}$ perturbed hypersphere along every geodesic are the same. Then invariance with respect to the reflection $T \in G$ implies that the only invariant Jacobi fields of $\doubleapprox$ are to be found in the following subspace.

\begin{defn}
	The subspace of approximate Jacobi fields of  $\handleapprox$ invariant under $G$ is the subspace of functions 
	\begin{equation}
	\label{eqn:doubleapker}
		\tilde{\mathcal K}_{\mathit{sym}} := \linspan_\R \left\{ \sum_{s=1}^{|G|} \tilde q^1_{sk} : k = 1, \ldots, N \right\} \, .
	\end{equation}
	Note that the dimension of $\mathcal K_{\mathit{sym}}$ is $N$.
\end{defn}

\subsection{The Proof of Main Theorem 4}

The proof of Main Theorem 4 is now in all respects identical to the proof of Main Theorem 3, for the following reasons.  One can set up normal deformations of $\doubleapprox$ using functions invariant under $G$.  Then the construction of the right inverse of the linearized mean curvature operator and its estimate follow as in Proposition \ref{thm:seclinest} and the analogous constructions performed in Butscher and Pacard's papers \cite{mepacard1,mepacard2} because the nature of the Jacobi fields of $\doubleapprox$ is unchanged; namely, there is one invariant Jacobi field on each perturbed hypersphere and none on the necks and on the perturbed Clifford hyperspheres by the nature of the symmetries and the choice of weighted function space.  This gives the linear estimate \eqref{eqn:iftestone} needed to invoke the inverse function theorem.  Then the estimates of the mean curvature of the Delaunay-like hypersurfaces and the estimates of the mean curvature of the generalized Clifford tori from \cite{mepacard1,mepacard2} can again be combined to yield the first of the non-linear estimates \eqref{eqn:iftesttwo} needed to invoke the inverse function theorem.  Finally, the second of the non-linear estimates holds as before.  Hence when $\tau$ is sufficiently small the inverse function theorem can be applied and yields a deformation of $\doubleapprox$ to a constant mean curvature hypersurface that is embedded if and only if $\doubleapprox$ is embedded. \hfill \qedsymbol

\bibliography{sphcmc}
\bibliographystyle{amsplain}

\end{document}